%Paper: alg-geom/9308004
%From: rf@shire.math.columbia.edu (Robert Friedman)
%Date: Mon, 23 Aug 93 14:58:28 -0400

%%Paper written in AMS-TeX
%%using the amsppt style

\input amstex

\define\scrO{\Cal O}
\define\Pee{{\Bbb P}}
\define\Zee{{\Bbb Z}}
\define\Cee{{\Bbb C}}
\define\Ar{{\Bbb R}}

\define\im{\operatorname{Im}}
\define\Hom{\operatorname{Hom}}
\define\Sym{\operatorname{Sym}}

\define\Id{\operatorname{Id}}
\define\Pic{\operatorname{Pic}}
\define\Supp{\operatorname{Supp}}
\define\ad{\operatorname{ad}}
\define\ch{\operatorname{ch}}
\define\Todd{\operatorname{Todd}}
\define\Ext{\operatorname{Ext}}
\define\Spec{\operatorname{Spec}}

\define\pt{\text{pt}}
\define\proof{\demo{Proof}}
\define\endproof{\qed\enddemo}
\define\endstatement{\endproclaim}
\define\theorem#1{\proclaim{Theorem #1}}
\define\lemma#1{\proclaim{Lemma #1}}
\define\proposition#1{\proclaim{Proposition #1}}
\define\corollary#1{\proclaim{Corollary #1}}
\define\claim#1{\proclaim{Claim #1}}

\define\section#1{\specialhead #1 \endspecialhead}
\define\ssection#1{\medskip\noindent{\bf #1}}

\documentstyle{amsppt}
\leftheadtext{Vector bundles  and $SO(3)$-invariants for elliptic surfaces
III}
\rightheadtext{The case of odd fiber degree}
\pageno=1
\topmatter
\title Vector bundles  and $SO(3)$-invariants for elliptic surfaces III:\\
The case of odd fiber degree
\endtitle
\author Robert Friedman
\endauthor
\address Columbia University, New York, NY 10027
\endaddress
\email rf\@math.columbia.edu \endemail
\thanks Research partially supported by NSF grant  DMS-9203940
\endthanks
\endtopmatter

\section{Introduction.}

Let $S$ be a simply connected elliptic surface with at most two multiple
fibers,  of multiplicities $m_1$ and $m_2$, where one or both of the $m_i$
are allowed to be 1. In this paper, the last of a series of three, we shall
study stable rank two vector bundles $V$ on $S$ such that $\det V \cdot f$ is
odd, where $f$ is a general fiber of $S$. Thus necessarily the multiplicities
$m_1$ and $m_2$ are odd as well. Bundles $V$ such  that $\det (V|f)$ has even
degree for a general fiber
$f$ have been studied extensively [3], [4, Part II], and as we shall see the
case of odd fiber degree is fundamentally different. Thus we shall have to
develop the analysis of the relevant vector bundles from scratch, and the
results in this paper are for the most part independent of those in [3] and
[4].  Our goal in this paper is to give a description of the moduli space of
stable rank two bundles with odd fiber degree and then to use this information
to calculate certain Donaldson polynomials. Before stating our main result,
recall that, for an elliptic surface
$S$, $J^d(S)$ denotes the elliptic surface whose general fiber is the set of
line bundles of degree $d$ on the general fiber of $S$.  We shall prove the
following two theorems:

\theorem{1} Let $\frak M_t$ be the moduli space of stable rank two bundles
$V$  on
$S$ \rom(with respect to a suitable ample line bundle\rom) with $\det V
\cdot  f = 2e+1$ and $4c_2(V) - c_1^2(V) - 3\chi (\scrO _S) = 2t$. Then
$\frak M_t$ is  smooth and irreducible and is birational to $\Sym
^tJ^{e+1}(S)$.
\endstatement
\medskip

\theorem{2} Let $\gamma _t$ be the Donaldson polynomial of degree $2t$
corresponding to the choice of moduli space $\frak M_t$. Let $\kappa \in
H_2(S;\Zee)$ be the primitive element such that $m_1m_1\kappa = f$. Then, for
all
$\Sigma \in H_2(S)$,   \roster
\item"{(i)}"  $\gamma _0 = 1$.
\item"{(ii)}" $\gamma _1(\Sigma, \Sigma) = (\Sigma ^2) + ((m_1^2m_2^2)(p_g(S)
+ 1)- m_1^2-m_2^2) (\Sigma \cdot \kappa)^2$.
\item"{(iii)}"  $\gamma _2(S)(\Sigma, \Sigma ,\Sigma, \Sigma ) = 3(\Sigma
^2)^2 + 6C_1(\Sigma ^2 )(\Sigma \cdot \kappa)^2 + (3C_1^2 - 2C_2)(\Sigma
\cdot
\kappa)^4$, where
$$\align C_1&= (m_1^2m_2^2)(p_g(S) + 1)-m_1^2-m_2^2;\\ C_2 &= (m_1^4m_2^4)
(p_g(S) + 1) -m_1^4-m_2^4.
\endalign$$
\endroster
\endstatement
\medskip

Let us outline the basic ideas behind the proof of Theorem 1. Standard
arguments show that, for a suitable choice of an ample line bundle $L$ on
$S$, a rank two vector bundle $V$ with $c_1(V)\cdot f= 2e+1$ is $L$-stable if
and only if its restriction to a general fiber $f$ is stable. A pleasant
consequence of the assumption of odd fiber degree is that there is a unique
stable bundle of a given determinant of odd degree on each smooth fiber $f$.
Using this, it is easy to show that there exists a rank two vector bundle
$V_0$ whose restriction to {\it every\/} fiber $f$ is stable, and that $V_0$
is unique up to twisting by a line bundle. The bundle $V_0$ is the progenitor
of all stable bundles on $S$, in the sense that every stable rank two vector
bundle is obtained from $V_0$ by making elementary modifications along fibers.
Generically, this involves choosing $t$ smooth fibers
$f_i$ and line bundles $\lambda _i$ of degree $e+1$ on $f_i$. These choices
define the birational isomorphism from the moduli space to  $\Sym
^tJ^{e+1}(S)$.

Given the above analysis of stable bundles, the main problem in computing
Donaldson polynomials is to fit together all of the various possible
descriptions of stable bundles into a universal family whose Chern classes
can be calculated. This is easier said than done! Even in the case where $S$
has a section, the construction of the universal bundle for the
four-dimensional moduli space, which just involves well-known techniques of
extensions and elementary modifications, is already quite involved. We shall
therefore proceed differently, and try to describe the moduli spaces and
Chern classes involved up to contributions which only depend on the analytic
type of a neighborhood of the multiple fibers. But we shall not try to
analyze these contributions explicitly. Instead we shall repeatedly use the
fact that an elliptic surface with $p_g=0$ and just one multiple fiber is a
rational surface, and thus its Donaldson polynomials are the same as those
for an elliptic surface with $p_g=0$ and with a section, or equivalently no
multiple fibers. Thus if we know these, we can try to interpolate this
knowledge into the general case. We shall use this idea twice. The first
application will be to calculate the invariant $\gamma _1$. Here the moduli
space is $J^{e+1}(S)$ and a lengthy calculation with the
Grothendieck-Riemann-Roch theorem identifies the divisor corresponding to the
$\mu$-map up to a rational multiple of the fiber, which depends only on the
multiplicities. Appealing to the knowledge of the invariant for a rational
surface enables us to determine this multiple. Of course, it is likely that
the exact multiple could also be determined by a direct calculation. In order
to calculate the polynomial
$\gamma _2$, we shall use a variant of this idea. In this case, the divisor
corresponding to the
$\mu$-map is essentially known from the corresponding calculation in the case
of
$\gamma _1$. However what changes is the moduli space itself: the presence of
multiple fibers means that the birational map from the moduli space to
$\operatorname{Hilb}^2J^{e+1}(S)$ is not a morphism, and the actual moduli
space differs from  $\operatorname{Hilb}^2J^{e+1}(S)$ in codimension two.
Thus while the divisors are known, their top self-intersection is not. Again
using the rational elliptic surfaces, we are able to determine the
discrepancy between the self-intersection of the $\mu$-divisors in
$\operatorname{Hilb}^2J^{e+1}(S)$ and in the actual moduli space. The methods
used here are in a certain sense the analogue in algebraic geometry of gluing
techniques for ASD connections.

Although the actual arguments are rather involved, the main point to
emphasize here is that the coefficients of the Donaldson polynomial are quite
formally determined by the knowledge of the polynomial for a rational
surface. It is natural to wonder if the techniques in this paper can be
pushed further to determine $\gamma _t$ for all $t$. I believe that this
should be possible, although one necessary and so far missing ingredient in
this approach is the knowledge of the multiplication table for divisors in
$\operatorname{Hilb}^2 J^{e+1}(S)$.

Here is a rapid description of the contents of the paper. In Section 1 we
describe some general results on rank two vector bundles on an elliptic
curve. In Section 2 these results are extended to cover the case of an
irreducible nodal curve of arithmetic genus one. In Section 3 we give the
classification of stable bundles on an elliptic surface $S$ and prove Theorem
1. In Section 4 we specialize to the case of a surface with a section. Our
purpose here is twofold: First, we would like to show how many of the results
of the preceding section take a very concrete form in this case. Secondly, we
shall make a model for a piece of the four-dimensional moduli space which we
shall need to use later. In Section 5 we calculate the two-dimensional
invariant $\gamma _1$ in case $S$ has a section. This calculation has already
been done by a different method in Section 4 and will be redone in full
generality. However it seemed worthwhile to do this special case in order to
make the general calculation more transparent. The next three sections are
devoted to calculating $\gamma _1$ in general. The outline of the argument is
given in Section 6. We construct a coherent sheaf which is an approximation
to the universal bundle over the moduli space, over a branched cover $T$ of
$S$. We determine its Chern classes via a lengthy calculation using the
Grothendieck-Riemann-Roch theorem, which is given in Section 7. The necessary
correction terms are identified via the results in Section 8. In Sections 9
and 10 we deal with the invariant $\gamma _2$. Once again the outline of the
calculation is given first and the technical details are postponed to Section
10. The paper concludes with an appendix which collects some general results
about elementary modifications.

\section{Notation, conventions, and preliminaries.}

All spaces are over $\Cee$, all sheaves are coherent sheaves in the classical
topology unless otherwise specified. We do not distinguish between a vector
bundle and its locally free sheaf of sections. Given s subvariety $Y$ of a
compact complex manifold $X$, we denote the associated cohomology class by
$[Y]$.

If $V$ is a rank two vector bundle on a complex manifold or smooth scheme
$X$, we shall frequently need to  consider the first Pontrjagin class of
$\operatorname{ad}V$, which is $c_1^2(V) - 4c_2(V)$. We will  denote this
expression by $p_1(\operatorname{ad}V)$. We shall occasionally and
incorrectly use the shorthand  $p_1(\operatorname{ad}V)$, for an arbitrary
coherent sheaf
$V$, to denote $c_1^2(V) - 4c_2(V)$.

Given a vector bundle $V$, we shall need to know how
$p_1(\operatorname{ad}V)$  changes under elementary modifications. Recall
that an  elementary modification is defined as follows.  Let $X$ be a smooth
scheme and let $D$ be an effective divisor on $X$, not  necessarily smooth,
with $i\: D\to X$ the inclusion. Let $L$ be a line bundle on $D$. Then $i_*L$
is a coherent sheaf on
$X$, which we shall frequently just denote by $L$. Suppose that $V_0$ is a
rank two vector bundle on $X$ and that $V_0 \to i_*L$ is a surjective
homomorphism. Let $V$ be the kernel of the map $V \to i_*L$. Then $V$ is again
a rank two vector bundle on $X$ (and in particular it is locally free). We
call $V$ an {\sl elementary modification\/} of
$V_0$. The change in $p_1$ is given as follows:

\lemma{0.1} Let $X$ be a smooth scheme and let $D$ be an effective divisor on
$X$,  not necessarily smooth. Let $L$ be a line bundle on $D$ and $V_0$ a rank
two vector bundle, and suppose that there is an exact sequence
$$0 \to V \to V_0 \to i_*L \to 0,$$ where $i\: D \to X$ is the inclusion. Then
$$p_1(\operatorname{ad}V) - p_1(\operatorname{ad}V_0) = 2c_1(V_0)\cdot [D] +
[D]^2 -4 i_*c_1(L).$$
\endstatement
\proof The proof follows easily from standard formulas for $c_1(V)$ and
$c_2(V)$,  cf\. [7] or [5].
\endproof

Next we will recall some properties of the scheme $\operatorname{Hilb}^2S$,
where
$S$ is an algebraic surface. In general, we denote by
$\operatorname{Hilb}^nS$ the smooth projective scheme parametrizing
$0$-dimensional subschemes of $S$ of length $n$. There is a universal
codimension two subscheme $\Cal Z
\subset S \times \operatorname{Hilb}^nS$. We may describe the case $n=2$ quite
explicitly. Let $\tilde H$ be the blowup of $S\times S$ along the diagonal
$\Bbb D$ and let $\tilde \Bbb D$ be the exceptional divisor. There is an
involution $\iota$ of $\tilde H$ whose fixed set is $\tilde \Bbb D$. We claim
that the quotient
$\tilde H /\iota$ is naturally $\operatorname{Hilb}^2S$. Indeed, if
$\tilde\Bbb D_{12}$ and $\tilde\Bbb D_{13}$ are the proper transforms in
$S\times \tilde H$ of the subsets $$\Bbb D_{1j} = \{\, p \in S\times S\times
S \mid \pi _1(p) = \pi _j(p)\,\},$$ then $\tilde \Cal Z = \tilde \Bbb D_{12}
+ \tilde \Bbb D_{13}$ is a codimension two subscheme of $S\times \tilde H$
which is easily seen to be a local complete intersection. Thus it defines a
flat family of subschemes of $S$ and so a morphism $\pi \:\tilde H \to
\operatorname{Hilb}^2S$. It is easy to see that the induced morphism $\tilde
H/\iota \to  \operatorname{Hilb}^2S$ is an isomorphism. The projection
$\Cal Z
\to  \operatorname{Hilb}^2S$ is a double cover which identifies $\Cal Z$ with
$\tilde H$.

Given $\alpha \in H_2(S)$, we can define the element $D_\alpha \in
H^2(\operatorname{Hilb}^2S)$ by taking slant product with $[\Cal Z] \in
H^4(S\times \operatorname{Hilb}^2S)$. If for example $\alpha = [C]$ where $C$
is an irreducible curve on $S$, then $D_\alpha$ is represented by the
effective divisor consisting of pairs $\{x, y\}$ of points of $S$ such that
either $x$ or
$y$ lies on $C$. The inverse image $\pi ^*D_\alpha\in H^2(\tilde H)$ is the
pullback of the class $1\otimes \alpha + \alpha \otimes 1\in H^2(S\times S)$.
There is also the class in $ H^2(\operatorname{Hilb}^2S)$ represented by the
divisor $E$ of subschemes of $S$ whose support is a single point. Since $\pi$
is branched over $E$, the class $[E]$ is divisible by $2$ and $\pi ^*[E]=
2\tilde \Bbb D$. Using this it is easy to check that the map $\alpha \mapsto
D_\alpha$ defines an injection $H_2(S) \to H^2(\operatorname{Hilb}^2S)$ and
that
$H^2(\operatorname{Hilb}^2S) = H_2(S) \oplus \Zee\cdot [E/2]$. Finally the
multiplication table in $H^2(\operatorname{Hilb}^2S)$ can be determined from
the fact that $\tilde H$ is the blowup of $S\times S$ along the diagonal and
that the normal bundle of the diagonal in $S\times S$ is the tangent bundle of
$S$: we have
$$\gather D_\alpha ^4 = 3(\alpha ^2)^2; \qquad D_\alpha ^3\cdot E = 0;
\qquad D_\alpha ^2\cdot E^2 = -8(\alpha ^2); \\
 D_\alpha \cdot E^3 = -8(c_1(S)\cdot \alpha); \qquad E^4 = 8(c_2(S) -
c_1(S)^2).
\endgather$$

Finally we need to say a few words about calculating Donaldson polynomials.
Let
$M$ be a closed oriented simply connected 4-manifold with a generic
Riemannian metric $g$, and let $P$ be a principal $SO(3)$-bundle over $M$
with invariants
$w_2(P)=w$ and $p_1(P) = p$. There is a Donaldson polynomial  $\gamma
_{w,p}(S)$ defiend via the moduli space of $g$-ASD connections on $P$,
together with a choice of orientation for this space. If $b_2^+(M)>1$, then
this polynomial is independent of $g$, whereas if $b_2^+(M)=1$ then it only
depends on a certain chamber in the positive cone of $H^2(M; \Ar)$. If $M=S$
is a complex surface,
$\Delta$ is a holomorphic line bundle such that $w= c_1(\Delta )\mod 2$ and
$g$ is a Hodge metric corresponding to an ample line bundle $L$, there is a
diffeomorphism of real analytic spaces from the moduli space of $g$-ASD
connections on $P$ to the moduli space of $L$-stable rank two vector bundles
$V$ on $S$ with $c_1(V) =
\Delta$ and $c_2(V) = (\Delta ^2-p)/4$. We denote this moduli space for the
moment by $\frak M$. We shall always choose the orientation of the moduli
space of
$g$-ASD connections which agrees with the natural complex orientation of
$\frak M$.

If $\frak M$ is smooth, compact, and of real dimension $2d$ and there is a
universal bundle $\Cal V$ over $S\times \frak M$, then slant product with
$-p_1(
\ad \Cal V)/4$ defines a homomorphism $\mu$ from $H_2(S)$ to $H^2(\frak M)$.
In general we can define the holomorphic vector bundle $\ad \Cal V$ even when
the universal bundle $\Cal V$ does not exist. To see this, note that there is
always a universal $\Pee ^1$-bundle $\pi \: \Pee (\Cal V)\to S\times \frak M$,
and taking $\pi _*$ of the relative tangent bundle gives $\ad \Cal V$. Thus
given a class $\Sigma \in  H_2(S)$, we can evaluate $\mu (\Sigma)^d$ on the
fundamental class of $\frak M$ and this gives the value $\gamma _{w,p}(\Sigma,
\dots,
\Sigma)$. For the applications in this paper, since the moduli spaces always
have the correct dimension and in particular are empty if
$-p-3\chi(\scrO_S)<0$, the moduli spaces of complex dimension zero and two
are compact. For the four-dimensional moduli space, we can calculate $\gamma
_{w,p}$ by choosing an appropriate compactification of $\frak M$. For the
purposes of gauge theory, there is the Uhlenbeck compactification. For the
purposes of algebraic geometry, there is the Gieseker compactification
$\overline{\frak M}$. Following O'Grady [11], the divisors $\mu (\Sigma)$
extend naturally to divisors $\nu(\Sigma)$ on
$\overline{\frak M}$, which we shall continue to denote by $\mu(\Sigma)$. If
there is a universal sheaf
$\Cal V$ on the Gieseker compactification, then the
$\mu$-map is again defined by taking slant product with $-p_1(\ad
\Cal V)/4$. In general, for holomorphic curves $\Sigma$ (which would suffice
for the applications in this paper) we can use determinant line bundles on
the moduli functor. For a general $\Sigma \in H_2(S)$, we can define
$\mu(\Sigma)$ for the moduli spaces that arise in this paper (where there are
no strictly semistable sheaves) as follows: there exists a universal coherent
sheaf
$\Cal E$  over $S\times U$, where $U$ is the open subset of an appropriate
Quot scheme corresponding to stable torsion free sheaves with the appropriate
Chern classes. Thus we can define an element of $H^2(U)$  by taking slant
product with
$p_1(\ad \Cal E)$. As $\overline{\frak M}$ is a quotient of $U$ by a free
action of $PGL(N)$ for some $N$, $H^2(\overline{\frak M})\cong H^2(U)$, and
this defines
$\mu(\Sigma)$ in general.

 We can now evaluate $\mu (\Sigma)^d$ on the fundamental class of
$\overline{\frak M}$.   By recent results of Li [9] and Morgan [10] the
result is again $\gamma _{w,p}$. Strictly speaking, their results are stated
with certain extra assumptions. However, the cases we will need in this paper
involve the following situation: all moduli spaces are smooth of the expected
dimension and there are no strictly semistable torsion free sheaves. Under
these assumptions, the proofs in e.g\. [10] go over essentially unchanged.

\section{1. Review of results on vector bundles over elliptic curves.}

We recall the following well-known result of Atiyah [1]:

\theorem{1.1} Let $V$ be a rank two vector bundle over a smooth curve $C$ of
genus
$1$. Then exactly one of the following holds:
\roster
\item"{\rm (i)}" $V$ is a direct sum of line bundles;
\item"{\rm (ii)}" $V$ is of the form $\Cal E \otimes L$, where $L$ is a line
bundle on $C$ and $\Cal E$ is the \rom(unique\rom) extension of $\scrO _C$ by
$\scrO _C$ which does not split into the direct sum $\scrO _C\oplus \scrO _C$;
\item"{\rm (iii)}" $V$ is of the form $\Cal F_p\otimes L$, where $L$ is a line
bundle on $C$, $p \in C$,  and $\Cal F_p$ is the unique nonsplit extension of
the form $$0 \to \scrO _C \to \Cal F_p \to \scrO _C(p) \to 0.\qed $$
\endroster
\endstatement

We shall not prove (1.1) but shall instead prove the analogous statement in
the slightly more complicated case of a singular curve in Section 2.

\corollary{1.2} Let $V$ be a stable rank two bundle over a smooth curve $C$
of  genus $1$. Then $\deg V$ is odd, say $\deg V = 2e +1$. Moreover, for every
line  bundle $\lambda$ of degree $e+1$ we have $\dim \operatorname{Hom}(V,
\lambda) = 1$,
$H^1(V\spcheck\otimes \lambda ) = 0$, and there is an exact sequence $$0 \to
\mu
\to V \to \lambda \to 0,$$ where $\mu$ is a line bundle of degree $e$ on $C$,
uniquely determined by the isomorphism
$$\mu \otimes \lambda = \det V,$$ and the surjection $V \to \lambda$ is unique
mod scalars. \endstatement
\proof Clearly, if $V$ is stable we must be in case (iii) of the theorem.
Conversely, suppose that $V$ is as in (iii). We shall show that $V$ is stable.
It suffices to show that $\Cal F_p$ is stable. Let $M$ be a line bundle on $C$
of degree at least $\det \Cal F_p /2 = 1/2$ such that there is a nonzero map
$M
\to \Cal F_p$. Clearly $\deg M \leq 1$ and $\deg M = 1$ if and only if
$M = \scrO _C(p)$. Since
$\deg M \geq 1/2$, $\deg M = 1$ and $M = \scrO _C(p)$. But then $\Cal F_p$ is
the split extension, contradicting the definition of $\Cal F_p$. Thus $\Cal
F_p$ is stable.

Now let $V$ be a stable bundle of degree $2e+1$, so that there exists a line
bundle $L$ of degree $e$ on $C$ with $V = \Cal F_p \otimes L$. Then, if
$\lambda$ is a line bundle of degree $e+1$, we have an exact sequence $$0 \to
\operatorname{Hom} (L\otimes \scrO_C(p), \lambda) \to  \operatorname{Hom} (V,
\lambda) \to  \operatorname{Hom}(L, \lambda) \to H^1(\lambda \otimes
L^{-1}\otimes
\scrO_C(-p)).$$ If $\lambda = L\otimes \scrO_C(p)$, then by assumption there
exists a surjection $V\to \lambda$. If $\varphi _1$ and $\varphi _2$ are two
nonzero maps from $V$ to $\lambda$, then for every $p\in C$ there is a scalar
$c$  such that $\varphi _1 -c\varphi _2$ vanishes at $p$, and thus defines a
map $V\to
\lambda \otimes \scrO_C(-p)$. By stability this map must be zero, so that
$\varphi _1 =c\varphi _2$. Thus the surjection is unique mod scalars.

If $\lambda \neq L\otimes
\scrO_C(p)$, then $\lambda\otimes L^{-1}\otimes \scrO_C(-p)$ is a line bundle
of degree zero on $C$ which is not trivial. Hence $H^1(\lambda \otimes
L^{-1}\otimes
\scrO_C(-p))=0$, and $\operatorname{Hom}(V, \lambda) \cong
\operatorname{Hom}(L,
\lambda) $. Moreover $\operatorname{Hom}(L, \lambda) = H^0(L^{-1}\otimes
\lambda)$ has dimension one since $\deg (L^{-1}\otimes \lambda) = 1$. Thus
there is a  nontrivial map $V\to \lambda$, which is unique mod scalars. If it
is not surjective, there is a factorization $V\to \lambda \otimes \scrO_C(-q)
\subset
\lambda$, and this contradicts the stability of $V$. Lastly we see that
$H^1(V\spcheck\otimes  \lambda) \cong  H^1(L^{-1}\otimes \lambda)$, and this
last group is zero since $\deg (L^{-1}\otimes \lambda) = 1$. \endproof

We can generalize the last statement of (1.2) as follows.

\lemma{1.3} Let $C$ be a smooth curve of genus one.
\roster
\item"{(i)}" Let $V$ be a stable rank two vector bundle over $C$ and suppose
that $\deg V = 2e +1$. Let $d \geq e+1$, and let $\lambda$ be a line bundle on
$V$ of degree $d$. Then $\dim \operatorname{Hom}(V, \lambda) = 2d-2e-1$, and
there  exists a surjection from $V$ to $\lambda$. Conversely, with $V$ as
above, let $\lambda$ be a line bundle such that there  exists a nonzero map
from $V$ to $\lambda$. Then $\deg \lambda \geq e+1$.
\item"{(ii)}" Suppose that $V=L_1\oplus L_2$ is a direct sum of line bundles
$L_i$ with $\deg V = 2e+1$ and $\deg L_1 \leq e<\deg L_2$. Let $\lambda$ be a
line bundle on $C$ with $d=\deg \lambda > \deg L_2$. Then $\dim
\operatorname{Hom}(V,
\lambda) = 2d-2e-1$, and there exists a surjection from $V$ to $\lambda$.
Conversely, if $\lambda$  is a line bundle and there exists a surjection from
$L_1\oplus L_2$ to $\lambda$, then either $\deg \lambda > \deg L_2$ or
$\lambda = L_2$ or $\lambda = L_1$. If $\lambda = L_2$, then $\dim
\operatorname{Hom}(V,
\lambda) = 2d-2e$, where $d=\deg L_2= \deg \lambda$. \endroster
\endstatement
\proof We shall just prove (i), as the proof of (ii) is simpler. Let $\lambda$
be  a  line bundle on $C$ of degree $d\geq e+1$. We may assume that $\deg
\lambda >e+1$, the case $\deg \lambda = e+1$ having been dealt with in (1.2).
There is an exact sequence $$0 \to L_1 \to V \to L_2 \to 0,$$ where $\deg L_1
= e$ and $\deg L_2 = e+1$. Thus there is an exact sequence $$0 \to
H^0(L_2^{-1}\otimes \lambda)
\to \operatorname{Hom}(V, \lambda) \to  H^0(L_1^{-1}\otimes \lambda)\to
H^1(L_2^{-1}\otimes \lambda).$$ We have  $\deg (L_1^{-1}\otimes \lambda )=
d-e>0$ and $\deg (L_2^{-1}\otimes \lambda )= d-e-1>0$. Thus
$H^1(L_2^{-1}\otimes \lambda) =0$, $\dim H^0(L_1^{-1}\otimes \lambda) = d-e$,
and $\dim H^0(L_2^{-1}\otimes
\lambda) = d-e-1$. So  $\dim \operatorname{Hom}(V, \lambda) = 2d-2e-1$. To see
the last statement, let $Y$ be the set of elements $\varphi$ of
$\operatorname{Hom}(V,
\lambda)$ such that $\varphi$ is not surjective. Then $Y$ is the union over
$x\in C$ of the spaces  $\operatorname{Hom}(V, \lambda\otimes \scrO _C(-x))$,
each of which has dimension at most $2d-2e -3$. So the dimension of $Y$ is at
most $2d-2e -2$. Thus  $\operatorname{Hom}(V, \lambda)-Y$ is nonempty, and
every $\varphi \in
\operatorname{Hom}(V, \lambda)-Y$ is a surjection. The final statement of (i)
is then an immediate consequence of the stability of
$V$. \endproof

For future use let us also record the following lemmas:

\lemma{1.4} Let $C$ be a smooth curve of genus one and let $\xi$ be a  line
bundle on $C$ of degree zero such that $\xi ^{\otimes 2} \neq 0$. Let $V$ be
a stable rank two vector bundle on $C$. Then $\Hom (V, V\otimes \xi) = 0$.
\endstatement
\proof Since $\deg V = \deg (V\otimes \xi)$ and both are stable, a nonzero map
between them must be an isomorphism, by standard results on stable bundles.
However $\det (V\otimes \xi) = \det V \otimes \xi ^{\otimes 2} \neq \det V$,
and so the bundles cannot be isomorphic. Thus there is no  nonzero map from
$V$ to
$V\otimes \xi$.
\endproof

\corollary{1.5} Let $\bold F\subset X$ be a scheme-theoretic multiple  fiber
of odd multiplicity $m$ of an elliptic surface, and let $F$ be the reduction
of $\bold F$. Let $\bold V$ be a rank two vector bundle on $\bold F$ whose
restriction $V$ to
$F$ is stable. Then $\dim _\Cee\Hom (\bold V, \bold V)= 1$ and every nonzero
map from $\bold V$ to itself is an isomorphism.
\endstatement
\proof Let $\xi$ be the normal bundle of $F$ in $X$. Thus $\xi$ has order $m$.
For $a>0$, let $aF$ denote the subscheme of $X$ defined by by the ideal sheaf
$\scrO_X(-aF)$. Thus $\bold F=mF$ and there is in general an exact sequence
$$0 \to \xi ^{-a} \to \scrO_{(a+1)F} \to \scrO_{aF} \to 0.$$ Tensor the above
exact sequence by $Hom (\bold V, \bold V) = \bold V\spcheck
\otimes  \bold V$ and take global sections. This gives an exact sequence
$$0 \to \Hom (V, V\otimes \xi ^{-a}) \to \Hom (\bold V|(a+1)F, \bold V|(a+1)F)
\to
\Hom  (\bold V|aF, \bold V|aF).$$ For $a=1$ we have $\dim _\Cee\Hom  (\bold
V|F, \bold V|F)= \dim _\Cee\Hom (V,V) = 1$. For $1\leq a \leq m-1$, $\xi
^{-a}$ is a nontrivial line bundle of odd order. Thus by (1.4) $\Hom (V,
V\otimes \xi ^{-a}) = 0$. It follows that the map
$\Hom (\bold V|(a+1)F, \bold V|(a+1)F) \to
\Hom  (\bold V|aF, \bold V|aF)$ is an injection, so that by induction
$\dim _\Cee \Hom (\bold V|(a+1)F, \bold V|(a+1)F)\leq 1$. On the other hand
multiplication by an element of $H^0(\scrO_{(a+1)F})=\Cee$ defines a nonzero
element of $\Hom (\bold V|(a+1)F, \bold V|(a+1)F)$. Thus
$\dim _\Cee \Hom (\bold V|(a+1)F, \bold V|(a+1)F) = 1$ for all $a\leq m-1$,
and in particular $\dim _\Cee \Hom  (\bold V|mF, \bold V|mF)=1$. \endproof

\section{2. The case of a singular curve.}

Our goal in this section will be to show that the statements of the previous
section hold for rank two vector bundles on singular nodal curves $C$. Let $C$
be an irreducible curve of arithmetic genus one,  which has one node $p$ as a
singularity. Locally analytically, then, $\hat {\Cal O}_{C,p} \cong
\Cee [[x,y]]/(xy)$.  Let   $a\: \tilde C \to C$ be the normalization map, and
let
$p_1$ and $p_2$ be the preimages of the singular point on $\tilde C$. We begin
by giving a preliminary discussion concerning torsion free sheaves on $C$.

\definition{Definition 2.1} A {\sl torsion free rank one sheaf\/} on $C$ is a
coherent sheaf which has rank one at the generic point of $C$ and has no local
sections vanishing on an open set. It is well known that every  torsion free
rank one sheaf on $C$ is either a line bundle or of the form $a_*L$, where $L$
is a line bundle on $\tilde C$. For example, the maximal ideal sheaf  of the
singular point $p$ of $C$ is $a_*L$, where $L= \scrO_{\tilde C}(-p_1-p_2)$,
where $p_1$ and $p_2$ are the preimages of $p$ in $\tilde C$. Here the line
bundle $L$ has degree $-2$ on $\tilde C$. We define the {\sl degree} of a
torsion free rank one sheaf $F$ on $C$ by $\deg F = \chi (F)$. By the
Riemann-Roch theorem on $C$,
$\deg F$ is the usual degree in case $F$ is a line bundle, whereas for $F=
a_*L$, an easy calculation shows that $\chi (F) = \deg L +1$. (Note that, in
case
$p_a(C)$ is arbitrary, we would have to correct by a term $p_a(C) -1$, which
is zero in our case, to get the usual answer for a line bundle.) It is easy to
check that, if $F$ is a rank one torsion free sheaf on $C$ and $L$ is a line
bundle, then $\deg (F\otimes L) = \deg F + \deg L$.  \enddefinition
\medskip

\lemma{2.2} If $F_1$ and $F_2$ are torsion free rank one sheaves on $C$, then
so is $Hom (F_1, F_2)$, and
$$\deg Hom (F_1, F_2) = \cases \deg F_2 - \deg F_1, &\text{if one of the
$F_i$ is a line bundle}\\
\deg F_2 - \deg F_1 + 1,&\text{if neither $F_1$ nor $F_2$ is a line bundle}.
\endcases$$ Finally if $\deg F_2>\deg F_1$ and neither is a line bundle, then
the natural map from $\Hom (F_1, F_2)\otimes \scrO_{C,p}$ to $\Hom
_{\scrO_{C,p}}(\frak m_p, \frak m_p)$ is surjective, where $\frak m_p$ is the
maximal ideal of $\scrO_{C,p}$.
\endstatement
\proof The proof is clear if $F_1$ is a line bundle. Thus we may assume that
$F_1$ is of the form $a_*L$ for a line bundle $L$ on $\tilde C$. First assume
that $F_2$ is a line bundle. An easy calculation shows that $Hom (F_1, F_2) =
a_*(L^{-1} \otimes \scrO_{\tilde C}(-p_1-p_2))\otimes F_2$. This is just the
local calculation $Hom _{\scrO _{C,p}}(a_*\widetilde{\scrO _{C,p}},
\scrO _{C,p}) = \frak m_p$, where $\frak m_p$ is the maximal ideal of $\scrO
_{C,p}$ and the isomorphism is canonical. Thus  $Hom (F_1, F_2)$ is again a
torsion free rank one sheaf and  $$\deg  Hom (F_1, F_2) = -\deg L + 1-2 +\deg
F_2 = \deg F_2 - \deg F_1.$$

Now assume that $F_1 = a_*L_1$ and $F_2 = a_*L_2$. Again using a local
calculation $Hom _{\scrO _{C,p}}(a_*\widetilde{\scrO
_{C,p}},a_*\widetilde{\scrO _{C,p}}) = a_*\widetilde{\scrO _{C,p}}$, where the
isomorphism is also canonical, it is easy to check that $Hom (a_*L_1,  a_*L_2)
= a_*(L_1^{-1}\otimes L_2)$, and so $Hom (F_1, F_2)$ is again a torsion free
rank one sheaf. Moreover
$$\deg Hom (F_1, F_2) = \deg L_2 -\deg L_1 +1= \deg F_2-\deg F_1 + 1.$$ To see
the final statement, again writing $F_1 = a_*L_1$ and $F_2 = a_*L_2$, we have
$Hom (a_*L_1,  a_*L_2) = a_*(L_1^{-1}\otimes L_2)$. Moreover the global
sections of $L_1^{-1}\otimes L_2$ separate the points
$p_1$ and $p_2$. It is then easy to see that the map $\Hom (F_1, F_2)\otimes
\scrO_{C,p} \to \Hom _{\scrO_{C,p}}(\frak m_p, \frak m_p) =
\widetilde{\scrO _{C,p}}$ is surjective.
\endproof

\lemma{2.3} Let $F$ be a torsion free rank one sheaf on $C$. If
$\deg F >0$, or if $\deg F = 0$ and $F$ is not trivial, then $h^0(F) = \deg F$
and $h^1(F) =0$.
 If $\deg F < 0$, or if
$\deg F = 0$ and $F\neq \scrO _C$, then   $h^0(F) = 0$ and $h^1(F) = \deg F$.
\endstatement
\proof If $\deg F \geq 0$ and $F$ is not $\scrO _C$, then the claim that
$h^0(F) = \deg F$ is clear if $F$ is a line bundle and follows from $h^0(a_*L)
=
\deg L +1$ in case $L$ is a line bundle of degree at least $-1$ on $\tilde C
\cong
\Pee ^1$. In this case, since by definition $\deg F = \chi (F) = h^0(F)$, we
must have $h^1(F) = 0$. The proof of the second statement is similar.
\endproof

Next let us consider extensions of torsion free sheaves. The maximal ideal
$\frak m_p$ has the following local resolution, where we set $R = \scrO
_{C,p}$:
$$\dots \to R \oplus R \to R \oplus R  \to \frak m_p \to 0,$$  where the maps
$R \oplus R \to R \oplus R$ alternate between $(\alpha, \beta) \mapsto
(x\alpha, y\beta)$ and $(\alpha, \beta) \mapsto (y\alpha, x\beta)$. A
calculation shows that
$\operatorname{Ext} ^1_R(\frak m_p, \frak m_p)$ has length two. More
intrinsically it is isomorphic to $\scrO _{\tilde C}(-p_1-p_2)/\scrO _{\tilde
C}(-2p_1-2p_2)$. Thus as an $R$-module, $\operatorname{Ext} ^1_R(\frak m_p,
\frak m_p) \cong \tilde R/\tilde{\frak m}_p$, where $\tilde R$ is the
normalization of
$R$ and $\tilde{\frak m}_p = \frak m_p\tilde R$. We can describe the $\tilde
R$-action on $\operatorname{Ext} ^1_R(\frak m_p, \frak m_p)$ more invariantly
as follows: multiplication by $r \in \tilde R$ gives an endomorphism $\frak
m_p\to
\frak m_p$, and hence an action of $\tilde R$ on $\operatorname{Ext} ^1_R
(\frak m_p,  \frak m_p)$. We leave to the reader the straightforward
verification that this action is the same as the action on $\operatorname{Ext}
^1_R(\frak m_p, \frak m_p)$ implicit in the isomorphism  $\operatorname{Ext}
^1_R(\frak m_p, \frak m_p) \cong \tilde R/\tilde{\frak m}_p$ given above.
There is an induced action of the invertible elements ${\tilde R}^*$ on
$(\operatorname{Ext} ^1_R(\frak m_p,
\frak m_p) -0)/\Cee ^* = \Pee ^1$. Since $R^*$ acts trivially, this induces an
action of ${\tilde R}^*/R^* \cong \Cee ^*$ on  $(\operatorname{Ext} ^1_R(\frak
m_p, \frak m_p) -0)/\Cee ^*$. It is easy to see that there are three orbits of
this action: an open orbit isomorphic to $\Cee ^*$ and two closed orbits which
are points in $\Pee ^1$, corresponding to the case of an element $e
\in \operatorname{Ext} ^1_R(\frak m_p, \frak m_p)$ such that $\tilde R\cdot e
\neq  \operatorname{Ext} ^1_R(\frak m_p, \frak m_p)$.

Given an element $e\in \operatorname{Ext} ^1_R(\frak m_p, \frak m_p)$, denote
the corresponding extension of $\frak m_p$ by $\frak m_p$ by $M_e$. Note that
two extensions $M_e$ and $M_{e'}$ such that $e$ and $e'$ lie in the same
$\tilde R^*$-orbit are abstractly isomorphic as $R$-modules, via a diagram of
the form
$$\CD 0 @>>> \frak m_p @>>> M_{re} @>>> \frak m_p @>>> 0 \\ @.    @V=VV
@VVV    @V{\times r}VV   @. \\ 0 @>>> \frak m_p @>>> M_e @>>> \frak m_p @>>>
0,
\endCD$$ where $r\in \tilde R^*$ is such that $re = e'$.

\lemma{2.4} $M_e$ is locally free if and only if the image of $e$ in
$$\Bigl(\operatorname{Ext} ^1_R(\frak m_p,  \frak m_p) -0\Bigr)\big/\Cee ^*$$
is not a closed orbit of $\tilde R$.
\endstatement
\proof Consider the long exact Ext sequence
$$\operatorname{Hom}_R(\frak m_p, \frak m_p) \to  \operatorname{Ext} ^1_R
(\frak m_p, \frak m_p) \to \operatorname{Ext} ^1_R(M_e, \frak m_p).$$ We see
that $\operatorname{Ext} ^1_R(M_e, \frak m_p)$ contains as a submodule
$\operatorname{Ext} ^1_R(\frak m_p, \frak m_p)/\tilde R\cdot e$. Thus if
$\tilde R\cdot e \neq \operatorname{Ext} ^1_R(\frak m_p, \frak m_p)$, then
$\operatorname{Ext} ^1_R(M_e, \frak m_p) \neq 0$ and so $M_e$ is not locally
free. Conversely suppose that the image of $e$ does not lie in one of the
closed orbits.  Since every two extensions in the same orbit are abstractly
isomorphic,  it will suffice to exhibit one locally free extension of $\frak
m_p$ by $\frak m_p$. However we have the obvious surjection $R\oplus R \to
\frak m_p$ given above, and its kernel is easily seen to be isomorphic to
$\frak m_p$ again.
\endproof

We leave as an exercise for the reader the description of the extensions
corresponding to the closed orbits.

Let us also note that, using the resolution above, a short computation shows
that $\operatorname{Ext}^1_R(\frak m_p, R) = 0$. Thus there is no locally free
$R$-module $M$ which sits in an exact sequence
$$0 \to R \to M \to \frak m_p \to 0.$$

Globally, we have the following:

\lemma{2.5} Let $n$ be a positive integer and let  $\delta$ be a line bundle
of  degree one on $C$.
\roster
\item"{(i)}" There is a unique rank two vector bundle $V_{n, \delta}$ on $C$
such  that $\det V_{n, \delta} = \delta$ and such that there is an exact
sequence
$$0 \to F \to V_{n, \delta} \to F' \to 0,$$ where $F$ and $F'$ are torsion
free rank one sheaves of degrees $n$ and $1-n$ respectively, and $F$ and $F'$
are not locally free.
\item"{(ii)}" Let $G$ be a torsion free rank one subsheaf of $V_{n, \delta}$.
Then  either $\deg G \leq -n$ or $G$ is contained in $F$.

\item"{(iii)}" The vector bundle $V_{n, \delta}$ is indecomposable for all $n$
and $\delta$ and $V_{n, \delta} \cong V_{n', \delta'}$ if and only if $n=n'$
and $\delta = \delta '$.
\endroster
\endstatement
\proof To see (i), let $F$ and $F'$ be the unique torsion free rank one
sheaves of degrees $n$ and $1-n$ respectively which are not locally free. Let
us evaluate
$\operatorname{Ext}^1(F', F)$. From the local to global Ext spectral sequence,
there is an exact sequence
$$0 \to H^1(Hom (F', F)) \to \operatorname{Ext}^1(F', F) \to H^0(Ext ^1(F',
F))
\to 0.$$ Now $\chi (Hom (F', F)) = \deg Hom (F', F) = h^0(Hom (F', F))$,
since
$$\deg Hom (F', F) = 2n>0$$ by (2.2) and (2.3). So $H^1(Hom (F', F)) = 0$.
Thus $\operatorname{Ext}^1(F', F) \cong H^0(Ext ^1(F', F))$. Moreover, the set
of all locally free extensions is naturally a principal homogeneous space over
$H^0(\scrO_{\tilde C}^*/\scrO _C^*) = \tilde R^*/R^*$. On the other hand, from
the exact sequence $$0 \to \scrO _C^* \to \scrO_{\tilde C}^*
\to \scrO_{\tilde C}^*/\scrO _C^* \to 0,$$ we have a natural isomorphism
$\operatorname{Pic}^0C \cong H^0(\scrO_{\tilde C}^* /\scrO _C^*)=\tilde
R^*/R^*$.  Let $\partial\: \tilde R^*/R^* \to \operatorname{Pic}^0C$ be the
coboundary map; it is an isomorphism.  Given $e \in \operatorname{Ext}^1(F',
F) \cong H^0(Ext ^1 (F', F))$, let $V_e$ be the extension corresponding to
$e$. A straightforward exercise in the definitions shows that, for $r\in
\tilde R$,
$$\det V_{r\cdot e} = \partial (r) \otimes \det V_e.$$ From this it is clear
that there is a unique extension $V_{n, \delta}$ with determinant $\delta$.

Next we prove (ii). Let $G$ be a torsion free rank one subsheaf, possibly a
line bundle, of $V_{n,
\delta}$ such that $\deg G >-n$. We have an exact sequence $$0 \to
\operatorname{Hom}(G,F) \to \operatorname{Hom}(G, V_{n, \delta}) \to
\operatorname{Hom}(G,F').$$  Moreover  $\operatorname{Hom}(G,F') =
H^0(Hom(G,F'))$. First suppose that either $\deg G > 1-n$ or that $G$ is
locally  free. The torsion free sheaf  $Hom(G,F')$ has degree either $1-n-\deg
G$ or $2-n -
\deg G$, depending on whether $G$ is or is not locally free. In any case it
has degree $\leq 0$ and is not locally free, so that $H^0(Hom(G,F')) =0$, by
(2.2) and (2.3). So every such $G$ is contained in $F$. In the remaining case
where
$\deg G =1-n$ and $G$ is not locally free, then $G=F'$. Since
$\operatorname{Hom}(F',F')
\cong  k^*$, every nonzero homomorphism from $F'$ to itself is an isomorphism.
Thus the exact sequence defining $V_{n, \delta}$ would be split, contrary to
assumption. Hence this last case is impossible.

To see (iii), let $G$ be a torsion free rank one subsheaf of degree at least
$1-n$  such that $V_{n, \delta}/G$ is torsion free. Then by (ii) $G=F$. This
clearly  implies that $V_{n, \delta} \cong V_{n', \delta'}$ if and only if
$n=n'$ and
$\delta = \delta '$ and that $V_{n, \delta}$ is indecomposable. \endproof

\theorem{2.6} Let $C$ be an irreducible curve of arithmetic genus one, which
has one node as a singularity. Let $V$ be a rank two vector bundle on $C$ and
suppose that $\deg \det V = 2e+1$. Then $V$ is one of the following:
\roster
\item"{(i)}" A direct sum of line bundles;
\item"{(ii)}" $L\otimes \Cal F_x$, where $L$ is a line bundle of degree $e$,
$x\in C$ is a smooth point, and
$\Cal F_x$ is the unique nontrivial extension $$0 \to \scrO _C \to \Cal F_x\to
\scrO _C(x)\to 0;$$
\item"{(iii)}" $L\otimes V_{n, \delta}$, where $L$ is a line bundle of degree
$e$  and $V_{n, \delta}$ is the rank two vector bundle  described in
\rom{(2.5)}. In this case, the subsheaf $L\otimes F$, where $F$ is the
subsheaf in the definition of  $V_{n, \delta}$, is the maximal destabilizing
subsheaf.
\endroster
\endstatement
\proof Clearly we may assume that $\deg V = 1$. By the Riemann-Roch theorem,
$h^0(V) \neq 0$. Thus there is a map $\scrO _C \to V$. If this map is the
inclusion of a subbundle, then $V$ is given as an extension
$$0 \to \scrO _C \to V \to \scrO _C(x)\to 0$$ for some smooth point $x\in C$.
Either this extension splits, in which case we are in case (i), or it does not
in which case we are in case (ii).

Now suppose that the map $\scrO_C \to V$ vanishes at some point. There is a
largest rank one subsheaf $F$ of $V$ containing the image of $\scrO_C$, and
$\deg F = n>0$. The quotient $V/F = F'$ is torsion free. If $F$ is a line
bundle, then so is $F'$, since locally $\operatorname{Ext}^1_R(\frak m_p, R) =
0$. In this case $(F')^{-1}\otimes F$ has degree $2n-1 >0$, so that the
extension splits and $V$ is the direct sum of $F$ and $F'$. Hence we are in
case (i). Otherwise $F$ and $F'$ are not locally free. It follows that $V=
V_{n, \delta}$  for $\delta = \det V$, and we are in case (iii). The last
statement in (iii) then follows from the last paragraph of the proof of
(2.5). \endproof

Finally let us show that a statement analogous to (1.3) continues to hold for
the case of a singular curve.

\lemma{2.7} Let $C$ be an irreducible nodal curve of arithmetic genus one.
\roster
\item"{(i)}" Let $V$ be a stable rank two vector bundle over $C$ and suppose
that $\deg V = 2e +1$. Let $d \geq e+1$, and let $\lambda$ be a torsion free
rank  one on $V$ of degree $d$. Then $\dim \operatorname{Hom}(V, \lambda) =
2d-2e-1$, and there exists a surjection from $V$ to $\lambda$. Moreover, if
$\lambda$ is a line bundle on $C$  such that there exists a nonzero map from
$V$ to $\lambda$, then $\deg \lambda \geq e+1$. Finally, if $d=e+1$, then
$H^1(V\spcheck \otimes
\lambda) = 0$.

\item"{(ii)}" Suppose that $V=L_1\oplus L_2$ is a direct sum of line bundles
$L_i$ with $\deg V = 2e+1$ and $\deg L_1 \leq e<\deg L_2$.  Let $\lambda$  be
a  rank one torsion free sheaf on $C$  with $d=\deg \lambda > \deg L_2$. Then
$\dim \operatorname{Hom}(V, \lambda) = 2d-2e-1$, and there exists a surjection
from $V$ to $\lambda$. Moreover, if $\lambda$ is a rank one torsion free sheaf
on $C$ such that there exists a surjection from $V$ to
$\lambda$,  then either $d=\deg \lambda > \deg L_2$ or $\lambda = L_2$ or
$\lambda = L_1$ and $\dim \operatorname{Hom}(V, \lambda) = 2d-2e$.

\item"{(iii)}" Suppose that $V = L\otimes V_{n, \delta}$ for some $n$, where
$L$ is a line bundle of degree $e$ and that $L_2$ is the subsheaf $L\otimes F$
of $V$ of degree $e+n$ corresponding to  the subsheaf $F$ of $V_{n, \delta}$
in the definition of $V_{n, \delta}$ and that
$L_1$ is the quotient $V/L_2$.  Let $\lambda$ be a rank one torsion free sheaf
on
$C$  with $d=\deg \lambda > \deg L_2 = e+n$. Then $\dim \operatorname{Hom}(V,
\lambda) = 2d-2e-1$. Moreover, if there exists a surjection from $V$ to
$\lambda$ then either $\deg \lambda > e+n$ or $\lambda = L_1$ and $\dim
\operatorname{Hom}(V, \lambda) = 1$. \endroster
\endstatement
\proof The proof of (i) and (ii) follows the same lines as the proof of
(1.3),  with minor modifications, given Lemmas 2.2 and 2.3. Let us prove
(iii) in the case where $\lambda$ is not locally free (the proof in the other
case is slightly simpler). By definition there is an exact sequence
$$0 \to L_2 \to V \to L_1 \to 0,$$ where $L_1$ and $L_2$ are not locally free
and $\deg L_2 = e+n$, $\deg L_1 =  1-n+e$. There is a long exact sequence
$$0 \to \operatorname{Hom}(L_1, \lambda ) \to \operatorname{Hom}(V, \lambda)
\to
\operatorname{Hom}(L_2, \lambda ) \to \operatorname{Ext}^1(L_1, \lambda).$$
Moreover, by the long exact sequence for $\operatorname{Ext}$, we have an
exact  sequence $$0 \to H^1((Hom (L_1, \lambda)) \to \operatorname{Ext}^1(L_1,
\lambda)
\to H^0(Ext ^1(L_1, \lambda)).$$ Since $\deg Hom (L_1, \lambda) = d-e+n +1
\leq 0$ and $Hom (L_1, \lambda)$ is not locally free,  $H^1((Hom (L_1,
\lambda)) = 0$ by (2.3). Moreover $H^0(Ext ^1(L_1, \lambda)) = \Cee ^2$. We
claim that the composite map
$\operatorname{Hom}(L_2, \lambda ) \to \operatorname{Ext}^1(L_1, \lambda) \to
H^0(Ext ^1(L_1, \lambda))$ is surjective. Since $\deg L_2>\deg \mu$, the map
$\operatorname{Hom}(L_2, \lambda)\to \operatorname{Hom}_R(\frak  m_p, \frak
m_p)\cong \tilde R$ is onto the quotient $\tilde R/\tilde \frak m_p$ by the
last statement in (2.2). Thus the image of the map $\operatorname{Hom}(L_2,
\lambda )
\to  H^0(Ext ^1(L_1, \lambda))$ contains the orbit $\tilde R\cdot \xi
\subseteq
\operatorname{Ext}^1_R(\frak m_p,\frak m_p)$, where $\xi$ is the extension
class. Since $V$ is locally free, this orbit is all of
$\operatorname{Ext}^1_R(\frak m_p,\frak m_p)$ by the proof of (2.4), and so
the map $\operatorname{Hom}(L_2, \lambda )
\to  H^0(Ext ^1(L_1, \lambda))$ is onto. It follows that  $$\align \dim
\operatorname{Hom}(V, \lambda) &= \operatorname{Hom}(L_1, \lambda) +
\operatorname{Hom}(L_2, \lambda) -2\\ &= d-(e+n) +1 + d-(1-n+e) + 1 - 2 = 2d
-2e-1. \endalign$$ Let us finally consider the case when there is a surjection
from $V$ to $\lambda$.  Let the degree of $\lambda$ be $d+e$. Thus there is a
surjection from $V_{n,  \delta}$ to $\lambda \otimes L^{-1}$, which is of
degree
$d$. Let $G$ be the kernel of the map $V_{n, \delta} \to \lambda \otimes
L^{-1}$. Then $\deg G = 1-d$. By (2.5)(ii), either $1-d \leq -n$ or
$G\subseteq F$. Thus either $d>e$ or $\lambda = L_1$. In the last case, there
is a unique surjection from $V$ to $L_1$ mod scalars, by the proof of
(2.5)(iii). \endproof

\section{3. A Zariski open subset of the moduli space.}

Let $\pi \:S \to \Pee ^1$ be an algebraic elliptic surface of geometric genus
$p_g(S)=p_g$. We shall always assume that the only singular fibers of $\pi$
are either irreducible nodal curves or multiple fibers with smooth reduction.
Denote the multiple fibers by $F_1$ and $F_2$ and suppose that the
multiplicity of $F_i$ is $m_i$. We  shall assume that the multiple fibers lie
over points where the
$j$-invariant of $S$ is unramified. We denote by $J(S)$  the associated
Jacobian elliptic surface or basic elliptic surface. For an integer $n$,
$J^n(S)$ denotes the relative Picard scheme of line bundles on the fibers of
degree $n$ (see for example Section 2 in Part I of [4]). Hence $J(S) = J^0(S)$
and $S = J^1(S)$. If $n$ is relatively prime to $m_1m_2$, then $J^n(S)$ again
has two multiple fibers of multiplicities $m_1$ and $m_2$. We always have
$p_g(J^n(S)) = p_g$. If $\Delta$ is a divisor on $S$, we let $f\cdot \Delta$
denote the {\sl fiber degree\/}, i.e\. the degree of the line bundle $\Delta$
on a smooth fiber $f$. Let
$\operatorname{Pic}^{\text {v}}S$ denote the set of {\sl vertical} divisor
classes, i.e\. the set of divisor classes spanned by the class of a fiber and
the classes of the reductions of the multiple fibers. With our assumptions
$\operatorname{Pic}^{\text {v}}S \cong \Zee \cdot \kappa$, where $m_1m_2\kappa
=f$ (see also [6] Chapter 2 Corollary 2.9). Clearly
$\operatorname{Pic}^{\text {v}}S$ is the kernel of the natural map from
$\operatorname{Pic} S$ to the group of line bundles on the generic fiber. In
general let $\eta = \Spec k(\Pee ^1)$ be the generic point of $\Pee ^1$ and
let $\bar \eta = \Spec
\overline{k(\Pee ^1)}$,  where $\overline{k(\Pee ^1)}$ is the algebraic
closure of $k(\Pee ^1)$. Let
$S_\eta$ be the restriction of $S$ to $\eta$ and $S_{\bar\eta}$ be the
pullback of $S_\eta$ to $\bar\eta$. Define $V_\eta$ to be the restriction of
$V$ to
$S_\eta$ and similarly for $V_{\bar\eta}$.

\definition{Definition 3.1} An ample line bundle $L$ on $S$ is {\sl
$(\Delta, c)$-suitable\/} if for all divisors $D$ on $S$ such that $-D^2 +
D\cdot \Delta \leq c$, either $f\cdot(2D - \Delta) = 0$ or
$$\operatorname{sign} f\cdot (2D - \Delta) =
\operatorname{sign}L\cdot (2D - \Delta).$$

Given the pair $(\Delta, c)$, we set $w = \Delta \bmod 2 \in H^2(S;
\Zee/2\Zee)$ and let $p = \Delta ^2 - 4c$. Thus $(\Delta, c)$ and
$(\Delta ',c')$ correspond to the same values of $w$ and $p$ if and only if
$\Delta ' = \Delta + 2F$ for some divisor class $F$ and $c' = c + \Delta
\cdot F + F^2$. An easy calculation shows that  the property of being
$(\Delta, c)$-suitable therefore only depends on the pair $(w,p)$, and  we
will also say that $L$ is
$(w, p)$-suitable.
\enddefinition
\medskip

We have the following, which is Lemma 3.3 in Part I of [4]:

\lemma{3.2} For all pairs $(\Delta, c)$, $(\Delta, c)$-suitable ample line
bundles exist. \qed
\endstatement
\medskip

\definition{Definition 3.3} Let $\Delta$ be a divisor on $S$ and $c$ an
integer. Fix a $(\Delta, c)$-suitable line bundle $L$. We denote by $\frak
M(\Delta, c)$ the moduli space of equivalence classes of $L$-stable rank two
vector bundles $V$  on $S$ with $c_1(V) = \Delta$ and $c_2(V) = c$. Here $V_1$
and $V_2$ are {\sl equivalent\/} if there exists a line bundle $\scrO _S(D)$
such that $V_1$ is isomorphic to $V_2 \otimes  \scrO _S(D)$. In particular,
since $\det V_1 = \det V_2$, the divisor $2D$ is linearly equivalent to zero,
and in fact $V_1$ and $V_2$ must be isomorphic since there is no 2-torsion in
$\Pic S$. As the notation suggests and as we shall shortly show, the scheme
$\frak M(\Delta, c)$ does not depend on the choice of the  $(\Delta,
c)$-suitable line bundle $L$.

Given a divisor $\Delta$ on $S$ and an integer $c$, we let $w = \Delta \mod 2$
and $p = \Delta ^2 - 4c$. The moduli space $\frak M(\Delta, c)$ only depends
on $w$ and $p$ and we shall also denote it by $\frak M(w,p)$.
\enddefinition

Now fix an odd integer $2e+1$. We shall consider rank two vector bundles
$V$ such that the line bundle $\det V$ has fiber degree $2e+1$. However, it
will be convenient not to fix the determinant of $V$.  In this section we
shall show that the moduli space $\frak  M(w,p)$ is smooth and irreducible,
and shall describe a Zariski open and dense subset of it explicitly. The
basic idea is to show first that there is a largest integer $p_0$ such that
$\frak M(w,p_0)$ in nonempty and that there is a unique element in $\frak
M(w,p_0)$, corresponding to the bundle $V_0$. For all other $p<p_0$, the
bundles in $\frak M(w,p)$ are obtained by elementary modifications of $V_0$
along fibers. Let us begin by recalling the following result  (Corollary 4.4
in Part I of [4]):

\theorem{3.4} Let $V$ be a rank two bundle with $\det V = \Delta$ and $c_2(V)
= c$. Suppose that $\det V$ has fiber degree $2e+1$. Let $L$ be a $(\Delta,
c)$-suitable ample line bundle, and suppose that $V$ is $L$-stable. Then there
exists a Zariski open subset $U$ of $\Pee ^1$ such that, if $f$ is a fiber of
$\pi$ corresponding to a point of $U$, then $f$ is smooth and  $V|f$ is
stable. Conversely, if there exists a smooth  fiber $f$ such that $V|f$ is
stable, then $V$ is $L$-stable for every $(\Delta,c)$-suitable ample line
bundle $L$. \qed\endstatement
\medskip

Next we show that there exist bundles satisfying the hypotheses of Theorem
3.4:

\lemma{3.5} Let $\delta$ be a line bundle on the generic fiber $S_\eta$ of odd
degree $2e+1$. Then there exists a rank two vector bundle $V$ such that the
restriction of  $\det V$ to $S_\eta$ is $\delta$ and such that there exists a
smooth fiber $f$ for which the restriction $V|f$ is stable.
\endstatement
\proof Let $\Delta _0$ be a line bundle on $S$ restricting to $\delta$ on
$S_\eta$. Fix a smooth fiber $f$. By (1.1) there exists a stable bundle $E$ on
$f$  with determinant equal to $\Delta _0|f$. Let $H$ be a line bundle on $S$
such that $\deg (H|f) \geq e+1$. Then by (1.3) there is a surjection from
$E$ to $H|f$, and thus $E$ is given as an extension
$$0 \to (H^{-1}\otimes \Delta _0)|f) \to E \to (H|f) \to 0.$$ This extension
corresponds to a class in $H^1(f; (H^{\otimes -2}\otimes \Delta _0)|f)$.  We
would like to lift this exact sequence to an exact sequence on
$S$. Of course, we can replace $\Delta _0$ by $\Delta _0 + Nf$ for an integer
$N$ and get the same restriction to $f$. It suffices to show that, for some
$N$, the map
$$H^1(S; H^{\otimes -2}\otimes \Delta _0 \otimes \scrO _S(Nf)) \to H^1(f;
(H|f)^{\otimes -2} \otimes \Delta _0)$$ is surjective. The cokernel of this
map is contained in
$$H^2(S;H^{-2}\otimes \Delta _0 \otimes \scrO_S((N-1)f)) = H^0(S; H^2\otimes
\Delta _0^{-1} \otimes \scrO _S((-N+1)f\otimes K_S)^*.$$  Clearly $H^0(S;
H^2\otimes
\Delta _0^{-1} \otimes \scrO _S((-N+1)f\otimes K_S) =0$ if $N\gg 0$, and thus
there is an extension on $S$ inducing $E$. \endproof

\noindent {\bf Note.} We could also have proved (3.5) by descent theory.
\medskip

Before we state the next lemma, recall that a stable vector bundle $V$ on $S$
is {\sl good} if $H^2(S; \operatorname{ad}V) =0$. This means that $V$ is a
smooth point of the moduli space, which has dimension
$-p_1(\operatorname{ad}V) - 3\chi  (\scrO _S)$. Thus the content of the next
lemma is that the moduli space is always smooth of the expected dimension.

\lemma{3.6} Let $V$ be a rank two bundle on $S$ such that the restriction of
$V$ to the generic fiber is stable. Then $V$ is good.
\endstatement
\proof By Serre duality, $H^2(S; \operatorname{ad}V) =0$ if and only if
$H^0(\operatorname{ad}V\otimes K_S) = 0$. A section $\varphi$ of
$H^0(\operatorname{ad}V\otimes K_S)$ gives a trace free endomorphism of
$V_{\bar \eta}$ (since $K_S$ has trivial restriction to the generic fiber).
But $V_{\bar \eta}$ is simple, so that $\varphi$ has trivial restriction to
the generic fiber. Hence $\varphi = 0$. \endproof

\lemma{3.7} Let $V_1$ and $V_2$ be rank two bundles on $S$ whose restrictions
to the generic fibers are stable and have the same determinant
\rom(as a line bundle on $S_\eta$\rom). Then there exists a divisor $D$ on
$S$, lying in $\operatorname{Pic}^{\text{v}}S$, and an inclusion $V_1\otimes
\scrO _S(D) \subseteq V_2$. Moreover for an appropriate choice of $D$ we have
an exact sequence $$0\to V_1\otimes \scrO _S(D) \to V_2 \to Q \to 0,$$ where
$Q$ is supported on fibers or reductions of fibers and the map
$V_1\otimes \scrO _S(D) \to V_2$ does not vanish on any fiber.
\endstatement
\proof By assumption
$V_1$ and $V_2$ have isomorphic restrictions to $S_\eta$. An isomorphism
between these extends to give a map $V_1\otimes \scrO_S(D_1) \to V_2 \otimes
\scrO_S(D_2)$, where $D_i$ have trivial restriction to the generic fiber.
Twisting  gives a map $\varphi \:V_1 \otimes \scrO_S(D') \to V_2$, where $D'$
has trivial restriction to the generic fiber. By construction $\varphi$ is an
isomorphism on the generic fiber, so $\varphi$ is an inclusion. The
determinant $\det \varphi$ is a nonzero section of $\det V_1^{-1} \otimes
\scrO _S(-2 D')\otimes \det V_2$, which restricts trivially to the generic
fiber. Thus $\det V_1^{-1} \otimes \scrO _S(-2 D')\otimes \det V_2 =
\scrO_S(\sum _in_iF_i + nf)$, where the $F_i$ are the multiple fibers, $f$ is
a general fiber and $n_i$, $n$ are $\geq 0$. Here $Q$ has support whose
reduction is the sum of the $F_i$ for which $n_i \neq 0$ plus some smooth
fibers. If
$\varphi$ vanishes identically on a fiber or fiber component $F$, then it
factors:
$$V_1 \otimes \scrO _S(D')\subset V_1 \otimes \scrO _S(D'+F) \to V_2.$$ So
after enlarging $D'$ to a new divisor $D$ we can assume that this doesn't
happen. Thus $D$ is as desired.
\endproof

\corollary{3.8} Let $V_1$ and $V_2$ be two rank two bundles on $S$ with the
following property: for every curve $F$ which is a reduced fiber or the
reduction of a multiple fiber, the restriction of $V_i$ to
$F$ is stable. Then there exists a divisor $D\in
\operatorname{Pic}^{\text{v}}S$ such that $V_2 = V_1\otimes \scrO _S(D)$.
\endstatement
\proof Find a nonzero map $\varphi \: V_1\otimes \scrO _S(D) \to V_2$ which
does not vanish on $F$ for every $F$ the reduction of a fiber, via (3.7). For
all $F$, $V_1\otimes \scrO _S(D)|F$ and $ V_2|F$ are stable bundles of the
same degree and $\varphi|F$ is a nonzero map between them. Thus $\varphi|F$
is an isomorphism for all $F$ and so $\varphi $ is an isomorphism as well.
\endproof

\corollary{3.9} Suppose that $V_0$ is a rank two vector bundle satisfying the
hypotheses of \rom{(3.8)}: the restriction $V|F$ is stable for every reduction
$F$ of a fiber component. Let $\Delta = \det V_0$ and $c = c_2(V_0)$. Then
$\frak M (\Delta, c)$ consists of a single reduced point corresponding to the
bundle $V_0$. Thus necessarily $p_1(\operatorname{ad}V_0) =  -3\chi (\scrO
_S)$. \endstatement
\proof If $V'$ is another such, $V' = V_0 \otimes \scrO _S(D)$, and so $V'$
and
$V_0$ are equivalent. By (3.6) $V_0$ is good. Thus  $\frak M (\Delta, c)$ is a
single reduced point. Moreover the dimension of  $\frak M (\Delta, c)$ is
$-p_1(\operatorname{ad}V_0) - 3\chi (\scrO _S) = 0$, and so
$p_1(\operatorname{ad}V_0) = -3\chi (\scrO _S)$. \endproof

Next we establish the existence of such a $V_0$. Before we do so let us pause
to  record the following lemma.

\lemma{3.10} Let
$$0 \to V_1 \to V_2 \to Q\to 0$$ be an exact sequence of coherent sheaves on
$S$, where $V_1$ and $V_2$ are rank  two vector bundles and $Q=i_*M$ where
$i\: F\to S$ is the inclusion of a reduced fiber or the reduction of a
multiple fiber, and $M$ is a torsion free rank one sheaf on $F$. Then: \roster
\item"{(i)}" We have the following formula for $p_1(\operatorname{ad}V_2)$:
$$\align p_1(\operatorname{ad}V_2) &= p_1(\operatorname{ad}V_1) + 4\Bigl(\deg
M -
\frac{\deg (V_2|F) }2\Bigr)\\ &= p_1(\operatorname{ad}V_1) + 4\Bigl(\deg M -
\frac{\deg (V_1|F) }2\Bigr). \endalign$$
\item"{(ii)}" If we define $Q'$ by the exact sequence $$0 \to V_2\spcheck \to
V_1\spcheck \to Q' \to 0,$$ then $V_i\spcheck \cong V_i \otimes (\det
V_i)^{-1}$ is a twist of $V_i$ and $Q' =  Ext ^1(Q, \Cal O_X)$ is of the form
$i_*M'$, where $M'$ is a torsion free rank one sheaf on $F$ with $\deg M' =
-\deg M$. Finally $M$ is locally free if and only if
$M'$ is locally free. \endroster
\endstatement
\noindent {\it Proof.} The first equality in (i) follows from (0.1) if $M$ is
locally free, with minor modifications in general. To see the second, since
$\det V_2 = \det V_1\otimes \Cal O_S(F)$ and $F^2=0$, we have  $$\deg (V_1|F)
= \deg (\det V_1|F) = \deg (\det V_2|F) = \deg (V_2|F).$$

To prove (ii), note that, after trivializing the bundles $V_i$ in a Zariski
open  set $U$, the map $V_1\to V_2$ is given by a $2\times 2$ matrix $A$ with
coefficients in $\Cal O_U$, and so the dual map corresponds to the matrix
${}^tA$. A local calculation shows that $Q'=i_*M'$, where $M'$ is a torsion
free rank one sheaf on $F$, where $F$ is locally defined by $\det A$, and
that $M'$ is locally free if and only if $M$ is locally free. To calculate
$\deg M'$, use the formula in (i) for $\deg M'$, noting that
$p_1(\operatorname{ad} V_i\spcheck) = p_1(\operatorname{ad} V_i)$ and that
$\deg (V_1\spcheck|F) = - \deg (V_1|F)$. Putting this together gives $$\align
4\deg M' &=p_1(\operatorname{ad}V_1) -p_1(\operatorname{ad}V_2) -2\deg
(V_1|F)\\ &= -4\deg M. \qed
\endalign$$
\medskip

Using the above, we shall show the following, which together with (3.9)
proves (i) of Theorem 2 of the introduction.

\proposition{3.11} Given a line bundle $\delta$ on $S_\eta$ of odd degree,
there exists a rank two bundle $V_0$ on $S$ such that the restriction of
$\det V_0$ to $S_\eta$ is $\delta$ and such that the restriction $V|F$ is
stable  for every reduction $F$ of a fiber component. The rank two bundle
$V_0$ is unique up to equivalence: if $V_1$ is any other bundle with this
property, then there exists a line bundle $\scrO _S(D)$ such that $V_1 \cong
V_0 \otimes \scrO _S(D)$. Moreover $p_1(\operatorname{ad} V_0) \geq
p_1(\operatorname{ad}V)$ for every rank two bundle $V$ such that the
restriction of $\det V$ to $S_\eta$ is $\delta$ and such that there exists a
smooth fiber $f$ for which the restriction $V|f$ is stable, with equality if
and only if $V = V_0 \otimes \scrO _S(D)$.
\endstatement
\proof Begin with $V$ such that $\det V|S_\eta=\delta$ and such that  there
exists a smooth fiber $f$ for which the  restriction $V|f$ is stable. Such
$V$ exist by (3.5). If there exists an $F$ such that $V|F$ is not stable,
then there is a torsion free quotient $Q$ of $V|F$ such that $\deg Q < (\deg
(V|F))/2$. Define $V'$ by the exact sequence
$$0 \to V' \to V \to Q \to 0,$$ where we abusively denote by $Q$ the sheaf
$i_*Q$, where $i$ is the inclusion of $F$ in $S$. Using (i) of (3.10),
$$p_1(\operatorname{ad}V) = p_1(\operatorname{ad}V') + 4\Bigl(\deg Q -
\frac{\deg (V|F) }2\Bigr).$$ Thus $p_1(\operatorname{ad}V') >
p_1(\operatorname{ad} V)$. If $V'$ satisfies the conclusions of (3.11), we are
done. Otherwise repeat this process. At each stage $p_1$ strictly increases.
But
$p_1$ is bounded from above, either from Bogomolov's inequality or using the
fact that the dimension of the moduli space is always $-p_1 -3\chi (\scrO _S)
\geq 0$, by (3.6). Hence this process terminates and gives a $V_0$ as
desired. By (3.8)
$V_0$ is unique up to twisting by a line bundle, and the final statement is
clear from the method of proof. \endproof

Next we shall interpret the proof of (3.11) as saying that every stable
bundle $V$ is obtained from $V_0$ by an appropriate sequence of elementary
modifications.

\definition{Definition 3.12} Let $V$ be a rank two vector bundle on $S$ whose
restriction to the generic fiber is stable. Let $F$ be a fiber on $S$ and $Q$
be a torsion free rank one sheaf on $F$, viewed as a sheaf on $S$. A
surjection $V \to Q$ is {\sl allowable\/} if  $$2\deg Q > \deg (V|F).$$ Thus
if $\deg (V|F) = 2e+1$, then $\deg Q \geq e+1$.  If $W$ is defined as an
elementary modification
$$0\to W \to V \to Q \to 0,$$ then we shall say that the elementary
modification $W$ is {\sl allowable\/} if the surjection $V \to Q$ is
allowable. It then follows from (3.10) that, if $W$ is an allowable
elementary modification of $V$, then $p_1(\operatorname{ad}W) <
p_1(\operatorname{ad}V)$. \enddefinition

Let $Q$ be a rank one torsion free sheaf on a fiber $F$, viewed also as a
sheaf on
$S$, and let $d = \deg Q$. It is an easy consequence of (1.3) and (2.7) that
if $V
\to Q$ is allowable and $\deg (V|f) = 2e+1$, then $d>e$ and either $\dim
\operatorname{Hom}(V,Q) = 2d-2e-1$ or  $\dim \operatorname{Hom}(V,Q) = 2d-2e$
and
$Q$ is a uniquely specified  rank one torsion free sheaf on $F$.

With this said, we have the following:

\proposition{3.13} Let $\delta$ be a line bundle on $S_\eta$ and let $V$ be a
stable   rank two bundle on $S$ such that the restriction of $\det V$ to
$S_\eta$ is $\delta$. Then there is a sequence $V_0$, $V_1$, \dots, $V_n=V$
such that
$V_{i+1}$ is an allowable elementary modification of $V_i$ for $i=1, \dots ,
n-1$. Moreover $2n \leq p_1(\ad V_0)- p_1(\ad V)$. Finally if $V$ is obtained
from $V_0$ from a sequence of allowable elementary modifications then $\dim
\Hom (V, V_0) =1$.
\endstatement
\proof The construction given in the proof of (3.11) is the following: Begin
with  $V$. If $V\neq V_0$, then there is a fiber $F$, a rank one torsion free
sheaf $Q$ on $F$, and an elementary modification $$0 \to V' \to V \to Q \to
0,$$ where $\deg Q < \deg (V|F)/2$. If $V' \neq V_0$, we repeat the process.
So, noting  that $V \cong V\spcheck\otimes \det V$ and using the notation of
(3.10)(ii) it will suffice to show that the dual elementary modification
$$0\to V\spcheck\otimes \det V \to (V')\spcheck\otimes \det V  \to Q'\otimes
\det V
\to 0$$  is allowable, since then we have obtained $V$ as an allowable
elementary modification of $(V')\spcheck\otimes \det V $. But we have $\deg
Q' = -\deg Q$ by (3.10)(ii), and so $$\align \deg ( Q'\otimes \det V) &=
-\deg Q + \deg (V|F)\\ &>
\frac{\deg (V|F)}{2}. \endalign$$ Thus the surjection $(V')\spcheck\otimes
\det V  \to Q'\otimes \det V $ is  allowable.  The  statement about the
number of elementary modifications follows since an  allowable elementary
modification always decreases $p_1$ by a quantity whose absolute value is at
least 2.

Finally let us show that $\dim \Hom (V, V_0) =1$. Since $\dim \Hom (V_0, V_0)
=1$, it is enough by induction on the number of elementary modifications to
show the following: suppose given an exact sequence $$0 \to V_2 \to V_1 \to Q
\to 0,$$ where $\deg Q > \deg (V_i|F)/2$. Then $\Hom (V_1, V_0) \to \Hom (V_2,
V_0)$  is an isomorphism. For simplicity we shall just give the argument in
case
$Q$ is locally free on $F$. In any case $\Hom(Q, V_0)=0$ since $Q$ is a
torsion sheaf and the cokernel of the map is  $$\Ext ^1(Q, V_0) = H^0(Ext
^1(Q, V_0)) = H^0(Q\spcheck \otimes (V_0|F)) = \Hom (Q, V_0|F).$$ Since
$V_0|F$ is stable and
$\deg Q > \deg (V_i|F)/2= \deg (V_0|F)/2$, this last group is zero.
\endproof

Putting all this together, we shall describe a Zariski open subset of the
moduli space. Let us first observe that the moduli space $\frak M(\Delta, c)$
is always good and of dimension
$$4c - \Delta ^2 -3\chi (\scrO _S) = -p - 3\chi (\scrO _S).$$ By the
canonical bundle formula for an elliptic surface,
$$K_S= \scrO _S((p_g-1)f + (m_1-1)F_1+(m_2-1)F_2),$$ where $F_i$ are the
reductions of the multiple fibers. As $m_1$ and $m_2$ are odd,
$$K_S\cdot \Delta \equiv p_g -1\mod 2.$$ By the Wu formula, $\Delta ^2 \equiv
p_g -1\mod 2$ as well. Hence
$$4c - \Delta ^2 -3\chi (\scrO _S) \equiv 0\mod 2,$$ and the dimension of the
moduli space is always an even integer $2t$. Now suppose that $\delta$ is a
line bundle on the generic fiber $S_\eta$ of odd degree. Then there exists a
divisor $\Delta$ on $S$ which restricts to $\delta$ and $\Delta$ is
determined up to a multiple of $\kappa$. Mod 2, the only possibilities are
$\Delta$ and $\Delta +\kappa$. Note that $(\Delta +\kappa)^2 =
\Delta ^2 + 2(\Delta \cdot \kappa) \equiv \Delta ^2 + 2\mod 4$. Thus if we
also fix $\Delta ^2\mod 4$, there is a unique choice of $w = \Delta \mod 2$.
Fix an integer $t\geq 0$ and let $-p = 2t + 3\chi (\scrO _S)$. There is then
a unique class $w\in H^2(S; \Zee/2\Zee)$ with $w^2 \equiv p \mod 4$ such that
$w$ is the mod two reduction of a divisor $\Delta$ which restricts to
$\delta$ on $S_\eta$. Given $\delta$ and $t$, we shall denote the
corresponding moduli space by $\frak M_t$. The following theorem is a more
precise version of Theorem 1 of the Introduction:

\theorem{3.14} In the above notation, $\frak M_t$ is nonempty, smooth and
irreducible, and is birational to
$\operatorname{Sym}^tJ^{e+1}(S)$. More precisely, there exists a Zariski open
and dense subset $U$ of $\frak M_t$ which is isomorphic to the open subset of
$\operatorname{Sym}^tJ^{e+1}(S)$ consisting of $t$ line bundles  $\lambda_1,
\dots , \lambda _t$ of degree $e+1$ lying on smooth \rom(and reduced\rom)
fibers of
$\pi$ such that the images $\pi (\lambda _i)$ are distinct points of $\Pee
^1$,  where we continue to denote by $\pi$ the projection from $J^{e+1}(S)$
to $\Pee ^1$.
\endstatement
\proof Let us describe the set $U$. Given the line bundle $\delta$ on
$S_\eta$, let
$\deg \delta = 2e+1$. If $f$ is a smooth fiber and $\lambda$ is a line bundle
of degree $e+1$ on $f$, the restriction $V_0|f$ sits in an exact sequence
$$0 \to
\mu \to V_0|f \to \lambda  \to 0,$$ where $\mu \otimes \lambda = \delta$.
Once we have fixed $\lambda$, the surjection $V_0|f \to \lambda$ is unique
mod scalars.

Now fix $t$ distinct smooth fibers $f_1, \dots, f_t$ and line bundles
$\lambda _i$ of degree $e+1$ on $f_i$. Let $Q_i$ be the sheaf $\lambda _i$
viewed as a sheaf on $S$ and let $Q = \bigoplus _iQ_i$. We shall consider the
set of vector bundles $V$ described by an exact sequence
$$0 \to V \to V_0 \to Q \to 0.$$ The set of all such $V$ is clearly
parametrized by the open subset $U$ of
$\operatorname{Sym}^tJ^{e+1}(S)$ consisting of $t$ line bundles  $\lambda_1,
\dots , \lambda _t$ lying on smooth (reduced) fibers of $\pi$ such that the
images $\pi (\lambda _i)$ are distinct points of $\Pee ^1$. For such a bundle
$V$,  we also have
$$p_1(\operatorname{ad}V) = p_1(\operatorname{ad}V_0) -2t.$$ We shall first
construct a family of bundles parametrized by $U$ (more precisely, we shall
construct ``universal" bundles  over the product of $S$ with a  finite cover
of $U$), thereby giving a morphism from $U$ to $\frak  M_t$ which is easily
seen to be an open immersion. Finally we shall show that $U$ is in fact dense
in $\frak  M_t$.

\medskip
\noindent {\bf Step I.}  Let $U$ be the open subset of
$\operatorname{Sym}^tJ^{e+1}(S)$ described above, and let $\tilde U$ be
defined as follows:
$$\tilde U = \{\,(\lambda_1, \dots, \lambda _t)\in \bigl(J^{e+1}(S)\bigr)^t\:
\{\lambda_1, \dots, \lambda _t\} \in U\,\}.$$  We shall try to construct a
universal bundle
$\Cal V$ over $S\times \tilde U$ as follows. Let $\Cal Z \subset S \times U$
be defined by
$$\Cal Z = \{\, (p, \{\lambda_1, \dots , \lambda_t\}) \in S \times U : \text{
for some $i$, $\pi (p) = \pi (\lambda_i)$}\,\}.$$  Thus given a point $u =
\{\lambda_1, \dots , \lambda_t\} \in U$,
$$(S\times \{u\}) \cap \Cal Z =
\coprod _{i=1}^t (f_i\times \{u\}),$$  where $f_i$ is the fiber corresponding
to $\lambda_i$. Clearly $\Cal Z$ is a smooth divisor in $S\times U$.
Analogously, we have the pulled back divisor $\tilde {\Cal Z} \subset S\times
\tilde U$. In fact, $\tilde {\Cal Z}$ breaks up into a disjoint union of
divisors $\tilde {\Cal Z}_i$, where for example
$$\tilde {\Cal Z}_1 = \bigl(S\times _{\Pee ^1}J^{e+1}(S)\bigr)\times
J^{e+1}(S)^{t-1},$$ and the other $\tilde {\Cal Z}_i$ are defined by taking
the fiber product over $\Pee ^1$ of $S$ with the $i^{\text {th}}$ factor of
$J^{e+1}(S)^t$.  Thus each $\tilde {\Cal Z}_i$ fibers over $\tilde U$ and the
fiber is an elliptic curve. Let $\rho _i\: \tilde {\Cal Z}_i\to S\times _{\Pee
^1}J^{e+1}(S)$ be the projection. Over $S\times _{\Pee ^1}J^{e+1}(S)$, there
is a relative Poincar\'e bundle $\Cal P_{e+1}$. Actually, $\Cal P_{e+1}$
really just exists locally around sufficiently small neighborhoods of smooth
nonmultiple fibers of $J^{e+1}(S)$, or in irreducible \'etale neighborhoods
$\psi\:\Cal U_0\to J^{e+1}(S)$ of smooth nonmultiple fibers, but we will
write out all the arguments as if there were a global bundle. We shall return
to this point in Section 7. So we should really replace $\tilde U$ by $\tilde
U_0$ defined by $$\tilde U_0 =
\{\,(x_1, \dots, x_t) \in \Cal U_0^t: (\psi(x_1), \dots, \psi(x_t))\in
U\,\}.$$ We can define the divisors $\tilde\Cal Z_i$ on $S\times \tilde U_0$
as well. Thus we have $\rho _i^*\Cal P_{e+1}$, which is a line bundle on
$\tilde {\Cal Z}_i$. By extension, we can view $\rho _i^*\Cal P_{e+1}$ as a
coherent sheaf on $S\times \tilde U_0$.

\lemma{3.15} For every $i$, there is a line bundle $\Cal L_i$ on $\tilde U_0$
with the following property: There is a surjection
$$\pi _1^*V_0\to  \bigoplus _{i=1}^t\Bigl(\rho _i^*\Cal P_{e+1}\otimes \pi
_2^*\Cal L_i \Bigr),$$ and the surjection is unique up to multiplying by the
pullback of a nowhere vanishing function on $\tilde U_0$.
\endstatement
\proof We have
$$\align
\operatorname{Hom}(\pi _1^*V_0 , \bigoplus _{i=1}^t\rho _i^*\Cal P_{e+1}) &=
H^0(\bigl(\pi _1^*V_0\bigr)\spcheck\otimes \Bigl[\bigoplus _{i=1}^t\rho
_i^*\Cal P_{e+1})\Bigr] \\ &= H^0\Bigl(\tilde U_0;  \bigoplus _{i=1}^tR^0\pi
_2{}_*\Bigl(\bigl(\pi _1^*V_0\bigr)\spcheck\otimes  \rho _i^*\Cal
P_{e+1}\Bigr)\Bigr).
\endalign$$  By base change and (1.2), the sheaf
$R^0\pi _2{}_*\Bigl(\bigl(\pi _1^*V_0\bigr)\spcheck\otimes \rho _i^*\Cal
P_{e+1}
\Bigr)$ is a line bundle on $\tilde U_0$, which we denote by $\Cal L_i^{-1}$.
Choosing a nowhere vanishing section of $\scrO _{\tilde U_0}$ gives an
element of
$$\align
\operatorname{Hom}(\pi_1^*V_0, \rho _i^*\Cal P_{e+1}\otimes \pi _2^*\Cal
L_i)&=
 H^0\Bigl(\tilde U_0; R^0\pi _2{}_*\Bigl(\bigl(\pi
_1^*V_0\bigr)\spcheck\otimes
\rho _i^*\Cal P_{e+1}\otimes \pi _2^*\Cal L_i)\Bigr)\Bigr)\\ &= H^0(\tilde
U_0;\Cal L_i^{-1} \otimes \Cal L_i) = H^0(\tilde U_0;  \scrO _{\tilde U_0}).
\endalign$$  Since the
$\tilde Z_i$ are disjoint, we can make such a choice for each $i$ to obtain
the  desired surjection. \endproof

\noindent {\bf Note.} We shall essentially calculate $\Cal L_i$ in Section 7.
\medskip

Making a choice of a surjection from $\pi _1^*V_0$ to $\bigoplus
_{i=1}^t\Bigl(\rho _i^*\Cal P_{e+1}\otimes  \pi _2^*\Cal L_i \Bigr)$ gives a
rank two vector bundle $\Cal V$ over $S\times \tilde U_0$ defined by the
exact  sequence $$0 \to \Cal V \to \pi _1^*V_0\to  \bigoplus
_{i=1}^t\Bigl(\rho _i^*\Cal P_{e+1}\otimes  \pi _2^*\Cal L _i\Bigr)\to 0.$$
Thus there is a morphism $\tilde U _0\to \frak  M_t$. It is easy to see that
this morphism descends to a morphism of schemes $U\to \frak  M_t$ whose image
is the set of bundles described at the beginning of the proof of (3.14).
Clearly the morphism $U\to \frak  M_t$ is injective. By Zariski's Main Theorem
it is an open immersion. This concludes the proof of Step I.

\medskip
\noindent {\bf Step II.} Now we must show that the open set $U$ constructed
above is Zariski dense. To do so, we shall make a standard moduli count which
essentially shows that the closed subset $\frak  M_t-U$ may be parametrized
by a scheme of dimension strictly smaller than $\dim \frak  M_t =2t$.
Consider the set of all allowable elementary modifications of a fixed vector
bundle $V'$ with $\deg (V'|F) = 2e+1$. Thus there is a  reduced fiber or the
reduction of a multiple fiber, say $F$, and a rank one torsion free  sheaf
$Q$ on
$F$ with $\deg Q =d\geq e+1$. By (1.3) and (2.7), there is a surjection from
$V'$ to $Q$, and the set of all such has dimension $2d-2e-1$ or $2d-2e$. Let
$V$ be the kernel of such a surjection. By (3.10),
$$p_1(\operatorname{ad}V') = p_1(\operatorname{ad}V) + 4d-4e-2.$$ Thus the
number of moduli of all $V$ is
$$\align -p_1(\operatorname{ad}V)-3\chi (\scrO _S) &= -p_1(\operatorname{ad}
V') - 3\chi (\scrO _S)+4d-4e-2. \endalign
$$ On the other hand, for $d$ and $V'$ fixed, the above construction depends
on
$2d-2e$ parameters. If $F$ is generic, there is one parameter to choose $F$.
Next, either $\dim \operatorname{Hom}(V', Q) = 2d-2e -1$ or $2d-2e$, and in
this last case $Q$ is fixed. Taking the homomorphisms mod scalars the number
of moduli is either $2d-2e -2$ or $2d-2e-1$. In the first case the choice of
$Q$ is one more parameter, but  not in the second case. Thus we always get
$2d-2e -1$ parameters for the choice of the sheaf $Q$ and the surjection $V'
\to Q$. Adding in the choice of $F$ gives $2d-2e$ moduli. For the above
construction to account for a Zariski open subset of the moduli space, we
clearly must have $V'$ a general point of its moduli space, $F$ a general
fiber, and $2d-2e\geq 4d-4e-2$. It follows that
$d\leq e+1$, and hence since $d>e$ that $d=e+1$. Arguing by induction, we may
assume that $V$ is obtained from $V_0$ by performing successive elementary
modifications along distinct fibers $F_i$ which are smooth and nonmultiple and
with respect to line bundles $\mu _i$ on $F_i$ of degree exactly $e+1$. In
this case $V$ is in the open set $U$ described above. \endproof

\noindent{\bf Notation 3.16.} Given a line bundle $\delta$ on $S_\eta$ and a
nonnegative integer $t$, we let $\frak M_t$ be the moduli space defined prior
to (3.14) of equivalence classes of stable bundles $V$ with $-p_1(\ad V) = 2t
+ 3\chi (\scrO_S)$, such that $w_2(V)$ is the mod two reduction of a divisor
$\Delta$ with $\Delta |S_\eta = \delta$. Thus $\frak M_t$ depends only on
$\delta$ and
$t$. Let $\overline{\frak M}_t$ denote the Gieseker compactification of $\frak
M_t$.

\section{4. The case where $S$ has a section.}

In this section, we shall assume that there is a section $\sigma$ on $S$, so
that
$m_1 = m_2 = 1$. In this case, $\sigma ^2 = -(1+p_g(S))$.  Our goal is to
give a very explicit description of the set of stable bundles on $S$ such
that $\det V$ has the same restriction to the generic fiber as $\sigma$. Thus
$\det V = \sigma +nf$ for some integer $n$. We begin with a lemma on various
cohomology groups which will be used often.

\lemma{4.1} Let $S$ be an elliptic surface with a section $\sigma$. Let  $p_g
= p_g(S)$. \roster
\item"{(i)}" For all integers $a$, $h^0(-\sigma + af) = 0$.
\item"{(ii)}" For all integers $a$,
$$h^1(-\sigma + (p_g+1-a)f) = \cases 0, & a>0\\
 -a+1, & a\leq 0.
\endcases$$
\item"{(iii)}" For all integers $a$,
$$h^2(-\sigma + (p_g+1-a)f) = \cases a-1, & a\geq 2\\
 0, & a\leq 1.
\endcases$$
\endroster
\endstatement
\proof Clearly $h^0(-\sigma + af) = 0$ for all integers $a$. Likewise $R^0\pi
_*
\scrO_S(-\sigma + af) = 0$ for all $a$. In addition $R^2\pi _*
\scrO_S(-\sigma + af) = 0$ for all $a$ since $\pi$ has relative dimension
one.  Thus, from the Leray spectral sequence, we see that
$$\gather H^1(\scrO_S(-\sigma + (p_g+1-a)f)) = H^0(R^1\pi _* \scrO_S(-\sigma
+ (p_g+1-a)f))\\ H^2(\scrO_S(-\sigma + (p_g+1-a)f)) = H^1(R^1\pi _*
\scrO_S(-\sigma + (p_g+1-a)f)).
\endgather$$ Thus we must determine the sheaf $R^1\pi _* \scrO_S(-\sigma +
(p_g+1-a)f)$ on $\Pee ^1$. Now $R^1\pi _* \scrO_S(-\sigma + (p_g+1-a)f) =
R^1\pi _* \scrO_S(-\sigma )\otimes \scrO _{\Pee ^1}(p_g+1-a)$. To calculate
$R^1\pi _* \scrO_S(-\sigma )$, we use the exact sequence
$$0 \to \scrO_S(-\sigma) \to \scrO _S \to \scrO _{\sigma} \to 0.$$ Taking the
long exact sequence for $R^i\pi _*$ gives $R^1\pi _*\scrO_S(-\sigma)
\cong R^1\pi _* \scrO _S$, and, by e.g\. [6] Chapter 1 (3.18), $R^1\pi _*
\scrO _S
\cong \scrO _{\Pee ^1}(-p_g-1)$. So $R^1\pi _* \scrO_S(-\sigma + (p_g+1-a)f)
\cong \scrO _{\Pee ^1}(-a)$, and (ii) and (iii)  follow from the usual
calculations for $\Pee ^1$.
\endproof

Next we shall determe the unique stable vector vector bundle $V_0$ (up to
equivalence) which satisfies  $-p_1(\ad V_0) =3\chi (\scrO_S) $.

\proposition{4.2} Let $S$ be a nodal elliptic surface with a section $\sigma$.
\roster
\item"{(i)}" If $p_g(S)$ is odd, set $k=(1+p_g(S))/2$. Then there is a unique
nonsplit extension
$$0 \to \scrO _S(kf) \to V_0 \to  \scrO _S(\sigma - kf)\to 0,$$ and $\det V =
\sigma$, $-p_1(\ad V_0) = 3\chi (\scrO _S)$, and the restriction of
$V_0$ to every fiber is stable.
\item"{(ii)}" If $p_g(S)$ is even, set $k = p_g(S)/2$. Then there is a unique
nonsplit extension
$$0 \to \scrO _S(kf) \to V_0 \to  \scrO _S(\sigma - (k+1)f)\to 0,$$ and $\det
V = \sigma -f$, $-p_1(\ad V_0) = 3\chi (\scrO _S)$, and the restriction  of
$V_0$ to every fiber is stable.
\endroster
\endstatement
\proof We shall just consider the case where $p_g$ is odd; the other case is
identical. First note that $H^1(S; \scrO _S(-\sigma +2kf)) = H^1(-\sigma +
(p_g+1)f)$ has dimension one, by (4.1)(ii). Thus there is a unique nonsplit
extension up to isomorphism. Clearly $\det V_0 = \sigma$ and $-p_1(\ad V_0) =
4k -\sigma ^2 = 3(1+p_g)$.  Finally we claim that the restriction of $V_0$ to
every fiber is stable. It suffices to show that the restriction of $V_0$ to
every fiber
$f$ is the nontrivial extension of $\scrO_f(p)$ by $\scrO_f$, where $p$ is the
point $\sigma \cdot f$. Thus we must consider the restriction map
$$H^1(S; \scrO _S(-\sigma +2kf) \to H^1(f; \scrO _S(-\sigma +2kf)|f).$$  Its
kernel is  $H^1(S;
\scrO _S(-\sigma +(2k-1)f))= H^1(S; \scrO _S(-\sigma +p_gf))$. Again by
(4.1)(ii) this group is zero, so that $H^1(S; \scrO _S(-\sigma +2kf) \to
H^1(f; \scrO _S(-\sigma +2kf)|f)$ is an injection and hence an isomorphism
since both spaces have dimension one. It follows that $V_0|f$ is stable for
every $f$ and is thus the unique bundle up to equivalence satisfying the
hypotheses of (3.8). \endproof

The bundle $V_0$ (with a slightly different normalization) has been described
independently by Kametani and Sato [8].

Let us now consider the case where $V$ is a stable bundle with $-p_1(\ad V) -
3\chi (\scrO _S) = 2t\geq 0$.

\proposition{4.3} With $S$ as above, let $V$ be a stable rank two vector
bundle  over $S$ such that $\det V=\sigma +nf$ for some $n$ and $-p_1(\ad V)
- 3\chi (\scrO _S) = 2t$. \roster
\item"{(i)}" If $p_g$ is odd and we set $k= (1+p_g)/2$, then, after twisting
by a  line bundle of the form $\scrO_S(af)$, there exist an integer $s$,
$0\leq s \leq t$, and an exact sequence
$$0 \to \scrO _S((k-s)f) \to V \to \scrO _S(\sigma +(-k+s-t)f)\otimes I_Z\to
0.$$ Here $Z$ is a codimension two local complete intersection subscheme of
length $s$. Moreover the inclusion of $\scrO _S((k-s)f)$ into $V$ is
canonically given by the  map $\pi ^*\pi _*V \to V$. If $\varphi\:\scrO_S(af)
\to V$ is a sub-line bundle, then $\varphi$ factors through the inclusion
$\scrO _S((k-s)f) \to V$.
\item"{(ii)}" If $p_g$ is even and we set $k= p_g/2$, then, after twisting by
a  line bundle of the form $\scrO_S(af)$, there exist an integer $s$, $0\leq
s \leq t$, and an exact sequence
$$0 \to \scrO _S((k-s)f) \to V \to \scrO _S(\sigma +(-k-1+s-t)f)\otimes
I_Z\to 0.$$ Here $Z$ is again a codimension two local complete intersection
subscheme of  length $s$. Finally the inclusion of $\scrO _S((k-s)f)$ into
$V$ is canonically given by the map $\pi ^*\pi _*V \to V$, and every nonzero
map $\scrO_S(af) \to V$ factors through $\scrO _S((k-s)f)$.
\endroster
\endstatement
\proof We shall just write down the argument in case $p_g$ is odd. By (3.13),
possibly after twisting, $V$ is obtained from $V_0$ by a sequence of $r\leq t$
allowable elementary modifications. In particular $V$ may be identified with a
subsheaf of $V_0$, and $\det V = \sigma -rf$. There is the map from
$\scrO_S(kf)$ to $V_0$, and clearly the image of the subsheaf
$\scrO_S((k-r)f)$ lies in $V$. Of course, the map  $\scrO_S((k-r)f) \to V$
may vanish along a divisor, but this divisor must necessarily be a union of
at most $r$ fibers. Thus there is an integer $u$ with $0\leq u\leq r$ and an
exact sequence for $V$ of the form
$$0 \to \scrO_S((k-r+u)f)\to V \to \scrO _S(\sigma +(-k-u)f)\otimes I_Z \to
0.$$ Using the condition that $-p_1(\ad V) - 3(p_g+1) =  2t$ gives
$$4\ell(Z) +4(k-r+u)  + (1+p_g) -2r -3(1+p_g)= 2t.$$ Solving, we get
$$-r +2u+2\ell(Z) = t.$$ Let $s = \ell(Z)$.  Twisting the exact sequence by
$\scrO_S(bf)$, where $b= u+\ell(Z) -t$, gives a new  exact sequence (where we
rename $V$ by $V\otimes \scrO_S(bf)$)
$$0 \to \scrO_S((k-s)f)\to V \to \scrO_S(\sigma + (-k+s-t)f)\otimes I_Z \to
0.$$ Clearly $s=\ell(Z) \geq 0$ and since $2s = t+r-2u$ with $u\geq 0$,
$r\leq t$, we  have $s\leq t$. This gives the desired expression of $V$ as an
extension. Since  the restriction of this extension to the generic fiber is
not split, the map
$$R^0\pi _* (\scrO_S(\sigma + (-k+s-t)f)\otimes I_Z) \to R^1\pi _*
\scrO_S((k-s)f)$$ is injective. Thus $\pi _*V = \pi _*\scrO_S((k-s)f) =
\scrO_{\Pee ^1}((k-s))$ and the map $\pi ^*\pi _*V \to V$ is just the
inclusion
$\scrO_S((k-s)f) \to V$. Finally if  $\scrO_S(af) \to V$  is nonzero then $\pi
_*\scrO_S(af) \to \pi _*V = \pi _*\scrO_S((k-s)f)$ is nonzero as well, and the
last assertion of the proposition is then clear. \endproof

There is an analogue of (4.3) for Gieseker semistable torsion free sheaves:

\proposition{4.3$'$} With $S$ and $k$ as above, suppose that $V$ is a  rank
two torsion free sheaf with $c_1(V) = \Delta = \sigma +nf$ for some $n$ and
$c_2(V) =c$ such that $V$ is Gieseker semistable with respect to a $(\Delta,
c)$-suitable line bundle. Suppose that $-p_1(\operatorname{ad}V)-3\chi
(\scrO_S) = 2t$.  Then the restriction of $V$ to  a general fiber of $S$ is
stable, and after twisting by $\scrO_S(af)$ for some $a$ there are
zero-dimensional subschemes $Z_1$  and $Z_2$ of $S$, not necessarily local
complete intersections,  an integer $s$ with $0\leq s\leq t$, and an exact
sequence
$$0 \to \scrO _S((k-s)f)\otimes I_{Z_1} \to V \to \scrO _S(\sigma +(-k+s-t)f)
\otimes I_{Z_2}\to 0,$$ if $p_g=2k-1$ is odd, and
$$0 \to \scrO _S((k-s)f)\otimes I_{Z_1} \to V \to \scrO _S(\sigma
+(-k-1+s-t)f)
\otimes I_{Z_2}\to 0$$ if $p_g=2k$ is even. Moreover $\ell(Z_1) +\ell (Z_2) =
s$.
\endstatement
\proof The double dual $V\spcheck{}\spcheck$ of $V$ is a semistable rank two
vector bundle. Thus it is stable and fits into an exact sequence as in (i) or
(ii) of (4.3). Thus (4.3$'$) follows from manipulations along the lines of
the proof of (4.3).
\endproof

Next let us consider when an extension as in (4.3) can be unstable. For
simplicity we shall just write out the case where $p_g$ is odd.

\proposition{4.4} Suppose that $p_g=2k-1$ is odd and that $V$ is an extension
of the form
$$0 \to \scrO _S((k-s)f) \to V \to \scrO _S(\sigma +(-k+s-t)f)\otimes I_Z\to
0,$$ where $\ell (Z)=s$. Let $s_0$ be the smallest integer such that
$h^0(\scrO_S(s_0f)\otimes I_Z)\neq 0$. Thus $0\leq s_0\leq s$, and $s_0 = 0$
if and only if $s=0$. If $V$ is unstable, then the maximal destabilizing
subbundle is equal to $\scrO_S(\sigma -af)$, where $$t+k-(s-s_0) \leq a \leq
t+k.$$ Thus if $s=s_0$ the only possibility is $\scrO_S(\sigma -(t+k)f)$.
\endstatement
\proof The maximal destabilizing subbundle has a torsion free quotient.
Clearly, it restricts to $\sigma$ on the generic fiber, and thus must be of
the form
$\scrO_S(\sigma -af)$ for some integer $a$. Using the exact sequence
$$0 \to \scrO_S(\sigma -af) \to V \to \scrO_S((a-t)f)\otimes I_{Z'} \to 0,$$
where $Z'$ is a codimension two subscheme, and the fact that
$$\align c_2(V) &= k-s+s = k\\ &= a-t+\ell (Z'),
\endalign$$ we see that $a\leq t+k$. On the other hand, there is a nonzero
map from
$\scrO_S(\sigma -af)$ to $\scrO _S(\sigma +(-k+s-t)f)\otimes I_Z$ and thus a
nonzero section of $\scrO _S((-k+s-t+a)f)\otimes I_Z$. Thus
$$-k+s-t+a \geq s_0,$$ or in other words $a\geq t+k-(s-s_0)$.
\endproof

\corollary{4.5} With assumptions as above, suppose that $Z=\emptyset$, so
that $V$ is an extension  $$0 \to \scrO _S(kf) \to V \to \scrO_S(\sigma -
(k+t)f) \to 0.$$ Then $V$ is stable if and only if it is not the split
extension. In this case we can identify the set of all nonsplit extensions
with $\Sym ^t\sigma$, and an extension $V$ corresponding to $\{p_1, \dots,
p_t\}\in \Sym ^t\sigma$ has unstable restriction to a fiber $f$ if and only
if $p_i\in f$ for some $i$.
\endstatement
\proof As we are in the case $s=0$ of (4.4), if $V$ is unstable then the
destabilizing line bundle is $\scrO_S(\sigma - (k+t)f)$, which splits the
exact sequence. Conversely, if the sequence is not split, then $V$ is stable.

The set of nonsplit extensions of $\scrO_S(\sigma - (k+t)f)$ by $\scrO
_S(kf)$ is parametrized by $\Pee H^1(\scrO_S(-\sigma +(2k+t)f))$. By (4.1)
$H^1(\scrO_S(-\sigma +(2k+t)f)) \cong H^0(R^1\pi _*\scrO_S(-\sigma +(2k+t)f))
= H^0(\Pee ^1; \scrO_{\Pee^1}(t))$. Moreover $\Pee H^0(\Pee ^1;
\scrO_{\Pee^1}(t)) =
\Sym ^t\sigma$ by associating to a section the set of points where it
vanishes. This says that the extension $V$ restricts to the split extension
on a fiber $f$ exactly when the corresponding section of $\scrO_{\Pee^1}(t)$
vanishes at the point of $\Pee ^1$ under $f$.
\endproof

Next we analyze the generic case where $\ell (Z)=t$.

\proposition{4.6} Suppose that $p_g=2k-1$ is odd and that $V$ is an extension
of the form
$$0 \to \scrO _S((k-t)f) \to V \to \scrO _S(\sigma -kf)\otimes I_Z\to 0,$$
where $\ell (Z)=t>0$.
\roster
\item"{(i)}" A locally free extension $V$ as above exists if and only if $Z$
has the Cayley-Bacharach property with respect to $|\sigma +(t-2)f|$.
\item"{(ii)}" Suppose that $s_0=t$ or $t-1$ in the notation of \rom{(4.4)},
and that $\Supp Z\cap \sigma = \emptyset$. Then $\dim \Ext ^1(\scrO _S(\sigma
-kf)
\otimes I_Z,\scrO _S((k-t)f))=1$. A locally free extension exists in this
case if
$s_0=t$.
\item"{(iii)}" Suppose that $Z$ consists of $t$ points lying in distinct
fibers, exactly one of which lies on $\sigma$. Then $\dim \Ext ^1(\scrO
_S(\sigma -kf)
\otimes I_Z,\scrO _S((k-t)f))=1$. A locally free extension exists in this
case if and only if $t=1$.
\item"{(iv)}" If $s_0 \leq t-1$, for example if $Z$ contains two distinct
points lying on the same fiber, then $V$ is unstable.
\item"{(v)}" If $s_0 =t$, then $V$ is stable if no point of $Z$ lies on
$\sigma$. Likewise if $t=1$ and $Z\subset \sigma$, then $V$ is not stable.
\endroster
\endstatement
\proof The long exact sequence for $\Ext$ gives
$$\gather H^1(-\sigma +2k-t)f) \to  \Ext ^1(\scrO _S(\sigma -kf)
\otimes I_Z,\scrO _S((k-t)f))\to \\
\to H^0(\scrO_Z) \to H^2(-\sigma +2k-t)f).
\endgather$$ By (4.1)(ii) $H^1(-\sigma +2k-t)f) =0$. The map $H^0(\scrO_Z)
\to H^2(-\sigma +2k-t)f)$ is dual to the map $H^0(\scrO_S(\sigma +(t-2)f))\to
H^0(\scrO_Z)$ defined by restriction. Thus (i) follows by definition. As for
(ii), since
$\Supp Z\cap \sigma = \emptyset$ and $H^0(\scrO_S(\sigma +(t-2)f)) =
H^0(\scrO_S((t-2)f))$ under the natural inclusion, clearly
$H^0(\scrO_S(\sigma +(t-2)f)\otimes I_Z) = H^0(\scrO_S((t-2)f)\otimes I_Z)$.
By assumption $H^0(\scrO_S((t-2)f)\otimes I_Z)=0$, so that the map
$H^0(\scrO_S (\sigma +(t-2)f))\to H^0(\scrO_Z)$ is an inclusion. But
$h^0(\scrO_S (\sigma +(t-2)f))=t-1$ and $h^0(\scrO_Z) = t$. Thus the cokernel
has dimension one. It is clear that if $s_0 = t$ and $Z$ is reduced, then it
has the Cayley-Bacharach property with respect to $|\sigma +(t-2)f|$. A more
involved argument left to the reader handles the nonreduced case. Thus a
locally free extension exists. This proves (ii), and the proof of (iii) is
similar.

To see (iv), note that if $s_0 \leq t-1$, then there is a section of
$\scrO_S((t-1)f)\otimes I_Z$. Consider the exact sequence
$$0 \to \scrO_S(-\sigma +(2k-1)f) \to Hom (\scrO_S(\sigma - (k+t-1)f), V) \to
\scrO_S((t-1)f)\otimes I_Z \to 0.$$ Since $H^1(\scrO_S(-\sigma +(2k-1)f))=0$
by (4.1), the section of
$\scrO_S((t-1)f)\otimes I_Z $ lifts to define a nonzero homomorphism from
$\scrO_S(\sigma - (k+t-1)f)$ to $V$. Thus $V$ is unstable.

Finally we must prove (v). The bundle $V$ is stable if and only if its
restriction to a general fiber $f$ is stable. Let $f$ be a fiber not meeting
$\Supp Z$. Then there is a natural map $\Ext ^1(\scrO _S(\sigma -kf)
\otimes I_Z,\scrO _S((k-t)f))\to \Ext^1(\scrO_f(p), \scrO_f) =
H^1(\scrO_f(-p))$. This fits into an exact sequence
$$\gather H^0(\scrO_f(-p))\to\Ext ^1(\scrO _S(\sigma -kf)
\otimes I_Z,\scrO _S((k-t-1)f))\to \\
\to \Ext ^1(\scrO _S(\sigma -kf)
\otimes I_Z,\scrO _S((k-t)f))\to H^1(\scrO_f(-p)).
\endgather$$ Since $H^0(\scrO_f(-p)) =0$ and  $h^1(\scrO_f(-p)= \dim \Ext
^1(\scrO _S(\sigma -kf) \otimes I_Z,\scrO _S((k-t)f))=1$, by (ii) and (iii),
it will suffice to show that $\dim \Ext ^1(\scrO _S (\sigma -kf) \otimes
I_Z,\scrO _S((k-t-1)f))\geq1$ if
$\Supp Z \cap \sigma \neq \emptyset$. Now since
$H^1(\scrO_S(-\sigma +(2k-t-1)f))=0$ by (4.1), $\Ext ^1(\scrO _S (\sigma -kf)
\otimes I_Z,\scrO _S((k-t-1)f))$ is dual to the cokernel of the restriction
map $H^0(\scrO_S(\sigma +(t-1)f))\to H^0(\scrO_Z)$. Since $s_0 =t$, by
definition $h^0(\scrO_S((t-1)f)\otimes I_Z)=0$. Thus if $\Supp Z \cap
\sigma =\emptyset$, then $H^0(\scrO_S(\sigma +(t-1)f))$ and
$H^0(\scrO_S((t-1)f))$ have the same image in $H^0(\scrO_Z)$ and
$H^0(\scrO_S((t-1)f))\to H^0(\scrO_Z)$ is injective. As both
$H^0(\scrO_S((t-1)f))$ and $H^0(\scrO_Z)$ have dimension $t$, the map between
them is an isomorphism and the cokernel is zero. It follows that
$V$ restricts to a stable bundle on $f$. Likewise, if $t=1$ and $ Z\subset
\sigma$, then clearly the map $H^0(\scrO_S(\sigma +(t-1)f))\to H^0(\scrO_Z)$
cannot be surjective, and so the cokernel is nonzero. Thus $V$ restricts on
$f$ to an unstable bundle for almost every fiber $f$, so that $V$ is unstable.
\endproof

Let us give another proof for (4.6)(v). Using (4.4) we know that the maximal
destabilizing line bundle, if it exists, must necessarily be of the form
$\scrO_S(\sigma -(t+k)f)$. There is an exact sequence
$$0 \to \scrO_S(-\sigma +2kf) \to Hom (\scrO_S(\sigma -(t+k)f), V) \to
\scrO_S(tf)\otimes I_Z \to 0,$$ and $V$ is unstable if and only if the
nonzero section of $\scrO_S(tf)\otimes I_Z$ lifts to a homomorphism from
$\scrO_S(\sigma -(t+k)f)$ to $V$. The nonzero section of $\scrO_S(tf)\otimes
I_Z$ defines an exact sequence
$$0 \to \scrO_S \to \scrO_S(tf)\otimes I_Z \to Q\to 0.$$ Here if $Z$ consists
of points $z_i$ on distinct fibers $f_i$, then $Q= \bigoplus
_i\scrO_{f_i}(-z_i)$. The coboundary map from $H^0(\scrO_S(tf)\otimes I_Z)$
to
$H^1(\scrO_S(-\sigma +2kf))$ is given by taking cup product of the nonzero
section with the extension class $\xi$ in $\Ext^1(\scrO_S(tf)\otimes
I_Z,\scrO_S(-\sigma +2kf))$ corresponding to $V$. It is easy to see by the
naturality of the pairing that this is the same as taking the image of $\xi$
in
$\Ext ^1(\scrO_S, \scrO_S(-\sigma +2kf))= H^1(\scrO_S(-\sigma +2kf))$ using
the above exact sequence. Taking the long exact Ext sequence and using the
fact that
$H^0(\scrO_S(-\sigma +2kf)) =0$, there is an exact sequence
$$\gather 0\to \Ext ^1(Q,\scrO_S(-\sigma +2kf))\to \Ext ^1(\scrO_S(tf)\otimes
I_Z,\scrO_S(-\sigma +2kf)) \to \\
\to H^1(\scrO_S(-\sigma +2kf)).
\endgather$$ Since $\dim \Ext ^1(\scrO_S(tf)\otimes I_Z,\scrO_S(-\sigma
+2kf))=1$, we see that $\xi \mapsto 0$ if and only if
$\Ext ^1(Q,\scrO_S(-\sigma +2kf)) \neq 0$. So we shall show that $\Ext ^1
(Q,\scrO_S(-\sigma +2kf)) = 0$ if and only if the support of $Z$ does not meet
$\sigma$.

First consider the case where $Z$ consists of points $z_i$ on distinct fibers
$f_i$. Then $Q= \bigoplus _i\scrO_{f_i}(-z_i)$, and standard arguments (cf\.
[6] Chapter 7 Lemma 1.27) show that  $\Ext ^1(\scrO_{f_i}(-z_i),
\scrO_S(-\sigma +2kf)) = H^0(\scrO_{f_i}(z_i-p_i))$, where $p_i = f_i\cap
\sigma$. This group is then zero unless $z_i = p_i$.

We shall briefly outline the argument in the case where $\Supp Z$ is a single
point $z$ supported on a fiber $f$ (the proof in the general case is then
just a matter of notation). In this case $Q= \bigl(\scrO_S(tf)\otimes
I_Z\bigr)/\scrO_S
\cong I_Z/I_{tf}$, where $I_{tf}$ is the ideal of the nonreduced subscheme
$tf$. Moreover the assumption that $s_0=t$ implies that $t$ is the smallest
integer
$s$ such that $x^s\in I_Z$, where $x$ is a local defining function for the
fiber $f$. Our goal now is again to prove that $\Ext ^1(Q, \scrO_S(-\sigma
+2kf)) = 0$.

Now the sheaf $Q$ has a filtration by subsheaves whose successive quotients
are
$$Q_n=I_Z\cap I_{nf}/I_Z\cap I_{(n+1)f}\cong \bigl(I_Z\cap
I_{nf}+I_{(n+1)f}\bigr)/I_{(n+1)f},$$ for $0\leq n\leq t-1$. It is easy to
see that each such quotient is a torsion free rank one
$\scrO_f$-module contained in $I_{nf}/I_{(n+1)f}\cong \scrO_f$. Thus it is a
line bundle on $f$ of strictly negative degree, necessarily of the form
$\scrO_f(-a_nz)$, unless  $(I_Z\cap I_{nf}+I_{(n+1)f})/I_{(n+1)f} =
I_{nf}/I_{(n+1)f}$, or in other words  $I_Z\cap I_{nf}+I_{(n+1)f} =  I_{nf}$.
In this case, in the local ring of $z$ we would have $x^n = h + gx^{n+1}$,
where $x$ is a local defining function for $f$ and $h\in I_Z$. But then $h=
x^n(1-gx)$, so that $x^n\in I_Z$, contradicting the fact that $x^t$ is the
smallest power of $x$ which lies in $I_Z$. Hence $Q_n \cong \scrO_f(-a_nz)$
with $a_n\geq 1$.

A standard argument with Chern classes shows that
$$c_2(Q) = t[z] = \sum _{n=0}^{t-1}c_2(Q_n)=-\sum _{n=0}^{t-1} \deg Q_n,$$
where $c_2(Q)$, $c_2(Q_n)$ are taken in the sense of sheaves on $S$ and $\deg
Q_n$ is in the sense of line bundles on $f$. Thus $\deg Q_n = -1$ for all $n$
and
$Q_n = \scrO _f(-z)$. It follows that $\Ext ^1(Q_n, \scrO_S(-\sigma +2nf))=
H^0(\scrO_f(z-p))$ where $p=\sigma \cap f$.  This group is zero if $z\neq p$
and is nonzero otherwise. Thus $\Ext ^1(Q, \scrO_S(-\sigma +2nf)) = 0$ if
$z\neq p$ and $\Ext ^1(Q, \scrO_S(-\sigma +2nf)) \neq 0$ if $z=p$.
\medskip

We shall now reverse the above constructions and try to find a universal
bundle  in the case where the  dimension of the moduli space is 2 or 4. For
simplicity we shall just consider the case where $p_g$ is odd.
\medskip
\noindent {\bf  The two-dimensional invariant.}
\medskip

 Let $\frak M_1$ denote the moduli space of equivalence classes of stable
rank two bundles  $V$ for which $-p_1(\ad V) -3\chi (\scrO_S) =2$. Thus
$\frak M_1$ is compact. Since $p_g$ is odd, we may fix the determinant of $V$
to be $\sigma -f$. Our goal is to show the following:

\theorem{4.7}  $\frak M_1 \cong S$. Moreover there is a universal bundle
$\Cal V$ over $S\times S$, and
$$p_1(\ad \Cal V)/ [\Sigma] = (2(\sigma \cdot \Sigma) -2p_g(f\cdot \Sigma))f
-4(f\cdot \Sigma)\sigma -4\Sigma.$$ Thus, as $-4\mu (\Sigma) = p_1(\ad \Cal
V)/ [\Sigma]$, we have
$$\mu (\Sigma)^2 = (\Sigma )^2 + (p_g-1)(f\cdot \Sigma)^2.$$
\endstatement
\proof It follows from (4.3), (4.5), and (4.6)(v) that if $V$ is stable
with  $-p_1(\ad V) -3\chi (\scrO_S) =2$ and $c_1(V)=\sigma -f$, then either
there is an exact sequence
$$0 \to \scrO _S((k-1)f) \to V \to \scrO _S(\sigma -kf)\otimes
\frak m_q \to 0$$ with $\frak m_q$ the maximal ideal of a point $q\notin
\sigma$ or there is an nonsplit exact sequence  $$0 \to \scrO _S(kf) \to V \to
\scrO _S(\sigma + (-k-1)f) \to 0. $$ In this case the set of all nonsplit
extensions is isomorphic to $\sigma$. Thus the moduli space $\frak M_1$ is
made up of $S-\sigma$, together with a copy of $\sigma$. To glue these two
pieces, we shall construct a universal bundle over $S\times S$ by taking
extensions and then making an elementary modification. To this end, let $\Bbb
D$ be the diagonal in $S\times S$. Consider the extension $\Cal W$ over
$S\times S$ defined as follows: $$0 \to \pi _1^*\scrO_S((k-1)f) \otimes \pi
_2^*\Cal L \to \Cal W \to \pi _1^*\scrO _S(\sigma -kf)\otimes I_{\Bbb D} \to
0.$$ Here, using the relative Ext sheaves and standard exact sequences we
should take
$$\align
\Cal L^{-1}&= Ext ^1_{\pi _2}(\pi _1^*\scrO _S(\sigma -kf)\otimes I_{\Bbb D},
\pi _1^*\scrO_S((k-1)f)) \\ &= \pi _2{}_*(\det N_{\Bbb D}\otimes \pi
_1^*(\scrO _S(-\sigma + (2k-1)f))\\ &= \pi _2{}_*(\scrO _{\Bbb D}(-(p_g-1)f)
\otimes \pi _1^*(\scrO _S(-\sigma + p_gf)
\\ &= \scrO _S(-\sigma + f).
\endalign$$ With this choice of $\Cal L$, we find that
$\operatorname{Ext}^1(\pi _1^*\scrO _S(\sigma -kf)\otimes I_{\Bbb D}, \pi
_1^*\scrO_S((k-1)f) \otimes \pi _2^*\Cal L)
\cong H^0(\scrO _{\Bbb D})$ and that the unique nontrivial extension is indeed
locally free. This defines $\Cal W$, and an easy computation gives
$$\align c_1(\Cal W) &= \pi _1^*(\sigma -f) + \pi _2^*c_1(\Cal L);\\ p_1(\ad
\Cal W) &= 2\pi _1^*(-\sigma + (2k-1)f)\cdot \pi _2^*(\sigma - f) -4[\Bbb D]
+ \cdots,
\endalign$$ where the omitted terms do not affect  slant product.

The restriction $W$ of $\Cal W$ to the slice $S\times \{q\}$ is the unique
nontrivial extension of $\scrO _S(\sigma -kf)\otimes \frak m_q$ by
$\scrO _S((k-1)f)$. By (4.6)(v) $W$ is stable if  and only if $q$ does not
lie on $\sigma$. To remedy this problem, we shall make an elementary
modification along $S\times \sigma$. Note that, if $W$ is given as an
extension  $$0 \to \scrO _S((k-1)f) \to W \to \scrO _S(\sigma -kf)\otimes
\frak m_q \to 0,$$ where $q\in \sigma$, then the maximal  destabilizing
sub-line bundle of $W$ must be $\scrO _S(\sigma + (-k-1)f)$ by (4.4) and thus
there is an exact sequence $$0 \to \scrO _S(\sigma + (-k-1)f) \to W \to \scrO
_S(kf) \to 0.$$

It follows that $\pi _2{}_*Hom(\Cal W|(S\times \sigma), \pi _1^*\scrO
_S(kf))$ is a line bundle $\Cal M$ and the natural map
$$\Cal W \to i_*(\pi _1^*\scrO _S(kf) \otimes \pi _2^*\Cal M)$$ is
surjective. Thus we can define $\Cal V$ by taking the associated elementary
modification. By construction there is an exact sequence
$$0 \to \Cal V \to \Cal W \to i_*(\pi _1^*\scrO _S(kf) \otimes \pi _2^*\Cal
M) \to 0.$$ Moreover for each $q\in \sigma$ there is an exact sequence
$$0 \to \scrO _S(kf) \to \Cal V|S\times \{q\} \to \scrO _S(\sigma +(-k-1)f)
\to 0.$$ Thus by (4.5) $\Cal V|S\times \{q\}$ is stable provided that this
extension does not split. We state this fact explicitly as a lemma, whose
proof will be deferred until later:

\lemma{4.8} In the above notation, the extension for $\Cal V|S\times \{q\}$
does  not split.
\endstatement
\medskip

Assuming the lemma, the restriction of $\Cal V$ to each slice is stable and
thus
$\Cal V$ defines a morphism from $S$ to the moduli space $\frak M_1$. It is
clear that this morphism is a bijection between two smooth surfaces and is
thus an isomorphism. Moreover, by (0.1),
$$p_1(\ad \Cal V) = p_1(\ad \Cal W) + 2c_1(\Cal W)\cdot [S\times \sigma] +
[S\times \sigma]^2-4i_*c_1(\pi _1^*\scrO _S(kf) \otimes \pi _2^*\Cal M).$$
Plugging in for $c_1(\Cal W)$ and $p_1(\ad \Cal W)$ gives
$$p_1(\ad \Cal V) = 2\pi _1^*(-\sigma + (2k-1)f)\cdot \pi _2^*(\sigma - f) -
4[\Bbb D] +2\pi _1^*(\sigma -f)\cdot \pi _2^*\sigma -4\pi _1^*(kf) \cdot \pi
_2^*\sigma + \cdots,$$ where as usual the omitted terms do not affect slant
product. Thus collecting  terms and taking slant product gives
$$-4\mu (\Sigma) = 2(\sigma\cdot \Sigma)f -2p_g(f\cdot \Sigma)f-4(f\cdot
\Sigma)
\sigma - 4\Sigma,$$ as claimed in the statement of Theorem 4.4. This
concludes the proof of Theorem 4.7. \endproof

\demo{Proof of \rom{(4.8)}} We shall use the criterion (A.4) of the Appendix
and the discussion following it to  see that the extension does not split.
Given
$q\in \sigma$, let $W= \Cal W|S\times \{q\}$. We need to show:
\roster
\item"{(i)}" $ \Hom (\scrO_S(\sigma +(-k-1)f), \scrO_S(kf))=0$.
\item"{(ii)}" The map (coming from the usual long exact sequences)
$$\gather R^0\pi _2{}_*\bigl(\pi _1^*\scrO_S(\sigma -kf) \otimes I_{\Bbb D}
\otimes
\pi _1^* \scrO _S(-\sigma +(k+1)f)\bigr) = R^0\pi _2{}_*\bigl(\pi
_1^*\scrO_S(f)
\otimes I_{\Bbb D}\bigr) \to \\
\to R^1\pi _2{}_*\pi _1^*(\scrO _S((k-1)f) \otimes
\scrO _S(-\sigma +(k+1)f)) = R^1\pi _2{}_*\pi _1^*(\scrO _S(-\sigma + 2kf))
\endgather$$ vanishes simply along $\sigma$.
\item"{(iii)}" $H^1(\scrO_S(f)\otimes \frak m_q)$ is independent of $q$, and
is nonzero  only if $p_g=0$. Moreover $H^2(\scrO_S(-\sigma +2kf))=0$.
\item"{(iv)}" At each point of $\sigma$, the  map $H^1(\scrO _S(-\sigma +
2kf))\to H^1(\scrO _S(kf)\otimes \scrO _S(-\sigma +(k+1)f)) = H^1(-\sigma +
(2k+1)f)$ induced by the map $H^1(\scrO _S(-\sigma + 2kf))\to H^1(W\otimes
\scrO _S(-\sigma +(k+1)f)$ followed by the natural map $H^1(W\otimes
\scrO _S(-\sigma +(k+1)f)\to H^1(\scrO _S(kf)\otimes \scrO _S(-\sigma
+(k+1)f)$ is injective.
\endroster

The statement (i) is clear. To prove (ii), we shall calculate $R^1\pi
_2{}_*(\Cal W \otimes \scrO _S(\sigma +(-k-1)f))$ by an argument similar to
the second proof of (4.6)(v). By base change
$R^0\pi _2{}_*(\pi _1^*\scrO _S(f) \otimes I_{\Bbb D})=\Cal L_1$ is a line
bundle on $S$. From the definition of $\Cal W$ and $\Cal L$ the sheaf
$Ext ^1_{\pi _2}(\pi _1^*\scrO _S(f)\otimes I_{\Bbb D} , \pi _1^*\scrO
_S(-\sigma +2kf))\otimes \Cal L$ is the trivial line bundle. A global section
induces the map
$\Cal L_1 \to R^1\pi _2{}_*\pi _1^*\scrO _S(\sigma +2kf))\otimes  \Cal L$. The
cokernel of this map is a subsheaf of $R^1\pi _2{}_*(\Cal W \otimes \scrO
_S(-\sigma +(k+1)f))$. To determine where the map vanishes, use the exact
sequence
$$0 \to \pi _2^*\Cal L_1 \to \pi _1^*\scrO _S(f)\otimes I_{\Bbb D}\to \Cal P
\to 0.$$  Here the map
$\pi _2^*\Cal L_1 \to \pi _1^*\scrO _S(f)\otimes I_{\Bbb D}$ is the natural
one and a calculation in local coordinates shows that it vanishes simply
along $D = S\times _{\Pee ^1}S\subset S\times S$. It follows that, up to a
line bundle pulled back from the second factor $\Cal P = \scrO_D(-\Bbb D)$.
Thus $\Cal P$ is up to sign a Poincar\'e bundle.

Now $\Ext ^2(\scrO _S(f)\otimes \frak m_q,\scrO _S(-\sigma +2kf))=0$ since
$H^2(\scrO _S (-\sigma +(2k-1)f))=0$. Thus $Ext ^2_{\pi _2}(\pi _1^*\scrO
_S(f)\otimes I_{\Bbb D} , \pi _1^*\scrO _S(-\sigma +2kf))=0$ and there is an
exact sequence
$$\gather Ext _{\pi _2}^1(\Cal P, \pi _1^*\scrO _S(-\sigma +2kf))\to  Ext
^1_{\pi _2}(\pi _1^*\scrO _S(f)\otimes  I_{\Bbb D} , \pi _1^*\scrO _S(-\sigma
+2kf)) \to \\
\to R^1\pi _2{}_*\pi _1^*\scrO _S(\sigma +2kf))\to Ext _{\pi _2}^2(\Cal P,
\pi _1^*\scrO _S(-\sigma +2kf))\to 0.
\endgather$$ It follows from the naturality of the pairings involved that the
image of the map
$$Ext ^1_{\pi _2}(\pi _1^*\scrO _S(f)\otimes  I_{\Bbb D} , \pi _1^*\scrO
_S(-\sigma +2kf)) \to R^1\pi _2{}_*\pi _1^*\scrO _S (\sigma +2kf))$$ is, up
to a twist by the line bundle $\Cal L$, the image of $\Cal L_1$.

Note that the restriction of $\Cal P\spcheck \otimes  \pi _1^*\scrO
_S(-\sigma +2kf)$ to $\pi _2^{-1}(q)\subset D$,  where $q$ is any point
except the singular point on a singular fiber, is $\scrO _f(p-q)$, where $f$
is the fiber containing $q$ and $p = f\cap \sigma$.  Ignoring the possible
double points of $D$, we have by standard arguments
$$\gather Ext ^1_{\pi _2}(\pi _1^*\scrO _S(f) \otimes I_{\Bbb D} , \pi
_1^*\scrO _S(-\sigma +2kf))=\pi _2{}_*(\Cal P\spcheck \otimes  \pi _1^*\scrO
_S(-\sigma +2kf))\\ Ext ^2_{\pi _2}(\pi _1^*\scrO _S(f) \otimes I_{\Bbb D} ,
\pi _1^*\scrO _S(-\sigma +2kf))= R^1\pi _2{}_*(\Cal P\spcheck \otimes  \pi
_1^*\scrO _S(-\sigma +2kf))
\endgather$$  (where $\Cal P\spcheck$ means that the dual is taken as a line
bundle on $D$). Thus  $\pi _2{}_*(\Cal P\spcheck \otimes  \pi _1^*\scrO
_S(-\sigma +2kf)) =0$ and  $R^1\pi _2{}_*(\Cal P\spcheck \otimes  \pi
_1^*\scrO _S(-\sigma +2kf))$ is supported on
$\sigma$. To calculate its length, we have (again ignoring the double points
of
$D$ which will not cause trouble) an exact sequence
$$0 \to \scrO_D(\Bbb D -\pi _1^*\sigma +\pi _1^*(2kf)) \to \scrO_D(\Bbb D +
\pi _1^*(2kf)) \to  \scrO_D(\Bbb D + \pi_1^*(2kf))|\pi _1^*\sigma \cap D\to
0.$$ Now $\pi _1^*\sigma \cap D \cong S$ via $\pi _2$ and  under this
isomorphism
$\scrO_D(\Bbb D + \pi_1^*(2kf))|\pi _1^*\sigma \cap D \cong \scrO_S(\sigma
+2kf)$. The map
$$R^0\pi _2{}_*\scrO_D(\pi_1^*(2kf))
\to R^0\pi _2{}_*\scrO_D(\Bbb D+\pi_1^*(2kf))$$ is an isomorphism, since the
induced map on $H^0$'s for the restriction to each fiber of $\pi _2$ is an
isomorphism. Using the exact sequence
$$0\to \scrO_D(\pi_1^*(-\sigma + 2kf)) \to \scrO_D(\pi_1^*(2kf)) \to
\scrO_S(2kf) \to 0,$$ it follows that the image of
$R^0\pi _2{}_*\scrO_D(\pi_1^*(2kf)) = R^0\pi _2{}_*\scrO_D(\Bbb
D+\pi_1^*(2kf))$ in $$R^0\pi _2{}_*\scrO_D(\Bbb D +
\pi_1^*(2kf))|\pi _1^*\sigma \cap D= \scrO_S(\sigma +2kf)$$ is just the image
of
$\scrO_S(2kf)$ in $\scrO_S(\sigma +2kf)$. Thus this map vanishes simply along
$\sigma$, and its cokernel, which is $$R^1\pi _2{}_*\scrO_D(\Bbb D -\pi
_1^*\sigma +\pi _1^*(2kf)) = R^1\pi _2{}_* (\Cal P\spcheck \otimes  \pi
_1^*\scrO _S(-\sigma +2kf)),$$ is a line bundle on $\sigma$. It follows that
the map of line bundles
$$Ext ^1_{\pi _2}(\pi _1^*\scrO _S(f)\otimes I_{\Bbb D} , \pi _1^*\scrO
_S(-\sigma +2kf))\otimes \Cal L \to R^1\pi _2{}_*\pi _1^*\scrO _S(\sigma
+2kf))\otimes \Cal L$$  vanishes simply along $\sigma$, so that we are in the
situation of (A.4): the cokernel contributes torsion of length one.

To see that the above is exactly the torsion in $R^1\pi _2{}_*(\Cal W\otimes
\scrO_S(-\sigma +(k+1)f))$ follows from (iii), as in the discussion after
(A.4). To see (iii), use the exact sequence
$$0\to \scrO_S(f)\otimes \frak m_q \to \scrO_S(f) \to \Cee _q \to 0.$$ The
long exact cohomology sequence shows that $H^1(\scrO_S(f)\otimes \frak m_q)
\cong H^1(\scrO_S(f))$. It is easy to see that this last group is zero if $p_g
>0$ and has dimension one if $p_g=0$ (and in any case its dimension is
obviously independent of $q$). Finally $H^2(\scrO_S(-\sigma +2kf))=0$ by
(4.1). Thus we have identified the torsion in $R^1\pi _2{}_*(\Cal W\otimes
\scrO_S(-\sigma +(k+1)f))$, compatibly with base change.

 We finally need to check that the induced map $$H^1(\scrO_S(-\sigma +2kf))
\to H^1(\scrO_S(-\sigma +(2k+1)f))$$ is injective. But this map is induced
from the composition of the map of sheaves $\scrO_S(-\sigma +2kf) \to W\otimes
\scrO_S(-\sigma +(k+1)f)$ together with the map $W\otimes  \scrO_S(-\sigma
+(k+1)f) \to \scrO_S(-\sigma +(2k+1)f)$. This composition is then a nonzero
map from $\scrO_S(-\sigma +2kf)$ to $\scrO_S(-\sigma +(2k+1)f)$ and so fits
into  an exact sequence
$$0 \to \scrO_S(-\sigma +2kf) \to \scrO_S(-\sigma +(2k+1)f) \to
\scrO_f(-p) \to 0.$$  Since $H^0(\scrO_f(-p))=0$, the map  the map
$H^1(\scrO_S(-\sigma +2kf)) \to H^1(\scrO_S(-\sigma +(2k+1)f))$ is injective.
Thus the extension class for $\Cal V|S\times \{q\}$ is nonzero, and we are
done.
\endproof

\medskip
\noindent {\bf The four-dimensional invariant.}
\medskip

We again assume that $p_g$ is odd and list the possible types of extensions
for a  stable bundle. The generic case (Type 1) is where there exists a
codimension two subscheme
$Z$ with $\ell(Z)=2$ and an exact sequence
$$0 \to \scrO _S((k-2)f) \to V \to \scrO _S(\sigma -kf) \otimes I_Z \to 0.\tag
Type 1$$  Other possibilities (Types 2 and 3 respectively) are
$$\gather 0 \to \scrO _S((k-1)f) \to V \to \scrO _S(\sigma +(-k-1)f) \otimes
\frak m_q \to 0;\tag Type 2\\ 0 \to \scrO _S(kf) \to V \to \scrO _S(\sigma
+(-k-2)f) \to 0.
\tag Type 3
\endgather$$ Here $\frak m_q$ is the maximal ideal of a point $q$. Finally
there is also the  case where $V$ is not locally free. In this case the
double dual of $V$ fits into an extension $$0 \to \scrO _S((k-1)f) \to
V\spcheck{}\spcheck  \to \scrO _S(\sigma +(-k-1)f)  \to 0$$ which must be
nonsplit if $V$ is to be stable, in which case
$V\spcheck{}\spcheck$ is just a twist of $V_0$. One possibility is that $V$ is
given as the unique non-locally free extension of $\scrO _S(\sigma +(-k-1)f)
\otimes  \frak m_q$ by $\scrO _S((k-1)f)$ as in the second exact sequence
above. The remaining possibility (Type 4) is that $V$  is given as an
extension
$$0 \to \scrO _S((k-1)f)\otimes \frak m_q \to V \to \scrO _S(\sigma +(-k-1)f)
\to 0.\tag Type 4$$  For a fixed $q$, the set of all such extensions is
parametrized by a $\Pee ^1$, one point of which correspond to a $V$ such that
$V\spcheck{}\spcheck$ is unstable.

Our goal here is to give a very brief sketch of the following, where we use
the notation of the introduction for divisors on $\operatorname{Hilb}^2S$:

\theorem{4.9} The moduli space $\frak M_2$ of dimension $4$ is isomorphic to
$\operatorname{Hilb}^2S$, and for all $\Sigma \in H_2(S)$,
$$\mu (\Sigma) = D_{\alpha _2} -\bigl((f\cdot \Sigma)/2\bigr)E,$$ where,
setting  $\alpha _1 = \mu _1(\Sigma)$ to be the class computed by the
$\mu$-map for the two-dimensional invariant,
$$\align
\alpha _2 &= \Sigma +\bigl(-(\sigma \cdot \Sigma) + (p_g+1)(f\cdot
\Sigma)/2\bigr)f + (f\cdot \Sigma)\sigma\\ &= \alpha _1 + \bigl((f\cdot
\Sigma)/2\bigr)f.
\endalign$$
\endstatement

Thus an easy calculation using the multiplication table for
$\operatorname{Hilb}^2S$ gives the following:

\corollary{4.10} In the above notation,
$$\mu (\Sigma )^4 = 3(\Sigma ^2)^2 + 6(p_g-1)(\Sigma ^2)(f\cdot \Sigma)^2 +
\bigl[3(p_g+1)(p_g-1)-8(p_g-1)\bigr](f\cdot \Sigma)^4.$$
\endstatement

We shall not give a complete proof of (4.9) here, but shall outline the
argument and prove some statements which will be used later. In Sections 9
and 10, we shall prove a more general statement which will imply (4.9).

We begin as before by analyzing the generic case, Type 1. Let $Z$ be a
codimension two subscheme of $S$ with $\ell (Z) =2$. Let $D_\sigma$ be the
effective divisor of $\operatorname{Hilb}^2S$ which is the closure of the
locus of pairs $\{z_1, z_2\}$ where $z_1\in \sigma$. Then arguing as in the
proof of (4.6)(i)--(iii), we see that
 $$\dim \Ext ^1(\scrO _S(\sigma -kf)\otimes I_Z, \scrO _S((k-2)f) = \cases 1,
&
\text{if $Z \notin \Sym ^2\sigma \subset \operatorname{Hilb}^2S$;}\\ 2,
&\text{otherwise.}
\endcases$$ In case $Z\notin D_\sigma$, the unique extension class mod
scalars corresponds to a locally free extension. If $Z\in D_\sigma - \Sym
^2\sigma$, then the unique nontrivial extension is not locally free. If $Z\in
\Sym ^2\sigma$, then there exist locally free extensions.

Next we must analyze when a locally free extension is stable. Let $\Cal D$ be
the irreducible divisor in $\operatorname{Hilb}^2S$ corresponding to the
divisor
$S\times _{\Pee ^1}S \subset S\times S$. Equivalently
$$\Cal D=\{\, Z\in \operatorname{Hilb}^2S\mid h^0(\scrO_S(f)\otimes I_Z)  =
1\,\}.$$  The divisor $\Cal D$ is smooth, although $S\times _{\Pee ^1}S$ is
singular at the finitely many pairs of points $(x,x)$, where $x$ is a  double
point of a singular fiber. One way to see this is as follows. The divisor
$S\times _{\Pee ^1}S$ has ordinary threefold double points at the
singularities. Moreover it contains the diagonal $\Bbb D\subset S\times S$,
which is smooth and passes through the double points. It is well-known (and
easy to check) that the blowup of a threefold double point  $(xy-zw)$ along a
subvariety of the form
$(x-z, y-w)$ gives a small resolution of the singularity. Thus the proper
transform of $S\times _{\Pee ^1}S$ in the blowup of $S\times S$ along $\Bbb
D$ is smooth, and $\Cal D$ is the quotient of this proper transform by an
involution whose fixed point set is smooth (it is $\Bbb D$). Thus $\Bbb D$
is  smooth.

\lemma{4.11} Let $V$ be a vector bundle given by an extension
$$0 \to \scrO _S((k-2)f) \to V \to \scrO _S(\sigma -kf) \otimes I_Z \to 0,$$
where $\ell (Z)=2$. Then $V$ is not stable if and only if either $Z\in  \Sym
^2\sigma$ or $Z\in \Cal D$. If $Z\in \Cal D$, then the maximal destabilizing
sub-line bundle is $\scrO_S(\sigma +(-k-1)f)$ and there is an exact sequence
$$0 \to \scrO_S(\sigma +(-k-1)f) \to V \to \scrO_S((k-1)f)\otimes
\frak m_q \to  0.$$  Here $q=z_1+z_2-p$, where $f$ is the unique fiber
containing $Z=\{z_1, z_2\}$,
$p=\sigma \cap f$, and the addition is with respect to the group law on $f$
\rom(if $f$ is singular and $\Supp Z$ meets the singular point then $q$ is the
singular point as well\rom). If $Z\in \Sym ^2\sigma-\Cal D$, then the maximal
destabilizing sub-line bundle is $\scrO_S(\sigma -(k+2)f)$. \endstatement
\proof If $Z\notin D_\sigma$, then we have seen in (4.6)(v) that $V$ is
unstable if and only if $Z\in \Cal D$, and in this case the destabilizing
sub-line bundle must be $\scrO_S(\sigma +(-k-1)f)$ by (4.4). The quotient is
torsion free and by a Chern class calculation it must be
$\scrO_S((k-1)f)\otimes
\frak m_q$ for some point $q$. To identify the point $q$, let us assume for
simplicity that $\Supp Z$ does not meet the singular point of a singular
fiber, we can restrict the two exact sequences for $V$ to the fiber $f$
containing $Z$. From these we see that there are surjective maps $V|f \to
\scrO _f(p-z_1-z_2)$ and $V|f
\to \scrO _f(-q)$. Since $\deg V|f= 1$, it splits and the unique summand of
negative degree is thus  $\scrO _f(p-z_1-z_2) \cong \scrO _f(-q)$. It follows
that
$q=z_1+z_2-p$. The case where $\Supp Z$ contains the singular point of a
singular fiber is similar.

If $Z\in D_\sigma$, then since $V$ is locally free $Z\in \Sym ^2\sigma$.
Arguments as in (4.6) then show that $V$ is unstable. If moreover $Z\notin
\Cal D$, then by (4.4) the maximal destabilizing sub-line bundle is
$\scrO_S(\sigma -(k+2)f)$.  \endproof

Our next task will be to construct a universal sheaf $\Cal V$ over
$\operatorname{Hilb}^2S-D_\sigma$. We begin by finding a sheaf $\Cal W$ as
follows: let $\Cal Z\subset S\times \operatorname{Hilb}^2S$ be the universal
subscheme, and consider the relative extension sheaf $Ext ^1_{\pi _2}(\pi
_1^*\scrO _S(\sigma -kf)\otimes I_{\Cal Z}, \pi _1^*\scrO_S((k-2)f)$. Since
$H^1(\scrO_S(-\sigma +(2k-2)f)) = 0$, there is an exact sequence
$$\gather 0 \to Ext ^1_{\pi _2}(\pi _1^*\scrO _S(\sigma -kf)\otimes I_{\Cal
Z}, \pi _1^*\scrO_S((k-2)f) \to\\
\to R^0\pi _2{}_*Ext ^1(\pi _1^*\scrO _S(\sigma -kf)\otimes I_{\Cal Z}, \pi
_1^*\scrO_S((k-2)f).
\endgather$$ Over the complement of $\Sym ^2\sigma$, $Ext ^1_{\pi _2}(\pi
_1^*\scrO _S(\sigma -kf)\otimes I_{\Cal Z}, \pi _1^*\scrO_S((k-2)f)$ is a line
bundle on $\operatorname{Hilb}^2S - \Sym ^2\sigma$ which we denote by $\Cal
L^{-1}$ and thus there is a coherent sheaf $\Cal W$ defined by
$$0 \to \pi _1^*\scrO_S((k-2)f) \otimes \Cal L \to \Cal W \to \pi _1^*\scrO
_S(\sigma -kf)\otimes I_{\Cal Z} \to 0.$$ However, if $Z\in \Cal D$, then
$\Cal W|S\times \{Z\}$ is not stable, and if
$Z\in D_\sigma$ then $\Cal W|S\times \{Z\}$ is neither locally free nor
stable. We shall first study $\Cal W|S\times (\operatorname{Hilb}^2S
-D_\sigma)$, and shall denote this for simplicity again by $\Cal W$. There is
a unique point $q= z_1+z_2-p$ such that $\Cal W|S\times \{Z\}$ maps
surjectively to
$\scrO_S((k-1)f)\otimes \frak m_q$, and so we expect to be able to make an
elementary transformation along $\Cal D$. Indeed, since $\dim \Hom
(\scrO_S(\sigma +(-k-1)f), \Cal W|S\times \{Z\})=1$ for all $Z\in \Cal D$,
there are line bundles
$\Cal L_1$ and $\Cal L_2$ on $\Cal D$ and an exact sequence
$$0 \to \pi _1^*\scrO_S(\sigma +(-k-1)f)\otimes \pi _2^*\Cal L_1 \to   \Cal
W|S\times \Cal D \to \pi _1^*\scrO_S((k-1)f)\otimes \pi _2^*\Cal L_2 \otimes
I_{\Cal Y} \to 0,$$ where $\Cal Y$ is  the set
$$\{\,(q, z_1, z_2)\in S\times _{\Pee ^1}\Cal D\mid q=z_1+z_2-p\,\}.$$ It is
easy to check from the definition that $\Cal Y$ is smooth and that the map
$\Cal Y \to \Cal D$ is an isomorphism. Thus we may define $\Cal V$ by the
exact sequence
$$0 \to \Cal V \to \Cal W \to i_*\pi _1^*\scrO_S((k-1)f)\otimes \pi _2^*\Cal
L_2
\otimes I_{\Cal Y} \to 0,$$ where $i$ is the inclusion of $S\times \Cal D$ in
$S\times  (\operatorname{Hilb}^2S -D_\sigma)$. We then have the following:

\proposition{4.12} The sheaf $\Cal V$ is a reflexive sheaf, flat over
$\operatorname{Hilb}^2S -D_\sigma$. The restriction of $\Cal V$ to each slice
$S\times \{Z\}$ is a stable torsion free sheaf, which is locally free if and
only if $Z\notin \Cal D$.
\endstatement
\proof By (A.2) of the Appendix, $\Cal V$ is reflexive and flat over
$\operatorname{Hilb}^2S -D_\sigma$. For each $Z\in \Cal D$, if $V_Z$ is the
restriction of $\Cal V$ to the slice $S\times \{Z\}$, there is an exact
sequence
$$0 \to \scrO_S((k-1)f)\otimes \frak m_q \to V_Z \to \scrO_S(\sigma
+(-k-1)f)\to 0,$$ by (A.2) again. If $Z\notin \Cal D$ then $V_Z =\Cal
V|S\times \{Z\}$ is locally free and stable. Thus we need only check that the
double dual of $V_Z$, for
$Z\in \Cal D$, is the unique nonsplit extension of  $\scrO_S(\sigma
+(-k-1)f)$ by
$\scrO_S((k-1)f)$, which will imply that $V_Z\spcheck{}\spcheck$ is up to a
twist the stable bundle $V_0$.

To verify that the double dual of $V_Z$ is a nonsplit  extension amounts to
the following: the extension class corresponding to $V_Z$ lives in $\Ext ^1(
\scrO_S(\sigma +(-k-1)f), \scrO_S((k-1)f)\otimes \frak m_q) =
H^1(\scrO_S(-\sigma +2kf)\otimes \frak m_q)$, and we must show that its image
in $H^1(\scrO_S(-\sigma +2kf)$ is nonzero. To do this we shall use the result
(A.4) of the Appendix. Let
$M= \scrO _S(\sigma -(k+1)f)$ and $L=\scrO_S((k-1)f)$. Clearly $\Hom (M,L) =
0$. By the definition of $\Cal W$ there is an exact sequence $$0 \to \pi
_1^*\scrO _S(-\sigma +(2k-1)f)\otimes \pi _2^*\Cal L \to \Cal W \otimes \pi
_1^*M^{-1} \to
\pi _1^*\scrO _S(f)\otimes I_{\Cal Z} \to 0.$$ By (4.1) $R^1\pi _2{}_*\pi
_1^*\scrO _S(-\sigma +(2k-1)f) =  R^2\pi _2{}_*\pi _1^*\scrO _S(-\sigma
+(2k-1)f) = 0$. Thus $R^1\pi _2{}_*\Cal W
\otimes \pi _1^*M^{-1} \cong R^1\pi _2{}_*\bigl(\pi _1^*\scrO _S(f)\otimes
I_{\Cal Z}\bigr)$. To analyze $R^1\pi _2{}_*\bigl(\pi _1^*\scrO _S(f)\otimes
I_{\Cal Z}\bigr)$, use the exact sequence
$$0 \to \pi _1^*\scrO _S(f)\otimes I_{\Cal Z} \to \pi _1^*\scrO _S(f) \to
\scrO_{\Cal Z}\otimes \pi _1^*\scrO _S(f)
\to 0.$$ It is easy to check that $R^1\pi _2{}_*\pi _1^*\scrO _S(f) =0$ if
$p_g>0$, and is a line bundle if $p_g=0$. Clearly $R^1\pi _2{}_*(\scrO _{\Cal
Z} \otimes \pi _1^*\scrO_S(f)) = 0$. Thus the torsion in $R^1\pi
_2{}_*\bigl(\pi _1^*\scrO _S(f)\otimes I_{\Cal Z}\bigr)$ is the cokernel of
the map between two rank two vector bundles on $\operatorname{Hilb}^2S$
$$R^0\pi _2{}_*\pi _1^*\scrO _S(f) \to R^0\pi _2{}_*\bigl(\scrO_{\Cal
Z}\otimes \pi _1^*\scrO _S(f)\bigr).$$ Sine $\Cal D$ is a smooth divisor, by
using elementary divisors for the vector bundle map we can describe this
cokernel by describing what it looks like at the generic point. It is a
simple exercise in local coordinates to identify the determinant of the
vector bundle map with a local equation for $\Cal D$ at the generic point.
Thus the torsion in $R^1\pi _2{}_*\Cal W
\otimes \pi _1^*M^{-1}$ is a line bundle on $\Cal D$, which is identified with
the torsion in $R^1\pi _2{}_*\bigl(\pi _1^*\scrO _S(f)\otimes I_{\Cal
Z}\bigr)$. Similar statements hold via standard base change results if we
restrict to a first order neighborhood of $\Cal D$, where torsion is to be
interpreted in the sense of (A.4)(ii) of the appendix.

Next, let $Z= \{z_1, z_2\}\in \Cal D$ and let $W$ be the extension
corresponding to the restriction of $\Cal W$ to the slice $S\times \{Z\}$, we
must identify the corresponding extension class, i.e\. the image of the
one-dimensional vector space
$H^1(\scrO_S(f)\otimes I_Z)$ in $H^1( M^{-1} \otimes L\otimes \frak m_q)$ and
its further image in $H^1( M^{-1} \otimes L)$. Using the two exact sequences
$$\gather 0\to \scrO _S \to  W \otimes M^{-1} \to M^{-1}\otimes L\otimes
\frak m_q\to 0;\\ 0\to\scrO_S(-\sigma + (2k-1)f) \to W \otimes M^{-1}\to
\scrO_S(f)\otimes I_Z \to 0,
\endgather$$ we see that the composite map $\scrO _S \to \scrO_S(f)\otimes
I_Z$ is nonzero and gives the nontrivial section. Now the quotient of
$\scrO_S(f)\otimes I_Z$ by
$\scrO_S$ is $\scrO_f(-z_1-z_2)$. Thus there is an induced map
$M^{-1}\otimes L\otimes \frak m_q \to \scrO_f(-z_1-z_2)$ which must factor
through the natural map $M^{-1}\otimes L\otimes \frak m_q = \scrO_S(-\sigma +
2kf)\otimes \frak m_q \to \scrO _f(-p-q)$. (Here as usual $p =\sigma \cap
f$.) As the induced map $\scrO _f(-p-q) \to \scrO_f(-z_1-z_2)$ is nonzero, it
is an isomorphism, and we recover the fact that $q= z_1+z_2 -p$. Using the
commutativity of
$$\CD 0 @>>> \scrO_S(f)\otimes I_Z @>>> \scrO_S(f) @>>> \scrO_Z @>>> 0\\
@.@VVV @VVV @| @.\\ 0 @>>> \scrO_f(-z_1-z_2) @>>> \scrO_f @>>>\scrO_Z @>>> 0,
\endCD$$ we also see that the image of $H^1(\scrO_S(f)\otimes I_Z)$ in
$H^1(\scrO_f(-z_1-z_2)$ is the same as the image of $H^0(\scrO_Z)$ in
$H^1(\scrO_f(-z_1-z_2))$.

There is a commutative diagram
$$\CD @.  H^1(\scrO_S(-\sigma +2kf)\otimes \frak m_q) @>>>
H^1(\scrO_S(-\sigma +2kf))\\ @. @VVV   @VVV\\ H^1(\scrO_f(-z_1-z_2))
@>{\cong}>> H^1(\scrO _f(-p-q)) @>>> H^1(\scrO _f(-p)).
\endCD$$ Moreover the map $H^1(\scrO_S(-\sigma +2kf)) \to H^1(\scrO _f(-p))$
is an isomorphism. So the problem is the following one: does the image of
$H^0(\scrO_Z)$ in $H^1(\scrO_f(-z_1-z_2))$ map to zero in $H^1(\scrO _f(-p))$?
The image of $H^0(\scrO_Z)$ in $H^1(\scrO_f(-z_1-z_2))$ is  dual to the image
of
$H^0(\scrO _f)$ in $H^0(\scrO_f(z_1+z_2))$, giving a section vanishing at
$z_1$ and $z_2$. On the other hand the kernel of the map $H^1(\scrO _f(-p-q))
\to H^1(\scrO _f(-p))$ is dual to the image of the map $H^0(\scrO_f(p)) \to
H^0(\scrO_f(p+q))$, and the corresponding section of $\scrO_f(p+q)$ vanishes
at
$p$ and $q$. So the only way that this can equal the image of $H^0(\scrO_Z)$
is for
$z_1$ or $z_2$ to equal $p$, i.e\. $Z\in D_\sigma$. Conversely, if $Z\notin
D_\sigma$, then the image of the extension class in $H^1( M^{-1} \otimes L)$
is not zero. Thus the double dual of the restriction of $\Cal V$ to $S\times
\{Z\}$ is a nonsplit extension and so it is stable. \endproof

This is as far as we shall go in this section in calculating the
four-dimensional invariant. But let us sketch here how to obtain the full
formula in (4.8). We will prove a more general statement in Section 10, where
we will use (4.12).

First, to deal with the fact that $\dim
\Ext ^1(\scrO _S(\sigma -kf)\otimes I_Z, \scrO _S((k-2)f)$ jumps along $\Sym
^2\sigma$, blow up $\Sym ^2\sigma$ inside $\operatorname{Hilb}^2S$. Let the
exceptional divisor be $G$. After blowing up, we can asume that the extension
is not locally trivial along $G$. There is thus a universal extension of
torsion free sheaves $\tilde \Cal W$ over $S\times \operatorname{Bl}_{\Sym
^2\sigma}
\operatorname{Hilb}^2S$. Now make an elementary modification along $\Cal D$,
replacing unstable Type 1 extensions with $Z\in \Cal D - \Sym ^2\sigma $ with
stable Type 4 extensions. Next make an elementary modification along
$D_\sigma$, replacing  unstable Type 1 extensions with $Z\in D_\sigma$ with
Type 2 extensions; this also fixes some of the unstable Type 4 extensions.
Finally make an elementary modification along $G$ to replace the remaining
unstable extensions with Type 3 extensions. At this point every member of the
family is a stable torsion free sheaf, and the induced morphism to
$\overline{\frak M}_2$ blows $G$ back down again to $\Sym ^2\sigma$. The
morphism $\operatorname{Hilb}^2S \to
\overline{\frak M}_2$ is then an isomorphism. Keeping track of the Chern
classes gives the formula in Theorem 4.8.

Finally, we state a general conjecture:

\medskip
\noindent{\bf Conjecture 4.13.} If $S$ has a section, then the map of (3.14)
extends to an isomorphism $\operatorname{Hilb}^tS \to \overline{\frak M}_t$.
\medskip

If the conjecture is true, then the method of test surfaces used in the proof
of Lemma 9.2 can be used to show that the $\mu$-map is given by the following
formula (where we use the notation of the Introduction for divisors in
$\operatorname{Hilb}^tS$ as well):
$$\mu (\Sigma) = D_{\alpha _t} - \bigl((f\cdot \Sigma)/2\bigr)E,$$ where
$$\align
\alpha _t &= \Sigma + \bigl((-(\sigma \cdot \Sigma) +(p_g-1+t)(f\cdot
\Sigma))/2\bigr)f +(f\cdot \Sigma)\sigma \\ &= \alpha _1 + (t-1)\bigl((f\cdot
\Sigma)/2\bigr)f.
\endalign$$

\section{5. Calculation of the invariant for dimension two and no multiple
fibers.}

Our goal in this and the following three sections will be a complete
calculation of the  Donaldson polynomial invariant $\gamma_{w,p}$ in case $-p
-3\chi (\scrO_S) = 2$. In this case, the moduli space is compact of real
dimension four and complex dimension two, and may be identified with the
algebraic surface $J^{e+1}(S)$.  We shall begin with the case where $S$ has a
section $\sigma$ and $e=-2$. We have  already described how to calculate the
invariant in this case in the last section. However, we shall give another
method for doing so here, since it will serve to explain the construction in
the general case. In fact, we shall reprove (in a slightly different guise)
the formula in (4.7): $$-4\mu (\Sigma ) = (2(\sigma
\cdot \Sigma) -2p_g(f\cdot \Sigma))f -4(f\cdot \Sigma)\sigma -4\Sigma.$$

To describe the $\mu$-map, we begin by describing a universal bundle over
$S$.  Recall that every bundle $V$ with $-p_1(\ad V) -3\chi (\scrO_S) = 2$ is
obtained from the fixed bundle $V_0$ by a single allowable elementary
modification. For convenience we will look at the case where $e=-2$. Thus we
shall normalize $V_0$ to have $\det V_0\cdot f = -3$ and
$$-p_1(\ad V_0) -p = c_1(V_0)^2 - 4c_2(V_0) = 3(1+p_g(S)).$$ As $V_0$ is
well-defined up to twisting, so that we can assume that $c_1(V_0) =
-3\sigma$, if $p_g(S)$ is odd, and $c_1(V_0) =  -3\sigma +f$ if $p_g(S)$ is
even. (Here we could use the explicit description of
$V_0$ from the preceding section, or use the congruence $p\equiv 1+p_g\mod 4$
to see that these choices always give $c_1(V_0)^2 \equiv p\mod 4$.) We shall
just consider the argument in case $p_g$ is odd. Setting $c= c_2(V_0)$, we
have
$$4c - (-3\sigma )^2 = 3(1+p_g)$$ and thus
$$c = -\frac32(1+p_g).$$

If $V$ is stable, with $-p_1(\ad V) =3(1+p_g) +2$, then there is an exact
sequence $$0 \to V \to V_0 \to Q \to 0,$$ where $Q$ is a rank one torsion
free sheaf on a fiber $f$ with $\deg Q = -1$ and $\det V_0 \cdot f = -3$, and
conversely every such $V$ is stable. We need to parametrize such sheaves $Q$
as a family over $S\times S$, where the first factor should be viewed as the
surface and the second as the moduli space. To do so,  let $\pi _1$ and $\pi
_2$ be the projections of $S\times S$ to  the first and second factors,  let
$\Bbb D$ denote the diagonal inside $S\times S$ and let $D= S\times
_{\Pee^1}S$ be the fiber product. Thus $D$ is a Cartier divisor, which is not
however smooth at the images of pairs of double points. At such a point $D$
has the local equation $xy = zw$, and thus $D$ has an ordinary double point in
dimension three. The diagonal $\Bbb D$ is of course contained as a
hypersurface in
$D$, but this hypersurface fails to be Cartier at the singular points of $D$.
Let $\Cal P = I_{\Bbb D}/I_D$. In local analytic coordinates, $\Cal P$ looks
like $$(x-z, y-w)R/(xy-zw)R$$ near the double point, where $R =
\Cee\{x,y,z,w\}$.  We claim that the sheaf $\Cal P$ is flat over $S$ (the
second factor). Indeed there is an exact sequence $$0 \to I_{\Bbb D}/I_D \to
\scrO _D \to \scrO_{\Bbb D} \to 0.$$ Moreover $\scrO_{\Bbb D}$ is obviously
flat over $S$ and $\scrO _D$ is flat over $S$ since $D$ is a local complete
intersection inside $S\times S$. Thus $\Cal P$ is flat over $S$ also. Given
$q \in S$ denote $\Cal P|\pi _2^{-1}(q)$ by $\Cal P_q$, where we shall
identify $\Cal P_q$ with the corresponding torsion sheaf on
$S$. If $q$ is not a singular point of a nodal fiber, then  $\Cal P_q =
\scrO_f(-q)$, where $f$ is the fiber containing $q$ and we have identified
$\scrO_f(-q)$ with its direct image on $S$ under the inclusion. If $q$ is the
singular point of a singular fiber, then in local analytic coordinates $\Cal
P_q$ is given by $$(x-z, y-w)R/(xy-zw, z,w)R \cong
(x,y)\Cee\{x,y\}/(xy)\Cee\{x,y\}.$$ Thus globally $\Cal P_q$ is the maximal
ideal of $q$, in other words it is the unique torsion free rank one sheaf of
degree $-1$ on the singular fiber which is not locally free.

Fix as above $V_0$ to be a stable rank two vector bundle on $S$ of fiber
degree
$-3$ such that the restriction of $V_0$ to every fiber is stable. Thus as we
have seen in (1.2) and (2.7)(i), $\dim \Hom (V_0, \Cal P_q) = h^0(V_0\spcheck
\otimes
\Cal P_q) = 1$ and $h^1 (V_0\spcheck \otimes \Cal P_q) = 0$. It follows via
flat base change as in the proof of (3.15) that  $\pi _2{}_*((\pi
_1^*V_0)\spcheck
\otimes \Cal P)$ is a line bundle on $S$. We let $\Cal L$ denote the dual line
bundle. Thus $$\align
\Hom (\pi _1^*V_0, \Cal P \otimes \pi _2^*\Cal L) &= H^0(S\times S;(\pi
_1^*V_0)
\spcheck \otimes  \Cal P \otimes \pi _2^*\Cal L) \\ &= H^0(S; \pi _2{}_*((\pi
_1^*V_0)\spcheck \otimes \Cal P) \otimes \Cal L)\\ &= H^0(S; \Cal
L^{-1}\otimes \Cal L) = H^0(S; \scrO_S).
\endalign$$ Thus there is a nonzero map $\pi _1^*V_0 \to \Cal P \otimes \pi
_2^*\Cal L$,  essentially unique, and its restriction to each fiber $\pi
_2^{-1}(q)$ is also nonzero. We may then define a universal bundle by the
exact sequence
$$0 \to \Cal V \to \pi _1^*V_0 \to \Cal P \otimes \pi _2^*\Cal L \to 0.$$

\lemma{5.1} The sheaf $\Cal V$ is locally free and its restriction to each
slice
$S\times \{q\}$ is a stable rank two vector bundle $V_q$ with $-p_1(\ad V_q) -
3\chi (\scrO_S) = 2$. The resulting morphism $S\to \frak M_1$ is an
isomorphism.
\endstatement
\proof There is an exact sequence
$$0 \to V_q \to V_0 \to \Cal P_q \to 0.$$ Thus $V_q$ is locally free for all
$q$ and so is $\Cal V$. By construction $V_q$ has stable restriction to every
fiber except the one containing $q$. Thus $V_q$ is stable. The statement
about $p_1(\ad V_q)$ is clear. Finally, examining the description of (3.13),
we see that the map $S\to \frak M_1$ is a bijection. Since
$\frak M_1$ is smooth, the map is therefore an isomorphism.
\endproof

We now turn to calculating the Chern classes of $\Cal V$. By (0.1)
$$p_1(\ad\Cal V) - p_1(\ad \pi ^*V_0) = 2c_1(V_0)\cdot D + [D]^2
-4i_*c_1(\Cal P
\otimes \pi _2^*\Cal L),$$ where $i\:D \to S\times S$ is the inclusion. Here
the sheaf $\Cal P \otimes  \pi _2^*\Cal L$ fails to be a line bundle exactly
at the singular points of $D$, which does not affect the Chern classes $c_1$
and $c_2$. Thus we can simply define
$i_*c_1(\Cal P \otimes  \pi _2^*\Cal L)$ to be the unique extension of the
class
$i_*c_1(\Cal P \otimes  \pi _2^*\Cal L|D_{\text{reg}})$. Next we claim:

\lemma{5.2} In $H^2(S\times S)$, we have $[D] = f\otimes 1 + 1\otimes f$.
\endstatement
\proof Let $C$ be a Riemann surface embedded in $S$, and consider
$([C]\otimes  [x])\cup [D]$, where $x$ is a point of $S$. This is the same as
$\#\bigl((C\times
\{x\})\cap D\bigr)$, where the points are counted with signs. Clearly this
intersection is the same as $\#(C\cap f)$. A similar argument holds for
$([x]\otimes [C])\cap [D]$. Thus $[D]$ and  $f\otimes 1 + 1\otimes f$ define
the same element of $H^2(S\times S)$.  \endproof

It follows that, up to a term not affecting slant product,
$$p_1(\ad\Cal V) - p_1(\ad \pi ^*V_0) = -6\sigma \otimes f + 2f\otimes f -4i_*
c_1(\Cal P  \otimes \pi _2^*\Cal L).$$ Next we must calculate the most
interesting term in the expression for $p_1(\ad
\Cal V)$ above, the term $c_1(\Cal P  \otimes
\pi _2^*\Cal L)$, viewed as a coherent sheaf on $D$. As far as $c_1$ is
concerned,  we can ignore the singularities of $D$. Thus
$$\align c_1(\Cal P \otimes  \pi _2^*\Cal L|D_{\text{reg}}) &= c_1(I_{\Bbb
D}/I_D|D_{\text{reg}}) + \pi _2^*c_1(\Cal L)\\ &= -[\Bbb D] + \pi _2^*
c_1(\Cal L).
\endalign$$ Here $[\Bbb D]$ is viewed as a divisor on $D_{\text{reg}}$.
However the unique  extension of $i_*[\Bbb D]$ to an element of $H^4(S\times
S)$ is clearly again
$[\Bbb D]$, where we now view $\Bbb D$ as a codimension two cycle on $S\times
S$. Now let $\alpha = c_1(\Cal L^{-1}) \in H^2(S)$. Then
$$\align i_*\pi _2^*c_1(\Cal L^{-1}) &= i_*i^*(1\otimes \alpha)\\ &=
i_*i^*(1)\cup (1\otimes \alpha) = [D]\cup (1\otimes \alpha)\\ &= f\otimes
\alpha + 1\otimes [f\cdot \alpha].
\endalign$$ Thus up to a term which does not affect slant product, $i_*\pi
_2^*c_1(\Cal L^{-1}) =   f\otimes \alpha$. To calculate this term, we shall
use the following lemma:

\lemma{5.3} $\alpha = c_1(\Cal L^{-1}) = -3\sigma -\frac52(p_g+1)f$.
\endstatement
\proof We shall apply the Grothendieck-Riemann-Roch theorem to calculate the
Chern  classes of  $$\Cal L^{-1}  = \pi _2{}_*((\pi _1^*V_0)\spcheck\otimes
\Cal P).$$  We have
$$\ch((\pi _2)_!((\pi _1^*V_0)\spcheck\otimes \Cal P)\Todd S = \pi _2{}_*
\Bigl[\ch ((\pi _1^*V_0)\spcheck\otimes \Cal P))\cdot \Todd (S\times
S)\Bigl].$$ Now $H^1(V_0 \spcheck \otimes Q) = 0$ for all $Q$ a torsion free
rank one sheaf  on a fiber $f$, so that $(\pi _2)_!((\pi
_1^*V_0)\spcheck\otimes \Cal P) = \pi _2{}_*((\pi _1^*V_0)\spcheck\otimes
\Cal P) = \Cal L^{-1}$ and the left hand side above is just $c_1(\Cal
L^{-1})\Todd S$. Now we can also multiply by $(\Todd S)^{-1}$ to get $$\align
c_1(\Cal L^{-1}) &= \pi _2{}_*\Bigl[\ch ((\pi _1^*V_0)\spcheck\otimes \Cal
P))\cdot \Todd (S\times S)\Bigl]\cdot  (\Todd S)^{-1}\\ &= \pi _2{}_*\Bigl[\ch
((\pi _1^*V_0)\spcheck\otimes \Cal P))\cdot \Todd (S\times S)\cdot \pi
_2^*(\Todd S)^{-1}\Bigl]\\ &= \pi _2{}_*\Bigl[\ch ((\pi
_1^*V_0)\spcheck\otimes \Cal P))\cdot \pi _1^*\Todd S\Bigl]\\ &=\pi
_2{}_*\Bigl[\ch (\pi _1^*V_0)\spcheck\cdot \pi _1^*\Todd S\cdot \ch(\Cal
P)\Bigl],
\endalign$$ using the multiplicativity of the Todd class.   Moreover
$$\ch (V_0\spcheck) = 2-c_1(V_0) + \frac{c_1(V_0)^2 - 2c_2(V_0}2 = 2 + 3\sigma
-6(1+p_g)[\pt]$$ and
$$\Todd S = 1 - \frac{(p_g-1)}2f + (p_g+1)[\pt].$$ So
$$\pi _1^*\ch (V_0\spcheck)\cdot \pi _1^*\Todd S = 2 + 3\sigma \otimes
1-(p_g-1)f
\otimes 1 + N[\pt]\otimes 1,$$ where
$$N= \frac{(3\sigma)^2 -2c}2 +2(p_g+1) -\frac32(p_g-1) =
\frac{-5p_g +1}2,$$ using the fact that $c = -\frac32(1+p_g)$.

Next we compute $\ch\Cal P= \ch (I_\Bbb D/I_D) = \ch I_\Bbb D - \ch I_D$. Now
$I_D = \scrO_{S\times S}(-D)$, so that $\ch I_D = 1 -[D] + [D]^2/ 2 -\cdots$.
As for $\ch I_\Bbb D$, we have $\ch I_\Bbb D = 1- \ch \scrO _\Bbb D$.
Applying the Grothendieck-Riemann-Roch formula to the inclusion $j\: \Bbb D
\to S\times S$  gives $\ch \scrO _\Bbb D = j_*((\Todd N_{\Bbb D/S\times
S})^{-1})$, where
$N_{\Bbb D/S\times S}$ is the normal bundle of $\Bbb D$ in $S\times S$, and
so is equal to the tangent bundle $T_S$ on $\Bbb D$. Thus
$$\align
\ch\scrO _\Bbb D &= j_*\Bigl((1-\frac{(p_g-1)}2f +(1+p_g)[\pt])^{-1}\Bigr)\\
&=j_*(1+\frac{(p_g-1)}2f -(1+p_g)[\pt])\\ &= [\Bbb D] + \frac{(p_g-1)}2j_*f-
(1+p_g)j_*[\pt].
\endalign$$ Collecting up the terms through degree 3 (which are the only ones
which will contribute) gives
$$\ch\Cal P = [D] - \frac{[D]^2}2 - [\Bbb D] - \frac{p_g-1}2j_*f+\cdots.$$
Putting this together, we see that $\alpha$ is the degree one term in
$$\pi _2{}_*\Bigl[(2 + 3\sigma \otimes 1-(p_g-1)f\otimes 1 + N[\pt]\otimes
1)\cdot ([D] - \frac{[D]^2}2 - [\Bbb D] -
\frac{p_g-1}2j_*f)\Bigr].$$  Recalling that $D=f\otimes 1 +1\otimes f$ and
that
$[D]^2/2 = f\otimes f$, we must apply $\pi _2{}_*$ to
$$\pi _2{}_*((p_g-1)j_*f -3(\sigma \otimes 1)\cdot [\Bbb D] +
(p_g-1)(f\otimes 1)\cdot [\Bbb D] -3(\sigma \otimes 1)\cdot (f\otimes f)
+N[\pt]
\otimes f).$$  The result is then
$$-(p_g-1)f -3f -3\sigma  +(p_g-1)f +Nf=-3\sigma + (N-3)f,$$ as claimed.
\endproof

The above lemma thus implies that
$$-4i_*c_1(\Cal P
\otimes \pi _2^*\Cal L)= 4[\Bbb D] -12f\otimes \sigma -10(p_g+1)f\otimes f.$$
Putting this together gives (neglecting all terms which do not affect slant
product) $$p_1(\ad \Cal V) = -6(\sigma \otimes f) + (-10p_g-8)f\otimes f
-12f\otimes \sigma + 4[\Bbb D] +\cdots.$$ We may finally summarize our
calculations as follows:

\lemma{5.4} In the above notation,
$$-4\mu (\Sigma) = \Bigl[-6(\sigma\cdot\Sigma) +(-10p_g-8)(f\cdot
\Sigma)\Bigr]f -12(f\cdot \Sigma)\sigma + 4\Sigma.$$ Thus $\mu (\Sigma )^2 =
(\Sigma)^2 + (p_g-1)(f\cdot \Sigma)^2$. \qed
\endstatement
\medskip

At first glance, this formula looks quite different from the previous formula
$$-4\mu (\Sigma ) = (2(\sigma \cdot \Sigma) -2p_g(f\cdot \Sigma))f -4(f\cdot
\Sigma)\sigma -4\Sigma.$$ However, the surface $S$ (viewed as the moduli
space) has an involution $\iota$,  coming from taking $x\mapsto -x$ on each
fiber using $\sigma$ as the identity section. This involution corresponds to
viewing $S$ as the double cover of a rational ruled surface as in [6] Chapter
1. Since $S$ has only nodal singular fibers, it follows that on $H^2(S)$,
$\iota$ fixes $\sigma$ and $f$ and is equal to $-\Id$ on the orthogonal
complement $\{f, \sigma\}^\perp$. It is then an easy exercise to see that for
a general $\Sigma$ we have $$\iota ^*(\Sigma) = -\Sigma +  2\Bigl[(\sigma
\cdot \Sigma)+(p_g+1)(f\cdot\Sigma)\Bigr]f + 2(f\cdot
\Sigma)\sigma.$$ Applying $\iota$ then exchanges the above two expressions for
$\mu (\Sigma)$. Clearly this discrepancy arose as follows. In the general
scheme for identifying the moduli space implicit in (3.14) and (4.7) we used
not $\Cal P$ but its dual. However it was technically slightly simpler not to
make this choice in the Riemann-Roch calculation above. Thus the
identifications of the moduli space differ by $-\Id$.

\section{6. The case of multiple fibers.}

Having done the rather tedious calculation in the preceding section in case
$S$  has a section, we must now move on to deal with the case where $S$ has
multiple fibers. Fortunately, it will turn out that much of the calculation
in this case exactly follows the pattern of the previous calculation. Before
getting into the nitty-gritty, let us fix notation. Let $\pi\:S\to \Pee ^1$
be a nodal surface with at most two multiple fibers of odd multiplicity. Fix
a divisor on the generic fiber $S_\eta$ of odd degree $2e+1$. Let $V_0$ be a
rank two vector bundle on $S$ with $c_1(V_0) = \Delta$ and $c_2(V_0) = c$,
whose restriction to the reduction of every fiber is stable. Thus $4c -
\Delta ^2 = 3(p_g+1)$ and so $$2c = \frac{\Delta ^2 +3(p_g+1)}2.$$ We would
like to construct a universal bundle using
$J^{e+1}(S)$. Unfortunately, this is not in general possible, and we shall
instead use a finite cover. Thus we fix an elliptic surface $T$ together with
a map $T\to S$, such that $T$ has a section. We may further assume that $T$
is obtained as follows: choose a smooth multisection $C$ of $\pi$, for
example a general hyperplane section of $S$ in some projective embedding. For
$C$ sufficiently general, we may assume that $C$ meets the multiple fibers
transversally and that the map $C\to \Pee ^1$ is not branched at any points
corresponding to singular nonmultiple fibers of $\pi$. Then set $T$ to be the
normalization of $S\times _{\Pee ^1}C$. It follows that the only singular
fibers of $T$ lie over singular nonmultiple fibers of $S$, and that $T$ has a
section $\sigma$. If $d$ is the degree of $C\to \Pee ^1$, then at the point
of $\Pee ^1$ lying under the multiple fiber $F_i$ of multiplicity $m_i$,
$C\to \Pee ^1$ is branched to order $m_i$ at exactly $d/m_i$ points.

Let $\varphi\:T\to S$ be the natural map and $\rho\:T\to C$ be the elliptic
fibration, so that we have a commutative diagram
$$\CD T@>{\varphi}>>S\\ @V{\rho}VV  @VV{\pi}V \\ C@>>> \Pee ^1.
\endCD$$

Now we can state the main result of this section:

\theorem{6.1} There exists a vector bundle $\tilde \Cal V$ over $S\times T$
with  the following properties:
\roster
\item"{(i)}" The restriction of $\tilde \Cal V$ to each slice $S\times \{p\}$
is a  stable rank two vector bundle $V$ with $\det V = \Delta -f$ and
$-p_1(\ad V) -3\chi(\scrO _S) =2$. \item"{(ii)}" The morphism $T\to \frak
M_1$ induced by $\tilde \Cal V$ has degree $d$.
\item"{(iii)}" If $\tilde \mu\: H_2(S) \to H^2(T)$ is the map induced by
slant  product with the class $-p_1(\ad \tilde \Cal V)/4$, then, setting
$\delta =[\Delta]$,
$$\align -4\tilde \mu (\Sigma) &= \Bigl[\delta ^2 -(1+p_g) -4(e+2)^2(1+p_g)
+2 + c(e,m_1) + c(e,m_2)\Bigr] (f\cdot \Sigma)df \\ &-4(e+2)(\varphi ^*\delta
\cdot \sigma)(f\cdot \Sigma)f - 4(e+2)(\varphi ^*
\Sigma \cdot \sigma)f +2d(\delta \cdot \Sigma)f\\ & +4(f\cdot \Sigma)\varphi
^*\delta-8(e+2)(f\cdot \Sigma)\sigma +4\varphi ^*\Sigma,
\endalign$$ where $c(e,m_i)$ depends only on $m_i$ and $e$ and on an analytic
neighborhood of  the multiple fiber, and not on $S$ or $p_g$, and where
$c(e,1) =0$. \endroster
\endstatement
\medskip

We shall defer the proof of Theorem 6.1 to the next two sections. The
constant
$c(e, m_i)$ in fact might depend {\it a priori\/} on the particular choice of
the multiple fiber. However, as we shall see from Theorem 6.3, the choice of
the fiber and of $e$ does not matter. Let us begin with a calculation of $\mu
(\Sigma)^2$:

\lemma{6.2} With notation as in \rom{(6.1)}, we have
$$16\tilde \mu (\Sigma)^2 = 16d(\Sigma)^2 +16d(p_g-1-c(e,m_1)
-c(e,m_2))(f\cdot
\Sigma)^2.$$ Thus as $\tilde \mu (\Sigma)^2 = d\mu (\Sigma)^2$, we have
$$\align
\mu (\Sigma)^2 &=(\Sigma)^2 +(p_g-1-c(e,m_1) -c(e,m_2))(f\cdot \Sigma)^2\\
&=(m_1m_2)^2 (p_g-1-c(e,m_1) -c(e,m_2))(\kappa\cdot \Sigma)^2,
\endalign$$ where $\kappa$ is the primitive class such that $m_1m_2\kappa =
f$.
\endstatement
\proof This is a tedious calculation.
\endproof

\theorem{6.3} With notation as in the statement of \rom{(6.1)}, we have
$$c(e,m_i) = -1 + \frac{1}{m_i^2}.$$
\endstatement
\proof By symmetry it suffices to consider $i=1$. Choose a general nodal
rational  elliptic surface $S_0$ with a single multiple fiber of multiplicity
$m_1$. We can assume that an analytic neighborhood of the multiple fiber in
$S_0$ is analytically isomorphic to a neighborhood of $F_1$ in $S$, which is
possible since we assumed that the multiple fibers did not lie over branch
points of the
$j$-function of $S$. Since $m_1|2e+1$, there exists  a divisor $\Delta$ on
$S_0$ with $\Delta \cdot f = 2e+1$. Thus we may use $S_0$ to calculate $c(e,
m_1)$. Now setting $p_g =0$ and $m_2=1$ in the formula of (6.2) gives the
coefficient of
$(\kappa \cdot \Sigma)^2$ in the Donaldson polynomial: it is $(m_1)^2(-1
-c(e,m_1))$. On the other hand, $S_0$ is orientation-preserving diffeomorphic
to a rational elliptic surface $S_1$ with a section, by a diffeomorphism
$\psi$ which carries $\kappa$ to the class of a fiber. Using (3.5) of Part I
of [4], this diffeomorphism must then carry a $(w,p)$-suitable chamber for
$S_1$ to a
$(\psi ^*w,p)$-suitable chamber for $S_0$. The Donaldson polynomial for $S_1$
and a $(w,p)$-suitable chamber is then sent under $\psi ^*$ to the $\pm$ the
Donaldson polynomial for $S_0$ and a $(\psi ^*w,p)$-suitable chamber.
Normalizing the orientations so that the leading coefficients agree (these
are both $(\Sigma ^2)$), the coefficients of $(\kappa \cdot \Sigma)^2$ must
agree also. We have already calculated the coefficient of $(\kappa \cdot
\Sigma)^2$ for $S_1$ (by two different methods): it is $-1$. Thus
$$(m_1)^2(-1 -c(e,m_1)) = -1.$$ Hence
$c(e,m_1) = -1 + 1/m_1^2$, as claimed. \endproof

Thus we get the formula for $\mu (\Sigma)^2$ stated in (ii) of Theorem 2 of
the Introduction:

\corollary{6.4} The two-dimensional Donaldson polynomial is given by the
formula
$$\align
\mu (\Sigma)^2 &= (\Sigma)^2 +(m_1m_2)^2 (p_g-1 + 1-\frac{1}{m_1^2} +
1-\frac{1}{m_2^2}) (\kappa\cdot \Sigma)^2\\ &= (\Sigma)^2 +\bigl[(m_1m_2)^2
(p_g+1) - m_1^2 -m_2^2\bigr](\kappa\cdot
\Sigma)^2. \qed \endalign$$
\endstatement

For future reference we note the following lemma:

\lemma{6.5} If $f$ denotes the general fiber of $\frak M_1= J^{e+1}(S) \to
\Pee^1$, then  $$\mu (\Sigma)\cdot f = 2(f\cdot \Sigma).$$
\endstatement
\noindent {\it Proof.} It suffices to calculate $\tilde \mu (\Sigma )\cdot
f$,  where $\tilde \mu$ is as defined in (6.1)(iii) and here $f$ denotes a
general fiber on $T$. But using the formula in (6.1)(iii) gives
$$\tilde \mu (\Sigma )\cdot f = -(f\cdot \Sigma)(2e+1) +2(e+2) (f\cdot
\Sigma) -  (f\cdot \Sigma) = 2(f\cdot \Sigma). \qed$$

\section{7. Proof of Theorem 6.1: a Riemann-Roch calculation.}

We return to the notation of the preceding section. Our goal in this section
will  be  to approximate the universal bundle by a coherent sheaf which is
essentially an elementary modification of $\pi _1^*V_0$, where $V_0$ is as
described at the beginning of the preceding section and $\pi _i$ denotes the
$i^{\text{th}}$ projection now on $S\times T$. We have the map $\varphi \:
T\to S$  of elliptic surfaces covering the map $\rho \:C\to \Pee ^1$ of the
base curves. Let $\Gamma$ be the graph of $\varphi$ in $S\times T$ and let
$H$ be the graph of the composition $\psi\:
T@>{\rho}>>C@>{\sigma}>>T@>{\varphi}>>S$, where we view
$\sigma$ temporarily not as a curve in $T$ but rather as a morphism. Let $D =
S\times _{\Pee^1}T\subset S\times T$ and let $\tilde D$ be the normalization
of
$D$. The singularities of $D$ are of two types. The first type consists of
points
$(p,q)$ where $\varphi(q) = p$ and $p$ and $q$ are the singular points on a
nodal fiber. At such points $D$ has an ordinary double point as in the case
where $S$ has a section. The second type of singularity is along a multiple
fiber $F_i$. At a point of $\Pee ^1$ lying under $F_i$, the map $C\to \Pee
^1$ is branched to order $m_i$. Thus, in local analytic coordinates $x,y,z,w$
on $S\times T$ the divisor $D$ has the local equation   $x^{m_i} = z^{m_i}$.
If $R$ is the local ring of $D$ at such a point and $\tilde R$ is its
normalization, then the inclusion
$R\subseteq \tilde R$ is given by $$\Cee\{x,y,z,w\}/(x^{m_i} - z^{m_i})
\hookrightarrow \bigoplus _k\Cee\{x,y,w\},$$ where the map from $R$ to the
$k^{\text{th}}$ factor in the direct sum is given by setting $z = \zeta ^kx$
for
$\zeta = e^{2\pi \sqrt{-1}/m_i}$.

Let $\tilde F_i = \varphi ^{-1}(F_i)$ and let $E_i$ be a component of $\tilde
F_i$. There is thus an induced map $\nu _i \: E_i \to F_i$ which is \'etale
of degree
$m_i$. We also have maps $D\to T$ and $\tilde D\to D$. Clearly $D$ and
$\tilde D$ are flat over $T$ (note that $\tilde D$ is smooth away from the
images of pairs of double points). The calculations above for $R$ and $\tilde
R$ show that the scheme-theoretic fiber of $D$ at a point $q\in E_i$ is $F_i$
as a multiple fiber and that $i_*\scrO _{\tilde D}$ restricted to this fiber
is $\nu _i{}_*\scrO _{E_i}$.

Since a section cannot pass through a singular point of a fiber, the graph
$H$  avoids the double point singularities of $D$. Denote also by $H$ the
pullback of
$H$ to $\tilde D$. Then $H$ is a Cartier divisor on $\tilde D$. Define
$$\Cal P = i_*\scrO _{\tilde D}(-\Gamma + (e+2)H).$$ This notation does not
define $\Cal P$ near the double points of $D$, but as $H$  does not pass
through the double points and $\tilde D=D$ in a neighborhood of the double
points we can just glue $\Cal P$ to $I_\Gamma/I_D$ at the double points.
Equivalently we could just take the push-forward of the restriction of
$i_*\scrO _{\tilde D}(-\Gamma + (e+2)H)$ to $D_{\text{reg}}$. Finally we
shall let $\pi _1$ and $\pi _2$ denote the first and second projections on
$S\times T$.

\lemma{7.1} The sheaf $\pi _2{}_*((\pi _1^*V_0)\spcheck \otimes \Cal P)$ is a
line  bundle on $T$, whose dual is denoted $\Cal L$. Moreover
$$R^1\pi _2{}_*((\pi _1^*V_0)\spcheck \otimes \Cal P) = 0.$$
\endstatement
\proof Letting $h\: \tilde D \to T$ and $j\: \tilde D \to S\times S@>{\pi
_1}>> S$  be the natural maps, it is clear that
$$\pi _2{}_*((\pi _1^*V_0)\spcheck \otimes \Cal P) = h_*((j^*V_0)\spcheck
\otimes
\scrO _{\tilde D}(-\Gamma + (e+2)H)).$$ So we must check that the restriction
of $(j^*V_0)\spcheck \otimes
\scrO _{\tilde D}(-\Gamma + (e+2)H)$ to each fiber of $h$ has $h^0 =1$ and
$h^1=0$. The only new case is the case corresponding to a multiple fiber. In
this case the  restriction to the fiber is $(\nu _i^*V_0)\spcheck \otimes L$,
where $L$ is a line bundle of degree $e+1$ on $E_i$. The degree of $\nu
_i^*V_0$ is $m_i(\deg V/m_i) = 2e+1$ and $\nu _i^*V_0$ is stable since it is
the pullback of the stable bundle
$V_0|F_i$. Thus by (1.2),  $H^0(E_i;(\nu _i^*V_0)\spcheck \otimes L)$ has
dimension one and  $H^1(E_i;(\nu _i^*V_0)\spcheck \otimes L) = 0$.
\endproof

Thus arguing as in the case of a section there is a unique nonzero map  (mod
scalars)
$$\pi _1^*V_0 \to \Cal P \otimes \pi _2^*\Cal L.$$ Unfortunately, if there
are multiple fibers this map is no longer surjective. We  shall return to
this point in the next section. Our remaining goal in this section is to
calculate $\Cal L$:

\lemma{7.2} With $\Cal L^{-1}= \pi _2{}_*((\pi _1^*V_0)\spcheck \otimes \Cal
P)$  and $\delta =[\Delta]$, we have $$c_1(\Cal L^{-1}) =\Big[\frac{\delta
^2}4-\frac{1+p_g}4 -(e+2)^2(1+p_g) \Big]df  -(e+2)(\varphi ^*\delta \cdot
\sigma)f + \varphi ^*\delta -2(e+2)\sigma.$$ \endstatement
\proof As before we shall apply the Grothendieck-Riemann-Roch theorem to find
$c_1(\pi _2{}_*((\pi _1^*V_0)\spcheck \otimes \Cal P))$: it is the degree one
term in
$$\pi _2{}_*(\pi _1^*\ch V_0\spcheck \cdot\pi _1^*\Todd S \cdot \ch i_*\scrO
_{\tilde D}(-\Gamma + (e+2)H)).$$ We have
$$\ch V_0\spcheck = 2-\delta +\fracwithdelims(){\delta ^2-2c}2[\pt],$$ where
$\delta =[\Delta]$, and
$$\Todd S = 1 + \frac{r}2f + (1+p_g)[\pt],$$ where
$$-r = (p_g+1) - \frac{1}{m_1} - \frac{1}{m_2}.$$  Thus the product of the
first two terms above is
$\pi _1^*(2-\delta + rf +M[\pt])$, where
$$M= \frac{\delta ^2-2c}2 + 2(1+p_g) -\frac{r}2(2e+1).$$ Since we have
$$\delta ^2-2c = \frac{\delta ^2 -4c}2 + \frac{\delta ^2}2 = -\frac32(1+p_g)
+ \frac{\delta ^2}2,$$ we can rewrite this as
$$M= \frac{\delta ^2}4 + \frac54(1+p_g) -\frac{r}2(2e+1).$$ Next we must
calculate $\ch i_*\scrO _{\tilde D}(-\Gamma + (e+2)H))$. Again using the
Grothendieck-Riemann-Roch theorem, and setting  $G= -\Gamma + (e+2)H$ for
notational simplicity, we have
$$\ch i_*\scrO _{\tilde D}(G) = i_* \bigl[\ch \scrO _{\tilde D}(G)\cdot(\Todd
N_i)^{-1}\bigr],$$ where $N_i$ is the normal bundle to the immersion $i$. Now
$\ch \scrO _{\tilde D}(G) = 1 + G + G^2/2 + \cdots$. As for $N_i$, locally at
the multiple fiber $F_i$ $D$ is the union of $m_i$ sheets, and so
$$N_i = \scrO _{\tilde D}(D - (m_1-1) B_1 - (m_2-1)B_2),$$ where $B_i = F_i
\times \tilde F_i$.  It follows that
$$(\Todd N_i)^{-1} = 1-\frac{D - (m_1-1) B_1 - (m_2-1)B_2}2 +\cdots$$ and so
$$\align
\ch i_*\scrO _{\tilde D}(G) = D + G &-i_*\fracwithdelims(){D - (m_1-1) B_1 -
(m_2-1)B_2}2 + i_*\fracwithdelims(){G^2}2 \\ &-i_*\fracwithdelims(){G\cdot (D
- (m_1-1) B_1 - (m_2-1)B_2)}2 +\cdots .
\endalign$$ So we must take the degree three term in the product of the above
expression  with  $\pi _1^*(2-\delta + rf +M[\pt])$ and then apply $\pi
_2{}_*$. First, a  calculation along the lines of (5.2) shows that
$$[D] = f\otimes 1 + d(1\otimes f),$$ where $f$ denotes either the class of a
fiber in $S$ or $T$, depending on the  context. The degree three term above
is then a sum of three terms: $T_1+T_2+T_3$, where $$\align T_1 &=
M([\pt]\otimes 1)\cdot D \\ T_2 &= -G\cdot (\delta \otimes 1) + G\cdot
(rf\otimes 1) \\ &-\frac12i_*( D - (m_1-1) B_1 - (m_2-1)B_2)\cdot (-\delta
\otimes 1 + rf\otimes 1)\\ T_3 &= i_*(G^2 -G\cdot i^*D + (m_1-1)(G\cdot B_1)
+ (m_2-1)(G\cdot B_2)).
\endalign$$ Let us now apply $\pi _2{}_*$ to these terms. First
$$\pi _2{}_* T_1= \pi _2{}_*(Md)[\pt]\otimes f= (Md)f.$$ To calculate $\pi
_2{}_*T_2$, first note the following, whose proof is an easy verification:

\lemma{7.3} For every $\alpha \in H^2(S)$,
$$\align
\pi _2{}_*(\Gamma \cdot\alpha \otimes 1) &= \varphi ^*\alpha;\\
\pi _2{}_*(H\cdot \alpha \otimes 1) &= (\varphi ^*\alpha \cdot \sigma)f. \qed
\endalign$$
\endstatement
\medskip

So the terms involving $G$ in $\pi _2{}_*T_2$ give
$$\align  -(e+2)(\varphi ^*\delta \cdot \sigma)f &+ \varphi ^*\delta +
(e+2)r(\varphi ^*f
\cdot \sigma)f- r\varphi ^*f\\ =-(e+2)(\varphi ^*\delta \cdot \sigma)f &+
\varphi ^*\delta + (e+1)rdf,
\endalign$$ where we have used $\varphi  ^*f = df$.

To handle the terms involving $B_i$, note that $i_*[B_i] = m_i[F_i\times
\tilde F_i]$. Also $[F_i] = (1/m_i)f$ and $\tilde F_i$ consists of $d/m_i$
copies of $f$ (the fiber on $T$) so that
$$[F_i\times \tilde F_i] = \fracwithdelims(){d}{m_i^2}f\otimes f; \qquad
i_*[B_i]  =  \frac{d}{m_i}f\otimes f.$$ Also $i_*D = D^2 = 2d(f\otimes f)$.
Thus
$$-\frac12 i_*(D - (m_1-1) B_1 - (m_2-1)B_2) = -\frac{d}2\left(\frac1{m_1} +
\frac1{m_2}\right)(f\otimes f).$$ The product of this term with
$f\otimes 1$ is zero, and we are left with the product with $-\delta \otimes
1$,  which contributes
$$ \frac{d(2e+1)}2\left(\frac1{m_1} +
\frac1{m_2}\right)f.$$ Combining these, we see that
$$\pi _2{}_*T_2 =-(e+2)(\varphi ^*\delta \cdot \sigma)f + \varphi ^*\delta +
(e+1)rdf+ \frac{d(2e+1)}2\left(\frac1{m_1} +
\frac1{m_2}\right)f.$$ We turn now to the term $\pi _2{}_*T_3$. We have $G^2
= (e+2)^2H^2 -2(e+2)H \cdot
\Gamma + \Gamma ^2$. To calculate $\pi _2{}_*$ applied to these terms, we
shall use the following lemma:

\lemma{7.4}
\roster
\item"{(i)}"  $\pi _2{}_*i_*H^2 = \pi _2{}_*i_*\Gamma ^2 = -d(1+p_g)f$.
\item"{(ii)}"  $\pi _2{}_*i_*H\cdot \Gamma = \sigma$.
\endroster
\endstatement
\proof To see (i), note that we have an exact sequence
$$0 \to N_{\Gamma /\tilde D} \to N_{\Gamma /S\times T} \to N_i \to 0.$$ Also
$(\Gamma ^2)_{\tilde D} = \phi_*c_1(N_{\Gamma /\tilde D})$, where $\phi\:
\Gamma \to \tilde D$ is the inclusion. Now
$$\align c_1(N_{\Gamma /\tilde D}) &= c_1(N_{\Gamma /S\times T}) - c_1(N_i)\\
&= c_1(\pi _1^*T_S|\Gamma) -(D-(m_1-1)B_1-(m_2-1)B_2)|\Gamma\\ &= \varphi
^*(rf) - ((f\otimes 1 +d(1\otimes f)-(m_1-1)(B_1\cdot
\Gamma)-(m_2-1)(B_2\cdot\Gamma))\\ &= \Bigl[-d\left(p_g+1-\frac1{m_1}
-\frac1{m_2}\right)-\left(2d -\frac{d(m_1-1)} {m_1} -
\frac{d(m_2-1)}{m_2}\right)\Big]f\\ &= -d(p_g+1)f.
\endalign$$ Thus $\pi _2{}_*i_*\Gamma ^2 = -d(1+p_g)f$.  A similar
calculation handles $\pi _2{}_*i_*H^2$. The proof of (ii) is an easy
calculation.
\endproof

Thus
$$\pi _2{}_*G^2 = -d(p_g+1)((e+2)^2 +1)f - 2(e+2)\sigma.$$ The remaining term
is $-\pi _2{}_*(G\cdot(D-(m_1-1)B_1-(m_2-1)B_2))$. We have  seen in the
course of the proof of Lemma 7.4 that
$$\align
\pi _2{}_*\Gamma \cdot (D-(m_1-1)B_1-(m_2-1)B_2) &=
\pi _2{}_*H\cdot (D-(m_1-1)B_1-(m_2-1)B_2)\\ &= \left(\frac1{m_1} +
\frac1{m_2}\right)df.
\endalign$$ Thus
$$\pi _2{}_*G\cdot (D-(m_1-1)B_1-(m_2-1)B_2)= d(e+1)\left(\frac1{m_1} +
\frac1{m_2}\right)f.$$ In all then,
$$\pi _2{}_*T_3 = d\Big[ -(p_g+1)((e+2)^2 +1)- (e+1)\left(\frac1{m_1} +
\frac1{m_2}\right)\Big]f - 2(e+2)\sigma.$$

Combining terms, we have
$$c_1(\Cal L^{-1}) =\Big[\frac{\delta ^2}4-\frac{1+p_g}4 -(e+2)^2(1+p_g)
\Big]df  -(e+2)(\varphi ^*\delta \cdot \sigma)f + \varphi ^*\delta
-2(e+2)\sigma,$$ as claimed. This concludes the proof of Lemma 7.2. \endproof

\section{8. Proof of Theorem 6.1: Conclusion.}

We keep the notation of the two previous sections. We begin by constructing a
``universal bundle" $\tilde \Cal V$ over $S\times T$.  Begin with the morphism
$\pi _1^*V_0 \to \Cal P \otimes \pi _2^*\Cal L$ defined in the previous
section, and let $\tilde \Cal V$ be the kernel. By construction $\tilde \Cal
V$ is locally free away from $(F_1\times \tilde F_1)\amalg (F_2\times \tilde
F_2)$. There is an  exact sequence:  $$ 0 \to \tilde \Cal V \to \pi _1^*V_0
\to \Cal P \otimes \pi _2^*\Cal L \to  \Cal Q_1\oplus \Cal Q_2 \to 0 .$$
where $\Cal Q_i$ is supported on $F_i\times \tilde F_i$. Now
$\tilde F_i$ is a disjoint union of $d/m_i$ fibers of $T$. Let $\Cal Q = \Cal
Q_1 \oplus \Cal Q _2$ and let $c$ denote the total Chern polynomial. Then
$$c(\tilde \Cal V) = \pi _1^*c(V_0) \cdot c(\Cal P \otimes \pi _2^*\Cal
L)^{-1} \cdot c(\Cal Q).$$ Thus if we let $\pi _1^*c(V_0) \cdot c(\Cal P
\otimes \pi _2^*\Cal L)^{-1} = 1 +  x_1 + x_2+ \cdots$, then
$$\align c_2(\tilde \Cal V) &= x _2 +  c_2(\Cal Q);\\ c_1(\tilde \Cal V)^2
-4c_2(\tilde \Cal V) &= x_1 ^2 -4x _2 -4 c_2(\Cal Q).
\endalign$$ Now we claim that Theorem 6.1 is a consequence of the following
two results:

\theorem{8.1} There exist integers $q(e, m_i)$ such that
$$c_2(\Cal Q_i)= dq(e, m_i)[F_i\times f].$$ Here the integer $q(e, m_i)$
depends only on an analytic neighborhood of
$F_i$ and $e$ but not on $S$ or $p_g(S)$.
\endstatement
\medskip

\theorem{8.2} The coherent sheaf $\tilde \Cal V$ is locally free.
\endstatement
\medskip

\demo{Proof that \rom{(8.1)} and \rom{(8.2)} imply Theorem \rom{(6.1)}} Let
us  consider the restriction of $\tilde \Cal V$ to a slice $S\times \{q\}$.
In all cases this restriction is a vector bundle $V$ whose restriction to
every smooth fiber $f$ of  $S$ not equal to the fiber containing $\varphi(q)$
is $V_0|f$. Thus the restriction of $V$ to such a fiber $f$ is stable, and so
$V$ is stable by (3.4). Now if $\varphi(q)$ does not lie on a multiple fiber,
there is an exact sequence  $$0 \to V \to V_0 \to Q \to 0,$$ where $Q$ is the
direct image of the line bundle on $f$ corresponding to the  divisor
$(e+2)\psi(q) - \varphi(q)$, which has degree $e+1$. Thus $c_1(V) = \Delta
-f$ and $p_1(\ad V) = p_1(\ad V_0) -2$. This establishes (i) of Theorem 6.1.
Note also that the map $q\mapsto (e+2)\psi(q) - \varphi(q)$ defines a
rational map from $T$ to $J^{e+1}(S)$ (which in fact is a morphism) and the
map
$T \to \frak M_1$ factors through the map $T\to J^{e+1}(S)$, compatibly with
the identification of a dense open subset of $J^{e+1}(S)$ with a dense open
subset of
$\frak M_1$ given in (3.14).

Next let us calculate the degree of the induced morphism $T\to \frak M_1$.
Fix a general smooth fiber $f$ of $S$, a line bundle $L$ on $f$ of degree
$e+1$ and a  vector bundle $V$ which is uniquely specified by  an exact
sequence
$$0 \to V \to V_0 \to i_*L \to 0,$$ where $i\: f \to S$ is the inclusion. We
shall count the preimage of $V$ in $T$.  If $f$ is general, then $T \to S$ is
unbranched over $f$ and the preimage of $f$ consists of $d$ distinct fibers
$f_1, \dots, f_d$. Moreover $\varphi$ restricts to an isomorphism from $f_i$
to $f$ for each $i$. The image of $f_i$ under $\psi$ is a single point $p_i
\in f$ corresponding to the point $\sigma \cap f_i$. Now clearly there is a
unique point $q_i \in f_i$ such that  $$L = \scrO _f((e+2)p_i -
\varphi (q_i)).$$ Thus the preimage of $V$ consists of $d$ distinct points,
and so the map $T\to  \frak M_1$ has degree $d$.

Lastly we must calculate $p_1(\ad \tilde \Cal V)$. We begin by calculating
$\pi _1^*c(V_0) \cdot c(\Cal P \otimes \pi _2^*\Cal L)^{-1}$. Here
$\pi _1^*c(V_0) = 1 + \pi _1^*\delta + \pi _1^*c[\pt]$. As for the term
$c(\Cal P \otimes \pi _2^*\Cal L)$, we clearly have $c_1(\Cal P \otimes \pi
_2^*
\Cal L) = D = (f\otimes 1) + d(1\otimes f)$. On the other hand, with the
notation of Section 7 we may apply the Grothendieck-Riemann-Roch theorem to
the immersion
$i\: \tilde D \to S\times T$ to obtain
$$\ch (\Cal P \otimes \pi _2^*\Cal L) = i_*\bigl[\ch \scrO _{\tilde
D}((e+2)H-\Gamma)(\Todd N_i)^{-1}\cdot\pi _2^*\ch \Cal L)\bigr] .$$ A
calculation similar to those in Section 7 shows that this is equal to
$$i_*\Bigl[1 + (e+2)H - \Gamma - \pi _2^*\alpha - \frac{D-(m_1-1)B_1 -
(m_2-1)B_2}2+\cdots \Bigr],$$ where $\alpha = c_1(\Cal L^{-1})$ has been
calculated in Lemma 7.2, and further manipulation gives $$-2c_2(\Cal P\otimes
\pi _2^*\Cal L) = -2[D]^2 +  2\Bigl[(e+2)H - \Gamma - \pi _2^*\alpha \cdot
[D] -d\left(1-\frac1{m_1} + 1-\frac1{m_2}\right)(f\otimes f)\Bigr].$$
Recalling that $\pi _1^*c(V_0) \cdot c(\Cal P \otimes \pi _2^*\Cal L)^{-1} = 1
+ x_1 + x_2+ \cdots$, we have
$$(1+ x_1 + x_2+\cdots )(1+ [D] + c_2(\Cal P \otimes \pi _2^*\Cal L)) =  1 +
\pi _1^*\delta + \pi _1^*c[\pt].$$ Thus $x_1 = \pi _1^*\delta - [D]$ and
$$x_2 = \pi _1^*c[\pt] - \pi _1^*\delta \cdot [D] +[D]^2 -c_2 (\Cal P
\otimes  \pi _2^*\Cal L).$$ A calculation then shows that
$$\align x_1^2 -4x_2 &= \pi _1^* p_1(\ad V_0)+2 \pi _1^*\delta \cdot [D]
+[D]^2 -4(e+2)[H] + 4[\Gamma] + 4\pi _2^*\alpha
\cdot [D] \\ &+ 4\left(1-\frac1{m_1} + 1-\frac1{m_2}\right)d(f\otimes f).
\endalign$$ There are correction terms $b(m_i) = 1-1/m_i$ depending on the
multiple fibers. Now
$$\align p_1(\ad \tilde \Cal V) &=x_1^2-4x_2 -4c_2(\Cal Q)\\ &= \pi _1^*
p_1(\ad V_0)+2 \pi _1^*\delta \cdot [D] +[D]^2 -4(e+2)[H] + 4[\Gamma] + 4\pi
_2^*\alpha
\cdot [D] \\ &+ 4(b(m_1) - q(e, m_1)/m_1  + b(m_2) - q(e, m_2)/m_2
)d(f\otimes f),
\endalign$$ where the terms $b(m_i)$, $q(e, m_i)$ depend only on an analytic
neighborhood of the multiple fiber and are both 0 if $m_i =1$. Let $c(e, m_i)
=  4(b(m_i) - q(e, m_i)/m_i)$. Taking  slant product of this expression with
$[\Sigma]$, using  the fact that $[\Gamma ]
\backslash [\Sigma] = \varphi ^*\Sigma$ and $[H] \backslash [\Sigma] =
(\varphi ^*\Sigma\cdot \sigma)f$, and plugging in the expression for $\alpha$
given by Lemma 7.2 gives the final formula in Theorem 6.1(iii). \endproof

\demo{Proof of \rom{(8.1)}} Choose an analytic neighborhood $X$ of $F_i$. We
may  assume that $X$ fibers over the unit disk in $\Cee$. Then $\varphi
^{-1}(X)$ consists of $d/m_i$ copies of $\tilde X$, which is the
normalization of the pullback of $X$ by the map from the disk to itself
defined by $z = w^{m_i}$. Restrict $\varphi$ and $V_0$ to this local
situation, and let $D$ now denote the fiber product inside $X\times \tilde X$
and $\tilde D$ its normalization. We can similarly define the codimension two
subsets $\Gamma$ and $H$. Let us examine the dependence of the terms $V_0$
and $\scrO _{\tilde D}((e+2)H-\Gamma)$ on the various choices.

First, suppose that $V_0$ and $V_0'$ are two different choices of a bundle
over
$X$ whose determinants have fiber degree $2e+1$ and whose restrictions to the
reduction of every fiber are stable. Then $\det V_0 \otimes (\det V_0')^{-1}$
has fiber degree zero. On the other hand, from the exponential sheaf sequence
$$H^1(X; \scrO _X) \to \Pic X \to H^2(X; \Zee)$$ and the identification
$H^2(X; \Zee) \cong H^2(F_i;\Zee ) \cong \Zee$, it follows  that the group of
line bundles of fiber degree zero is divisible. Thus there is a line bundle
$L$  on $X$  such that $\det V_0' = \det (V_0\otimes L)$. The proof of
Corollary 3.8 shows that $V_0'$ and $V_0\otimes L$ differ by twisting by a
line bundle pulled back from the  disk, which is necessarily trivial. Thus
$V_0' \cong V_0\otimes L$.

The remaining choice was the choice of a section $\sigma$ of $\tilde X$.
Given two  such choices $\sigma _1$ and $\sigma_2$, we have two divisors
$H_1$ and $H_2$ on
$\tilde D$, and two line bundles $\scrO _{\tilde D}((e+2)H_1-\Gamma)$ and
$\scrO _{\tilde D}((e+2)H_2-\Gamma)$. Their difference is the line bundle
$\scrO  _{\tilde D}((e+2)(H_1-H_2)$. The restriction of $\scrO _{\tilde
D}((e+2)H_i-\Gamma)$ to each fiber $f$ of the map $\tilde D \to \tilde X$ over
$q\in \tilde X$ is the line bundle $\scrO _f((e+2)p_i - q)$, where $p_i =
\sigma _i \cap f$ and we can identify the fiber over $q$ with the fiber on
$\tilde X$ containing $q$ via $\varphi$. Let $\Psi \: \tilde X \to \tilde X$
be the inverse of the map given by translation by the divisor of fiber degree
zero
$(e+2)(\sigma _1-\sigma _2) - c_1(\tilde L)$, where $\tilde L$ is the
pullback to
$\tilde X$ of $L$. Thus $\Psi ^{-1}(q) = q + (e+2)(p_1-p_2)-\lambda$, where
$q\in f$ and $\lambda$ is the line bundle $L|f$. Now
$\operatorname{Id}\times \Psi$ acts on $X \times \tilde X$, preserving the
divisor
$D$ and acting as well on the normalization $\tilde D$. Clearly the line
bundles
$\scrO _{\tilde D}((e+2)H_2-\Gamma)\otimes \pi _1^*L$ and
$(\operatorname{Id}\times \Psi)^*\scrO _{\tilde D}((e+2)H_1-\Gamma)$ have
isomorphic restrictions to each fiber of the map $\tilde D \to \tilde X$. Thus
they differ by the pullback of a line bundle $L'$ on $\tilde X$.  Thus we
have an isomorphism
$$ (\Id \times \Psi )^*(\pi _1^*V_0\spcheck\otimes i_*\scrO _{\tilde
D}((e+2)H_1-\Gamma)) \cong (\pi _1^*V_0)\spcheck \otimes i_*\scrO _{\tilde D}
((e+2)H_2-\Gamma)\otimes \pi _1^*L\otimes \pi _2^*L' $$ and a similar
isomorphism when we apply $R^0\pi _2{}_*$. Lastly every map $\pi ^*V_0 \to
i_*\scrO _{\tilde D}((e+2)H_1-\Gamma)$ which corresponds to an everywhere
generating section of the line bundle $\pi _2{}_*Hom (\pi _1^*V_0, i_*\scrO
_{\tilde D}((e+2)H_1-\Gamma)$ under the natural map
$$\pi _2^*\pi _2{}_*Hom (\pi _1^*V_0, i_*\scrO _{\tilde D}((e+2)H_1-\Gamma)
\to Hom (\pi _1^*V_0, i_*\scrO _{\tilde D}((e+2)H_1-\Gamma)$$  is determined
up to multiplication by a nowhere vanishing function on $\tilde X$. It now
follows that, up to twisting by the pullback of the line bundle $L'$ on
$\tilde X$, we may identify the map $\pi _1^*V_0'\to  i_*\scrO _{\tilde
D}((e+2)H_2-\Gamma)$, up to a nowhere vanishing function on $\tilde X$ and up
to twisting by the pullback of a line bundle on $\tilde X$,  with the
pullback under  $(\operatorname{Id}\times \Psi)^*$ of the corresponding map
from
$\pi _1^*V_0$ to $i_*\scrO _{\tilde D}((e+2)H_1-\Gamma)$. In particular the
cokernels of these maps, viewed as sheaves supported on $F_i \times E_i$,
have the same length. But the lengths of the cokernels  are exactly what is
needed to calculate $c_2(\Cal Q_i)$, in the notation of  the beginning of
this section. Thus we have established (8.1). \endproof

\noindent {\bf Remark.} We could easily show directly by a slight
modification of  the proof above that the integers $q(e, m_i)$ defined above
are independent of
$e$. \medskip

\demo{Proof of \rom{(8.2)}} We begin with the following (see also (A.2)(i)):

\lemma{8.3} The sheaf $\tilde \Cal V$ is reflexive.
\endstatement
\proof Since $\tilde \Cal V$ is a subsheaf of the locally free sheaf  $\pi
_1^*V_0$, it is torsion free. Thus it will suffice to show that every section
$\tau$ of  $\tilde \Cal V$ defined on an open set of the form $W-Z$, where
$W$ is an open subset of $S\times T$ and $Z$ is a  closed subvariety of $W$
of codimension at least two, extends to a section of
$\tilde \Cal V$ over $W$. Now locally (after possibly shrinking $W$) $\tilde
\Cal V$ is given by an exact sequence $$0 \to \tilde \Cal V|W \to \scrO _W^2
\to i_*\scrO _{\tilde D}|W .$$ Now viewing the section $\tau$ as a section of
$\scrO _W^2$ over $W-Z$, it extends as a section of $\scrO _W^2$  by Hartogs'
theorem. Let $\tilde \tau$ be the unique extension. Then the image of
$\tilde \tau$ in $i_*\scrO _{\tilde D}|W$ vanishes on $D-Z$, which is
nonempty. Clearly then it is zero. Thus the extension $\tilde \tau$ defines a
section of
$\tilde \Cal V$ extending $\tau$, so that $\tilde \Cal V$ is reflexive.
\endproof

Returning to the proof of (8.2), let $U = S\times T -(F_1\times
\tilde F_1) - (F_2 \times \tilde F_2)$. By Lemma 8.3,  $\tilde
\Cal V$ is a reflexive sheaf which is locally free on $U$. We claim that
$\tilde \Cal V$ is everywhere locally free. The problem is local around each
point
$(x,y)$ of $F_i\times \tilde F_i$. Since $\tilde \Cal V$ is reflexive, it
will suffice to show the following:  each point $y$ of
$\tilde F_i$ has a neighborhood $\Cal N$  such that $\tilde \Cal V|(S\times
\Cal N)\cap U$ has an extension to a locally free sheaf over $S\times \Cal N$.

Let $T_0 = T- \tilde F_1-\tilde F_2$. Clearly $T_0$ is the inverse image of
$J^{e+1}(S)-F_1-F_2$ under the natural morphism from $T$ to $J^{e+1}(S)$. The
restriction of $\Cal V_U$ to  $S\times T_0$ is a bundle over $S\times T_0$ in
the sense of schemes since it is the restriction of a coherent sheaf over
$S\times T$. Thus it induces a morphism of schemes from $T_0$ to $\frak M_1$.
If we denote the points of $\frak M_1$ corresponding to multiple fibers by
$F_1$ and $F_2$ again, then it is easy to see that the map of (3.14) extends
to an embedding
$J^{e+1}(S)-F_1-F_2 \to \frak M_1-(F_1\cup F_2)$. Thus the map $T_0 \to \frak
M_1-(F_1\cup F_2)$ is proper. This map extends to a rational map from $T$ to
$\frak M_1$. After blowing up $T$, there is a morphism from the blowup
$\tilde T$ to $\frak M_1$. The image of $\tilde T-T_0$ must clearly lie
inside the two elliptic curves in $\frak M_1$ corresponding to elementary
modifications along
$F_1$ or $F_2$. Since there are no nonconstant maps from $\Pee ^1$ to an
elliptic curve,  every exceptional curve on $\tilde T$ is mapped to a point,
and the map
$T_0 \to \frak M_1$ extends to a morphism $\Phi \:T\to \frak M_1$. Clearly the
morphism $\Phi \:T\to \frak M_1$ identifies $\frak M_1$ with $J^{e+1}(S)$.

Given $y\in \tilde F_i$, choose a neighborhood $N_0$ of $\Phi(y)$ in $\frak
M_1$  such that there exists a universal bundle over $S\times N_0$, and let
$\Cal N$ be the component of $\Phi ^{-1}(N_0)$ containing $y$. Thus there is
a universal vector bundle $\Cal W$ over $S\times \Cal N$. By construction
$\Cal W|S\times (\Cal N -\tilde F_i)$ and $\tilde\Cal V|S\times (\Cal N
-\tilde F_i) $ have isomorphic restrictions to every slice $S\times \{z\}$.
Thus $\pi _2{}_*Hom (\Cal W, \tilde \Cal V)$ is a torsion free rank one sheaf
on $\Cal N$, which is thus an ideal sheaf on $\Cal N$ if $\Cal N$ is small
enough. We may assume that  $\pi _2{}_*Hom (\Cal W, \tilde \Cal V)|\Cal
N-\{y\}$ is just the structure sheaf. Choosing an everywhere generating
section of $\pi _2{}_*Hom (\Cal W, \tilde \Cal V)|\Cal N-\{y\}$ gives a
homomorphism $\Cal W|S\times (\Cal N-\{y\}) \to \tilde
\Cal V|S\times (\Cal N-\{y\})$. This homomorphism is an isomorphism over
$S\times (\Cal N -\tilde F_i)$ and is nonzero at a general point of $S\times
((\Cal N -\{y\})\cap \tilde F_i)$. As both $\Cal W$ and $\tilde \Cal V$ are
vector bundles away from $F_i\times (\Cal N \cap \tilde F_i)$ whose
restrictions to every smooth fiber of $S$ in every slice are stable, it
follows that $\Cal W|S\times (\Cal N-\{y\}) \to \tilde \Cal V|S\times (\Cal
N-\{y\})$ is an isomorphism in codimension one. Since both sheaves are
reflexive, they are isomorphic. Finally
$\Cal W$ and $\tilde \Cal V$ are two reflexive sheaves which are isomorphic
on the complement of the codimension two set $S\times \{y\}\subset S\times
\Cal N$, so they are isomorphic. Thus $\tilde \Cal V$ is locally free.
\endproof

\section{9. The four-dimensional invariant.}

Our goal in this section will be to calculate the four-dimensional
invariant.   What follows is an outline of the calculation. Let
$\overline{\frak M}_2$ denote the moduli space of Gieseker stable torsion
free sheaves on $S$ of dimension four. As we have seen, $\overline{\frak
M}_2$ is smooth and irreducible and birational  to
$\operatorname{Hilb}^2J^{e+1}(S)$. In fact, we shall begin by establishing a
more precise statement. Let $Y_i\subset \operatorname{Hilb}^2J^{e+1}(S)$ be
the subset of codimension two consisting of subschemes of $J^{e+1}(S)$ whose
support has reduction contained in the multiple fiber $F_i$ on $J^{e+1}(S)$.
Clearly $Y_i$ has two components: one component is just
$\operatorname{Sym}^2F_i$, the closure of the locus of two distinct points
lying on $F_i$, and the other is a $\Pee ^1$-bundle over $F_i$ corresponding
to nonreduced subschemes whose support is a point on $F_i$. There is a
similar subscheme $Y_i'$ of $\overline{\frak M}_2$, consisting of torsion
free sheaves $V$ on $S$ such that either $V$ is not locally  free and the
unique point where $V$ is not locally free lies on $F_i$ or $V$ is a  bundle
obtained from $V_0$ up to equivalence by taking two elementary modifications
along line bundles on $F_i$. We claim:

\lemma{9.1} The isomorphism defined in \rom{(3.14)} from a Zariski open
subset of
$\Sym ^2J^{e+1}(S)$ to an open subset of $\frak M_2$ extends to an isomorphism
$\operatorname{Hilb}^2J^{e+1}(S) -  Y_1-Y_2 \to \overline{\frak M}_2 -
Y_1'-Y_2'$.
\endstatement \medskip

Let us remark that, in case there are multiple fibers, the birational map
above does not extend to a morphism. This follows from the identification of
the function $d(e, m_i)$ below, and can also be seen directly as follows. The
moduli space $\overline{\frak M}_2$ contains the set of nonlocally free
sheaves, which is a smooth $\Pee^1$-bundle over $S$. The corresponding subset
of
$\operatorname{Hilb}^2J^{e+1}(S)$ is the image of the blowup of
$J^{e+1}(S)\times _{\Pee ^1}J^{e+1}(S)$ along the diagonal (which is not a
Cartier divisor) under the involution. It is easy to see that this image is
not normal along the image of $F_i\times F_i$ if $m_i>1$.

There is an isomorphism $H^2(\operatorname{Hilb}^2J^{e+1}(S) -  Y_1-Y_2 )
\cong H^2(\overline{\frak M}_2 - Y_1'-Y_2')$, so that by restriction we can
view $\mu (\Sigma)$ as an element of  $H^2(\operatorname{Hilb}^2J^{e+1}(S) -
Y_1-Y_2 ) \cong H^2(\operatorname{Hilb}^2J^{e+1}(S))$. Denote this element of
$H^2(\operatorname{Hilb}^2J^{e+1}(S))$ by $\mu '(\Sigma)$. In fact, it is
easy to  identify this element: let $\alpha _1 = \mu _1(\Sigma ) \in
J^{e+1}(S)$ be given  by the $\mu$-map for the two-dimensional invariant, and
set
$$\alpha _2 = \alpha _1 + \frac{(f\cdot \Sigma)}2f.$$ Then we have the
following formula:

\lemma{9.2}
$$\mu '(\Sigma) = D_{\alpha _2} -\frac{(f\cdot \Sigma)}2E.$$
\endstatement

Now  $\alpha _1^2$ is just the value of the two-dimensional invariant, which
we  shall write as $(\Sigma ^2) + C_1(\kappa\cdot \Sigma)^2$, where $C_1 =
m_1^2m_2^2(p_g+1) -m_1^2-m_2^2$. Thus
$$\align
\alpha _2^2 &= \alpha _1 ^2 + (f\cdot \Sigma)(\alpha _1\cdot f) \\ &= \alpha
_1 ^2 + 2(f\cdot \Sigma)^2 \\ &= (\Sigma ^2) + (C_1+2m_1^2m_2^2)(\kappa\cdot
\Sigma)^2,
\endalign$$ where we have used Lemma 6.5 to conclude that $\alpha _1\cdot f =
2(f\cdot
\Sigma)$.

Thus a routine calculation with the multiplication table in
$\operatorname{Hilb}^2J^{e+1}(S)$ gives:

\lemma{9.3}
$$\align
\mu '(\Sigma)^4 =& 3(\Sigma ^2)^2 + 6C_1(\Sigma ^2)(\kappa \cdot \Sigma)^2 +
\\ &+ \Bigl[ 3C_1^2-(2(p_g+1) +12)m_1^4m_2^4 + 8(m_1^3m_2^4 +
m_1^4m_2^3)\Bigr]  (\kappa \cdot \Sigma)^4.\qed
\endalign$$
\endstatement

Of course, this is a calculation on  $\operatorname{Hilb}^2J^{e+1}(S)$, not
on
$\overline{\frak M}_2$. To get an answer on $\overline{\frak M}_2$, we shall
argue  that the above formula must be corrected by terms which only depend on
the multiplicities of the multiple fibers and not on $p_g$.

\lemma{9.4} There exist a function  $d(e, m_i)$, depending only on $e$ and an
analytic neighborhood of the multiple fiber $F_i$ in $S$, with the  following
properties:
\roster
\item"{(i)}" $d(e, 1)=0$.
\item"{(ii)}" $\mu (\Sigma )^4 - \mu '(\Sigma )^4 = m_1^4m_2^4(d(e,m_1) +d(e,
m_2))(\kappa \cdot \Sigma )^4.$
\endroster
\endstatement
\medskip

We can now complete the proof of (iii) of Theorem 2 in the Introduction. It
follows from (9.3) and (9.4) that the coefficient of $(\kappa \cdot
\Sigma)^4$ in
$\mu (\Sigma )^4$  is given by
$$3C_1^2-(2(p_g+1) +12)m_1^4m_2^4 + 8(m_1^3m_2^4 + m_1^4m_2^3) +
m_1^4m_2^4(d(e,m_1) +d(e, m_2)).$$ To calculate  $d(e, m_i)$, take as before
$S$ to be a rational surface with a  multiple fiber of multipicity $m_1$. In
this case, arguing as in the proof of (6.3), the coefficient of  $(\kappa
\cdot \Sigma )^4$ is the same as the coefficient of $(\kappa \cdot
\Sigma )^4$ for the rational surface with no multiple fibers. To calculate
this coefficient, we apply (9.4) and (9.3) with $m_1=m_2 =1$ and $p_g=0$, to
see that
$\mu (\Sigma )^4 = \mu' (\Sigma )^4$ and thus that the coefficient of $(\kappa
\cdot \Sigma )^4$ is $3-14 +16 =5$. Now taking $p_g =0$ and $m_1$ arbitrary
and
$m_2=1$ in the above formulas gives $C_1=-1$ and  $$5 = 3(-1)^2 -14m_1^4
+8(m_1^3+m_1^4) +m_1^4d(e, m_1).$$  Thus $m_1^4d(e, m_1) = 2-8m_1^3 +6m_1^4$,
or $$d(e, m_1) = \frac{2}{m_1^4}-\frac{8}{m_1} +6.$$ Plugging this into the
expression above for the coefficient for $(\kappa \cdot \Sigma)^4$ in the
general case gives $$\align 3C_1^2-(2(p_g+1)& +12)m_1^4m_2^4 +12m_1^4m_2^4
+2m_1^4 + 2m_2^4\\ = 3C_1^2-&2((p_g+1)m_1^4m_2^4 -m_1^4 - m_2^4).
\endalign$$ We may write this answer more neatly as $3C_1^2-2C_2$, where
$$\align C_1 &= m_1^2m_2^2(p_g+1)-m_1^2-m_2^2;\\ C_2 &=
m_1^4m_2^4(p_g+1)-m_1^4-m_2^4.\qed
\endalign$$

\section{10. Proof of Lemmas 9.1, 9.2, and 9.4.}

In this section we shall give a proof of the remaining results from the
previous  section.

\demo{Proof of Lemma \rom{9.1}} The lemma asserts the existence of an
isomorphism
 from  $\operatorname{Hilb}^2J^{e+1}(S) - Y_1-Y_2$ to
$\overline{\frak M}_2 - Y_1'-Y_2'$ extending the isomorphism given in (3.14).
The isomorphism of (3.14) is defined on the open set $U$ of
$\operatorname{Hilb}^2J^{e+1}(S)$ consisting of pairs of points $\{z_1, z_2\}$
such that $z_1$ and $z_2$ lie in distinct fibers, neither of which is
singular or multiple. We must show that the map extends over the set of pairs
$\{z_1, z_2\}$, where $z_1$ and $z_2$ lie in distinct fibers, one or both of
which may be singular or multiple, as well as over the set of pairs $Z$ where
either $Z$ is nonreduced but the support of $Z$ does not lie in a multiple
fiber or where $Z =\{z_1, z_2\}$ with $z_1$ and $z_2$ lying in the same
nonmultiple fiber.

Let us first consider the case where $z_1$ and $z_2$ lie in distinct fibers.
As  in Section 7, choose an elliptic surface $T \to C$  with a section such
that $C$ is a finite cover of $\Pee ^1$, generically branched except below
the multiple fibers and $T$ is the normalization of $S\times _{\Pee ^1}C$.
Let $\varphi \: T\to S$ be the natural map. There is also the map $\varphi
_{e+1}\: T\to J^{e+1}(S)$ defined by $\Cal P$, i.e\. if $q\in T$, $f$ is the
fiber containing
$q$ and $p = f\cap \sigma$, then $\varphi _{e+1}(q) = \scrO_f((e+2)p-q)$.  We
have constructed a universal bundle $\tilde \Cal V \to S\times T$ in Section
7 for the choices of $w$ and $p$ corresponding to the two-dimensional
invariant. Let $\tilde U\subset T\times T$ be the open set of pairs of points
$(y_1, y_2)$ such that
$\varphi (y_1)$ and $\varphi (y_2)$ lie in different fibers. Let $\tilde \Cal
V_1$ be the pullback of $\tilde \Cal V$ to $S\times \tilde U$ via the natural
projection of $S\times \tilde U\subset S\times T\times T$ onto the first and
second factors. We also have the coherent sheaf $\Cal P$ on $S\times T$
defined at the beginning of Section 7. Let $\Cal P'$ be the pullback of $\Cal
P$ to $S\times
\tilde U$ defined by the projection of $S\times T\times T$ to the first and
third factor. Thus given a point $(y_1, y_2) \in \tilde U$, the restriction
of $\tilde
\Cal V_1$ to the slice through $(y_1, y_2)$ is an elementary modification of
$V_0$ along the fiber containing $\varphi (y_1)$ and the  restriction of
$\Cal P'$ to the slice through $(y_1, y_2)$ is the direct image of a line
bundle of degree
$e+1$ on the fiber through $\varphi (y_2)$. Thus, leting $\pi _2$ denote the
projection $S\times \tilde U \to \tilde U$, $\pi _2{}_*\bigl(\tilde \Cal
V_1\spcheck \otimes \Cal P'\bigr)$ is a line bundle on $\tilde U$, whose
inverse we denote by $\Cal L'$. Define $\tilde \Cal V_2$ as the kernel of the
natural map
$\tilde \Cal V_1 \to \Cal P' \otimes \Cal L'$. The proof of (8.2) shows that
$\tilde \Cal V_2$ is a vector bundle whose restriction to each slice $S\times
\{(y_1, y_2)\}$ is stable. The induced map $\tilde U \to \frak M_2$ then
descends to a map from the open subset of $\operatorname{Hilb}^2J^{e+1}(S)$
consisting of points lying in distinct fibers to $\frak M_2$. (In fact, the
proof shows that this morphism extends to a morphism defined on the
complement of the divisor $E$ of nonreduced points together with the proper
transforms of $\Sym ^2F_1$ and $\Sym ^2F_2$.)

Next we must extend the morphism over the points of
$\operatorname{Hilb}^2J^{e+1}(S)$ corresponding to points lying in the same
nonmultiple fiber and nonreduced points whose support does not lie in a
multiple fiber. In order to do so, we will need the model for elliptic
surfaces with a section constructed in Section 4. Let $Z$ be a point of
$\operatorname{Hilb}^2J^{e+1}(S)$ such that $\Supp Z$ lies in a single
nonmultiple fiber $f$, and let $X$ be a small neighborhood of $f$ mapping
properly to a disk inside $\Pee ^1$. Thus there is a biholomorphic map from
$X$ to a neighborhood of the corresponding fiber in the Jacobian surface
$J(S)$, and we may further assume that the image of $\Supp Z$ does not meet
the identity section $\sigma$ under this map. Now the results of Section 4,
after tensoring by $\scrO_X(e\sigma)$, give a rank two vector bundle $V_0'$
over $X$ whose restriction to every fiber is stable of degree $2e+1$ and a
rank two reflexive sheaf $\Cal V_0$ over $X\times
(\operatorname{Hilb}^2X-D_\sigma)$, flat over  $H =
\operatorname{Hilb}^2X-D_\sigma$, whose restriction to each slice is an
elementary modification of $V_0'$. Let $V_0$ denote as usual the bundle on
$S$ whose restriction to every fiber is stable. Then as in the proof of (8.1)
there is a line bundle $L$ on $X$ such that $V_0|X \cong V_0'\otimes L$.

The sheaf $\Cal V_0 \otimes \pi _1^*L$ has the following property. Let $\Cal B
\subset X\times H$ be the set
$$\Cal B =\{\,(x, z_1, z_2)\mid \pi (x) = \pi (z_i) \text{\, for some
$i$}\,\}.$$ Let $p$ be a point of $H$ and $\Bbb U$ a small neighborhood of
$p$, which we can identify with a neighborhood of $Z\in
\operatorname{Hilb}^2J^{e+1}(S)$. We can assume that $\Bbb U$ is a polydisk.
There is a  proper map $\Pi \: (X\times \Bbb U) - \Cal B
\to (D_0 \times \Bbb U)-\Cal B'$ induced by $\pi \: X\to D_0$, where $D_0$ is
the disk which is the base curve of $X$ and $$\Cal B' = \{\,(t, z_1, z_2)\mid
t = \pi (z_i) \text{\, for some $i$}\,\}.$$  By construction the restrictions
of $\Cal V_0 \otimes \pi _1^*L$ and $\pi _1^*V_0$ to each fiber of $\Pi$ are
isomorphic stable bundles on the fiber, which is reduced (possibly nodal).
Thus $R^0\Pi _*Hom (\Cal V_0 \otimes \pi _1^*L, \pi _1^*V_0)$ is a line
bundle $\Cal F$ on $(D_0 \times \Bbb U)-\Cal B'$. Both $\Cal V_0 \otimes \pi
_1^*L$ and $\pi _1^*V_0$ extend to coherent sheaves on $X\times \Bbb U$.
Therefore
$R^0\Pi _*Hom (\Cal V_0 \otimes \pi _1^*L, \pi _1^*V_0) =\Cal F$ extends to a
coherent sheaf on $D_0 \times \Bbb U$, which we shall continue to denote by
$\Cal F$. Replacing $\Cal F$ by its double dual if necessary, we can assume
that it is reflexive, and thus since its rank is one that it is a line
bundle. Since by assumption every line bundle on $D_0\times \Bbb U$ is
trivial, $R^0\Pi _* Hom (\Cal V_0 \otimes \pi _1^*L, \pi _1^*V_0)$ is a
trivial line bundle on
$(D_0 \times \Bbb U)-\Cal B'$, and we can thus choose an everywhere generating
section. This section corresponds to a homomorphism from $\Cal V_0 \otimes L$
to
$\pi _1^*V_0$ over $(X\times \Bbb U) - \Cal B$ which is an isomorphism on
every fiber. It follows that we can glue $\Cal V_0 \otimes L$ to $\pi
_1^*V_0$ over
$(X\times \Bbb U) - \Cal B$. Since $\{X\times \Bbb U, (S\times \Bbb U)-\Cal
B\}$ is an open cover of $S\times \Bbb U$ whose intersection is $(X\times
\Bbb U) -
\Cal B$, we have constructed an coherent sheaf on $S\times \Bbb U$, flat over
$\Bbb U$. In this way we have extended the morphism from $\tilde U \cap \Bbb
U$ over all of $\Bbb U$. So the morphism $U\to \overline{\frak M}_2$ extends
over all the points $Z\in \operatorname{Hilb}^2J^{e+1}(S)$ such that
$A\notin Y_1\cup Y_2$. Clearly its image is exactly $\overline{\frak M}_2
-Y_1'-Y_2'$.
\endproof

\demo{Proof of Lemma \rom{9.2}} We shall show that the divisor $\mu '(\Sigma)$
which is  the natural extension of the restriction of $\mu (\Sigma)$ to
$\operatorname{Hilb}^2J^{e+1}(S) - Y_1-Y_2$ to a divisor on
$\operatorname{Hilb}^2J^{e+1}(S)$ is equal to $D_{\alpha _2} -\bigl((f\cdot
\Sigma)/2\bigr)E$. Recall that $H^2(\operatorname{Hilb}^2J^{e+1}(S)) \cong
H^2(J^{e+1}(S))\oplus \Zee \cdot [E/2]$. Also, given a point $y\in
J^{e+1}(S)$, there is an induced morphism $\tau _y\:
\operatorname{Bl}_yJ^{e+1}(S) \to
\operatorname{Hilb}^2J^{e+1}(S)$  defined on $J^{e+1}(S) -\{y\}$ by $\tau
_y(x) =
\{x,y\}$. If $E_y$ is the exceptional divisor on
$\operatorname{Bl}_yJ^{e+1}(S)$, then it is easy to see that $\tau
_y^*D_\alpha = \alpha$ for all $\alpha \in H^2(J^{e+1}(S))$ (where we have
identified $H_2(J^{e+1}(S))$ and
$H^2(J^{e+1}(S))$ and identified $H^2(J^{e+1}(S))$ with a subspace of
$H^2(\operatorname{Bl}_yJ^{e+1}(S))$). Also $\tau _y^*[E] = 2[E_y]$, which
can  easily be checked by going up to the double cover of
$\operatorname{Hilb}^2J^{e+1}(S)$ which is the blowup of $J^{e+1}(S) \times
J^{e+1}(S)$ along the diagonal. Similarly, suppose that $\varphi\:T\to S$ is
a finite cover as usual and consider the morphism $\varphi _{e+1}\: T\to
J^{e+1}(S)$ defined by $\Cal P$, i.e\. if $q\in T$, $f$ is the fiber
containing
$q$ and $p = f\cap \sigma$, then $\varphi _{e+1}(q) = \scrO_f((e+2)p-q)$.
Suppose that $y$ is a general point of $J^{e+1}(S)$ (and so does not lie on a
multiple or singular fiber) and let $\varphi _{e+1}^{-1}(y) = \{y_1, \dots ,
y_d\}$. Then there is an induced map $\tau \: T-\{y_1, \dots, y_d\}\to
\operatorname{Hilb}^2J^{e+1}(S)$, and clearly we have $\tau ^*D_\alpha =
\varphi ^*\alpha$. In  particular the map
$\tau ^*$ is injective on the subspace $H_2(J^{e+1}(S))$, and we can determine
$\mu '(\Sigma)$ provided that we know $\tau ^*\mu '(\Sigma)$ and $\mu
'(\Sigma)|E_y$.  Note finally that the image of $\tau$ and $E_y$ are
contained in
$\operatorname{Hilb}^2J^{e+1}(S) - Y_1-Y_2$, so that we can calculate the
$\mu$-map by  finding a universal family of coherent sheaves on $S\times
(T-\{y_1,
\dots, y_d\})$ and over $S\times E_y$.

To find such a family, begin with the bundle $\tilde \Cal V$ over $S\times
T$. We know that $\tilde \mu (\Sigma) = \varphi _{e+1}^*\alpha _1$, where
$\tilde \mu$ is the natural $\mu$-map defined on $T$ and $\alpha _1 = \mu
(\Sigma)$ is the
$\mu$-map for the two-dimensional invariant. Fix a general fiber $f$ of $S$
and a point $y\in J^{e+1}(S)$ corresponding to a line bundle $\lambda$ of
degree $e+1$ on $f$. Let $f_1, \dots, f_d$ be the fibers on $T$ lying above
$f$ and $y_1,
\dots, y _d$ the points of $T$ corresponding to $\lambda$. We shall perform an
elementary modification along the divisor $f\times T$ with respect to the line
bundle $\pi _1^*\lambda$. This will run into trouble along $y_1, \dots, y_d$,
so that we will restrict to $S\times (T-\{y_1, \dots, y_d\})$. The upshot
will be a family of stable torsion free sheaves on $S\times (T-\{y_1, \dots,
y_d\})$ such that the induced morphism $T-\{y_1, \dots, y_d\} \to \overline
{\frak M}_2$ is
$\tau$.

First let us calculate $\Hom (\tilde \Cal V|S\times (T-\{y_1, \dots, y_d\}),
\pi _1^*\lambda)$. If $V_t$ is the restriction of $\tilde \Cal V$ to the slice
$S\times \{t\}$, then $V_t$ is an elementary modification of $V_0$ either at a
fiber different from $\lambda$ or along $f$ with respect to a line bundle
$\lambda '$ of degree equal to $\deg \lambda$ but with $\lambda '\neq
\lambda$. It follows that the map $\Hom (V_0, \lambda) \to \Hom (V_t,
\lambda)$ defined by the inclusion $V_t\subset V_0$ is a map between two
one-dimensional spaces by (1.3)(i), and its kernel is $H^0((\lambda
')^{-1}\otimes \lambda)=0$. Thus
$\Hom (V_0, \lambda) \cong \Hom (V_t, \lambda)$ and the induced map $R^0\pi
_2{}_*\pi _1^*(V_0\spcheck \otimes \lambda) \to R^0\pi _2{}_*\bigl(\tilde
\Cal V|S\times (T-\{y_1, \dots, y_d\})\spcheck\otimes \pi _1^*\lambda\bigr)$
is an isomorphism. As $R^0\pi _2{}_*\pi _1^*(V_0\spcheck \otimes \lambda)$ is
the trivial line bundle, there is a unique homomorphism mod scalars from
$\tilde \Cal V|S\times (T-\{y_1, \dots, y_d\})$ to $\pi _1^*\lambda$ and its
restriction to each slice is the corresponding nonzero homomorphism on the
slice. Let $\tilde \Cal V_2$ be the kernel, so that there is an exact sequence
$$0 \to \tilde \Cal V_2 \to \tilde \Cal V|S\times (T-\{y_1, \dots, y_d\})\to
\pi _1^*\lambda. $$ Note that the right arrow fails to be surjective over the
slice $S\times \{t\}$ only if $\varphi (t) \in f$, and in this case it
vanishes at one point. Thus by (A.5) $\tilde \Cal V_2$ is reflexive and flat
over $T-\{y_1, \dots, y_d\}$, and is a family of torsion free sheaves
parametrized by $T-\{y_1, \dots, y_d\}$. The restriction of $\tilde \Cal V_2$
to a general fiber in every slice is stable, and thus $\tilde \Cal V_2$ is a
flat family of stable torsion free sheaves.  Clearly the corresponding
morphism to $\overline {\frak M}_2$ is $\tau$.

Next we claim that
$$p_1(\ad \tilde \Cal V_2) = p_1(\ad \tilde \Cal V) -2d(f\otimes f)
+\cdots,$$  where the omitted terms do not affect slant product. Indeed the
defining map $\tilde \Cal V|S\times (T-\{y_1, \dots, y_d\})\to \pi
_1^*\lambda$ is surjective in codimension two, so that in calculating
$p_1(\ad \tilde \Cal V_2)$ we can in fact apply the formula (0.1) as if the
map were surjective. Now (0.1) gives
$$p_1(\ad \tilde \Cal V_2) = p_1(\ad \tilde \Cal V) +2c_1(\tilde \Cal V)
\cdot (f\otimes 1) - 4i_*(\lambda \otimes 1 )\cdot (f\otimes 1).$$  Using
$c_1(\tilde \Cal V) = \pi _1^*c_1(V_0) -\bigl[(f\otimes 1) +d(1\otimes
f)\bigr]$ and plugging in gives the claimed formula for $p_1(\ad \tilde \Cal
V_2)$. Thus
$$\align -(p_1(\ad \tilde \Cal V_2)\backslash \Sigma)/4 &= -(p_1(\ad \tilde
\Cal V)\backslash \Sigma)/4 +d(f\cdot \Sigma)f/2\\ &= \varphi _{e+1}^*(\alpha
_1) + \bigl((f\cdot \Sigma)/2\bigr)\varphi _{e+1}^*(f)\\ &= \varphi
_{e+1}^*(\alpha _2),
\endalign$$ and the pullback of $\mu '(\Sigma)$ to $T-\{y_1, \dots, y_d\}$
under $\tau$ is just $\varphi _{e+1}^*(\alpha _2)$. It follows that
$\mu '(\Sigma) = D_{\alpha _2} + aE$ for some rational number $a$.

To determine the coefficient of $E$ in $\mu '(\Sigma)$, fix a general fiber
$f$ of $S$ and a line bundle $\lambda$ of degree $e+1$ on $f$, which
corresponds to a point $y\in J^{e+1}(S)$. The set of points of
$\operatorname{Hilb}^2J^{e+1}(S)$ whose support is $\{y\}$ is a curve
$E_y\cong
\Pee ^1$. We shall construct a universal sheaf $\Cal V_2$ over $S\times E_y$
and show that $-(p_1(\ad \Cal V_2)\backslash \Sigma)/4 = (f\cdot \Sigma)$.

Begin with $V$ which is obtained from $V_0$ by a single elementary
modification along $\lambda$. Thus $V|f = \lambda \oplus \mu$ with $\deg
\lambda = e+1$ and
$\deg \mu =e$. By (1.3)(ii) $\dim \Hom(V, \lambda) = 2$ and there is a unique
nonzero homomorphism from $V$ to $\lambda$ which is not surjective, indeed
which vanishes exactly at the point corresponding to the degree one line
bundle
$\lambda \otimes \mu^{-1}$. Identify $\Pee(\Hom(V, \lambda))$ with $\Pee ^1$
(and with $E_y$). There is a general construction [5] of a universal
homomorphism
$\Phi\: \pi _1^*V\otimes \pi _2^*\scrO_{\Pee ^1}(-1) \to \pi _1^*\lambda$.
Thus we can define $\Cal V_2$ to be its kernel:
$$0 \to \Cal V_2 \to \pi _1^*V\otimes \pi _2^*\scrO_{\Pee ^1}(-1) \to \pi
_1^*\lambda.$$ By (A.5) in the appendix, $\Cal V_2$ is reflexive and flat
over $\Pee ^1$, and is a family of torsion free sheaves, which are locally
free except for the point of
$\Pee ^1$ corresponding to the non-surjective homomorphism. The restriction of
$\Cal V_2$ to a general fiber in every slice is stable, so that the
restriction of $\Cal V_2$ to each slice is a stable torsion free sheaf. The
induced map to
$\overline{\frak M}_2$ is easily seen to be one-to-one with image $E_y$. We
may again calculate $p_1(\ad \Cal V_2)$ by the formula of (0.1), noting that
$c_1(
\pi _1^*V\otimes \pi _2^*\scrO_{\Pee ^1}(-1) = \pi _1^*c_1(V) + 2
\pi _2^*c_1(\scrO_{\Pee ^1}(-1))$:
$$p_1(\ad \Cal V_2)=
\pi _1^*p_1(\ad V)+ 2(\pi _1^*c_1(V) + 2\pi _2^*c_1(\scrO_{\Pee ^1}(-1)))\cdot
(f\otimes 1)-4i_*\pi _1^*c_1(\lambda).$$  The only term which matters for
slant product is the term $4\pi _2^*c_1(\scrO_{\Pee ^1}(-1))\cdot (f\otimes
1)$. Thus
$$\mu '(\Sigma)\cdot E_y = -(1/4)4(-1)(f\cdot \Sigma) = (f\cdot \Sigma).$$
Bearing in mind that $E\cdot E_y = -2$, it follows that the coefficient of
$E$ in
$\mu '(\Sigma)$ is $-(f\cdot \Sigma)/2$.  So putting this all together gives
the final answer for $\mu '(\Sigma)$ in (9.2). \endproof

\demo{Proof of \rom{(9.4)}} The basic idea of the proof is similar to the idea
of the proof of (9.1). Fix an analytic neighborhood $X$ of the multiple fiber
$F_i$ as usual. Let $X_e$ be the corresponding subset of $J^{e+1}(S)$. Then
$\Bbb X = \operatorname{Hilb}^2X_e$ may be identified with an analytic open
subset of  $\operatorname{Hilb}^2J^{e+1}(S)$ which is a neighborhood of
$Y_i$. Under the birational map  $\operatorname{Hilb}^2J^{e+1}(S) \dasharrow
\overline {\frak M}_2$, the open set $\Bbb X$ corresponds birationally to an
open set $\Bbb X'$ which is a neighborhood of $Y_i'$. Moreover $\Bbb X - Y_i
\cong \Bbb X' - Y_i'$.

Now let $S_0$ be another nodal elliptic surface containing a multiple fiber of
multiplicity $m_i$ and let $X_0$ be an analytic neighborhood of the multiple
fiber. Let $\Delta _0$ be a divisor on $S_0$ of fiber degree $2e+1$ and let
$V_0'$ be a rank two vector bundle whose restriction to every fiber is
stable. We suppose that $X_0$ is biholomorphic to $X$ and identify them. We
may then define
$\Bbb X_0$ and $\Bbb X_0'$ analogously. There are also  closed subsets of
$\Bbb X_0$ and $\Bbb X_0$ corresponding to $Y_i$ and $Y_i'$, which we shall
again denote by $Y_i$ and $Y_i'$. Of course $\Bbb X_0 \cong \Bbb X$ under the
identification
$X_0 \cong X$. The main claim is then the following:

\claim{} There is a biholomorphic map $\Bbb X_0' \cong \Bbb X'$ which is
compatible with the isomorphisms
$$\Bbb X_0' - Y_i' \cong \Bbb X _0 -Y_i \cong \Bbb X-Y_i \cong \Bbb X'
-Y_i'.$$
\endstatement
\proof For emphasis, we will write $\overline{\frak M}_2(S)$ for the moduli
space for $S$, and similarly $\overline{\frak M}_2(S_0)$ for the moduli space
for
$S_0$. We shall glue $\Bbb X_0'$ to $\overline{\frak M}_2(S)-Y_i'$ along $\Bbb
X-Y_i$, and show that the result maps to $\overline{\frak M}_2(S)$,
compatibly with the inclusion $\overline{\frak M}_2(S)-Y_i' \subseteq
\overline{\frak M}_2(S)$. This will define a proper morphism from $\Bbb X_0'$
to
$\Bbb X'$ of degree one which is an isomorphism in codimension one, and thus
is an isomorphism by Zariski's Main Theorem.

We must show that the inclusion $\overline{\frak M}_2(S)-Y_i' \subseteq
\overline{\frak M}_2(S)$ extends to a morphism from $\Bbb X_0'$ to
$\overline{\frak M}_2(S)$. It suffices to do so  locally around each point of
$\Bbb X_0'$. Given an arbitrary point $p\in \Bbb X_0'$, let $\Bbb U\subset
\Bbb X_0'$ be an open neighborhood of $p$ which is biholomorphic to a
polydisk, so that $\Pic \Bbb U = 0$, and such that there exists a universal
sheaf
$\Cal V_{\Bbb U}$ over $S_0\times \Bbb U$. Denote again the restriction of
$\Cal V_{\Bbb U}$ to $X_0 \times \Bbb U = X\times \Bbb U$ by  $\Cal V_{\Bbb
U}$. Letting as usual $V_0$ denote the rank two bundle on $S$ whose
restriction to  every fiber is stable, we have seen that there is a line
bundle $L$ on $X$ such that
$V_0'\otimes L\cong V_0$. Now view $\Bbb U - Y_i'$ as an open subset of $\Bbb
X-Y_i\subset \operatorname{Hilb} ^2X$. As in the proof of (9.1) we have the
locus
$\Cal B \subset X\times \Bbb U$ which is the closure of the set
$$\{\,(x, z_1, z_2)\mid \pi (x) = \pi (z_i) \text{\, for some
$i$}\,\}.$$  The set $\Cal B$ is a closed analytic subset both of $X\times
\Bbb U$ and of $S\times \Bbb U$. The two sets $(S\times \Bbb U)-\Cal B$ and
$X\times \Bbb U$ cover $S\times \Bbb U$ and their intersection is $(X\times
\Bbb U)-\Cal B$. We shall show that there is an isomorphism of the
restriction of $\Cal V_{\Bbb U}\otimes \pi _1^*L$ to $(X\times \Bbb U)-\Cal
B$ with $\pi _1^*V_0$.

Let $\Pi \: (X\times \Bbb U)-\Cal B \to (D_0\times \Bbb U)-\Cal B'$ be the
projection, where $D_0$ is the base of $X$ and, as in the proof of (9.1),
$$\Cal B' = \{\,(t, z_1, z_2)\mid t = \pi (z_i) \text{\, for some
$i$}\,\}.$$ By construction, the restriction of $\Cal V_{\Bbb U}\otimes \pi
_1^*L$ to the reduction of every fiber of $\Pi$ is stable, and hence
isomorphic to the restriction of $\pi _1^*V_0$ to the fiber. Consider
$R^0\Pi _*Hom (\Cal V_{\Bbb U}\otimes \pi _1^*L, \pi _1^*V_0)$. By base
change and (1.5) for the case of a multiple fiber, this is a line bundle on
$(D_0\times
\Bbb U)-\Cal B'$. On the other hand, both $\Cal V_{\Bbb U}\otimes \pi _1^*L$
and
$\pi _1^*V_0$ extend to coherent sheaves on $X\times \Bbb U$, so that $R^0\Pi
_*Hom (\Cal V_{\Bbb U}\otimes \pi _1^*L, \pi _1^*V_0)$ also extends to a
coherent sheaf on $D_0 \times \Bbb U$. Arguing as in the proof of (9.1),
$R^0\Pi _*Hom (\Cal V_{\Bbb U}\otimes \pi _1^*L, \pi _1^*V_0)|(D_0\times
\Bbb U)-\Cal B'$ is a trivial line bundle and we may choose a section of $Hom
(\Cal V_{\Bbb U}\otimes \pi _1^*L, \pi _1^*V_0)$ which generates the fiber at
every point. This section then defines an isomorphism from $\Cal V_{\Bbb
U}\otimes \pi _1^*L$ to $\pi _1^*V_0)$ over
$(X\times \Bbb U)-\Cal B$. Thus we may define a coherent sheaf over $S\times
\Bbb U$, flat over $\Bbb U$, which by construction is a family of stable
torsion free sheaves on $S$. This sheaf defines a morphism from $\Bbb U$ to
$\overline{\frak M}_2(S)$ which is the desired extension. Doing this for a
neighborhood of every point of $\Bbb X_0'$ defines the extension over all of
$\Bbb X_0'$. \endproof

We return to the proof of (9.4). The proof will now follow from standard
algebraic topology. We have the moduli spaces
$\operatorname{Hilb}^2J^{e+1}(S)$ and $\overline{\frak M}_2(S)$ for $S$ and
corresponding moduli spaces
$\operatorname{Hilb}^2J^{e+1}(S_0)$ and $\overline{\frak M}_2(S_0)$ for
$S_0$. There are also the divisors
$\mu'(\Sigma)$ on $\operatorname{Hilb}^2J^{e+1}(S)$ and $\mu (\Sigma)$ on
$\overline{\frak M}_2(S)$, as well as the corresponding divisors $\mu
_0'(\Sigma_0)$ and $\mu _0(\Sigma _0)$ for $S_0$. Here $\Sigma \in H_2(S)$ and
$\Sigma _0 \in H_2(S_0)$. Finally we have the open sets $\Bbb X =\Bbb X_0$ and
$\Bbb X' \cong \Bbb X_0'$.

\claim{} If $\Sigma \cdot f = \Sigma _0\cdot f$, then $\mu'(\Sigma )|H^2(\Bbb
X) =
\mu '(\Sigma _0)|H^2(\Bbb X)$.
\endstatement
\proof We have $H^2(\Bbb X) \cong H^2(Y_i)$ by restriction. Here $Y_i$ is the
total transform of $\Sym ^2F_i$ in $\operatorname{Hilb}^2X$ and consists of
two components. One of these is the proper transform of $\Sym ^2F_i$ and the
other is the $\Pee ^1$-bundle over $F_i$ consisting of nonreduced length two
subschemes whose support lies in $F_i$. Clearly, if $D_\alpha + aE$ is a
divisor in
$H^2(\operatorname{Hilb}^2J^{e+1}(S))$, then the restriction of $D_\alpha +
aE$ to
$Y_i$ depends only on $\alpha \cdot f$ and $a$. By Lemma 9.2, $\mu'(\Sigma ) =
D_{\alpha _2} -\bigl((f\cdot \Sigma)/2\bigr)E$, where $\alpha _2 \cdot f =
\alpha _1\cdot f = 2(f\cdot \Sigma)$. A similar statement holds for $\mu
'_0(\Sigma _0)$. Thus $\mu'(\Sigma )|H^2(\Bbb X) = \mu '(\Sigma _0)|H^2(\Bbb
X)$. \endproof

To compare $\mu '(\Sigma)^4$ with $\mu (\Sigma )^4$, shrink $\Bbb X$ slightly
so that it is a manifold with boundary $\partial$. Thus $\Bbb X'$ can also be
shrunk slightly so that its boundary is $\partial$. Form the closed oriented
8-manifold
$\Bbb Y$ which is $\Bbb X$ glued to $-\Bbb X'$ along $\partial$. Doing the
same construction with $\Bbb X_0$ and $\Bbb X_0'$ gives an 8-manifold $\Bbb
Y_0$ diffeomorphic to $\Bbb Y$. Given $\Sigma \in H_2(S)$, the divisor $\mu
'(\Sigma)$ induces a class $\xi(i) \in H^2(\Bbb Y)$, and likewise  $\Sigma _0
\in H_2(S_0)$ induces a class $\xi_0(i) \in H^2(\Bbb Y_0)$. It follows from
the above claim that if $\Sigma \cdot f = \Sigma _0\cdot f$, then the classes
$\xi(i)$ and $\xi_0(i)$ agree under the natural identification of $\Bbb Y$
with $\Bbb Y_0$. Next we claim

\claim{} $\mu '(\Sigma)^4-\mu (\Sigma )^4 = \xi(1)^4+\xi (2)^4$.
\endstatement
\proof After passing to a multiple, we can assume that $\mu '(\Sigma)$ is
represented by a submanifold $\bold M$ of $\operatorname{Hilb}^2J^{e+1}(S)$.
Perturbing slightly, we can find four submanifolds $\bold M_1, \dots,
\bold M_4$ of $\operatorname{Hilb}^2J^{e+1}(S)$ whose signed intersection is
$\mu '(\Sigma)^4$. The closure of the image of $\bold M_i$ in $\Bbb X'$ is
the image of a blowup of $\bold M_i$. Restricting
$\bold M_i$ to
$\Bbb X$ and glue it to its closure in  image in $\Bbb X'$ gives a stratified
subspace
$M_i$ of
$\Bbb Y$ representing $\xi(i)$. Taking  small general deformations of $\bold
M_i$ gives stratified subspaces $M_1', M_2', M_3', M_4'$ meeting
transversally in finitely many points whose signed intersection number
calculates
$\xi(i) ^4$. After a small perturbation, we may glue $M_i'\cap \Bbb X_0$ to
$\bold M_i|
\operatorname{Hilb}^2J^{e+1}(S)-\Bbb X$ and use these to calculate $\mu
(\Sigma )^4$. Clearly the discrepancy between $\mu '(\Sigma )^4$ and $\mu
(\Sigma )^4$ is counted by  $\xi (1)^4+\xi (2)^4$. \endproof

Now $\xi (i)^4$ depends only on $(f\cdot \Sigma)$, $e$, and an analytic
neighborhood $X$ of the multiple fiber and is homogeneous of degree four in
$\Sigma$. Thus we can write $\xi (i)^4 = d(e, m_i)(f\cdot \Sigma )^4$ for some
rational number $d(e,m_i)$ depending only on the  analytic neighborhood $X$,
where $d(e,1)=0$. Using the previous claim, we see that
$$\mu '(\Sigma)^4-\mu (\Sigma )^4 = m_1^4m_2^4(d(e, m_1) + d(e, m_2))(\kappa
\cdot \Sigma )^4,$$ as claimed in Lemma 9.4.
\endproof

\section{Appendix: Elementary modifications.}

In this appendix, we consider the following problem (and its
generalizations): let
$X$ be a smooth projective scheme or compact complex manifold, let $T$ be
smooth and let $D$ be a smooth divisor on $T$. Suppose that $\Cal W$ is a
rank two vector bundle  over $X\times T$, and that $L$ is a line bundle on
$X$. Let $i\: X\times D
\to X\times  T$ be the inclusion, and suppose that there is a surjection
$\Cal W
\to i_*\pi _1^*L$ defining $\Cal V$ as an elementary modification:
$$0 \to \Cal V \to \Cal W \to i_*\pi _1^*L \to 0.$$ For $t\in T$, let $W_t =
\Cal W|X\times \{t\}$ and $V_t= \Cal V|X\times \{t\}$.  If $0$ is a reference
point of $D$, then there are two extensions
$$\align 0\to M \to &W_0 \to L \to 0; \\ 0 \to L \to &V_0 \to M \to 0.
\endalign$$ In particular the second exact sequence defines an extension
class $\xi \in H^1(M^{-1}\otimes L)$. We want a formula for $\xi$ and in
particular we want to  know some conditions which guarantee that $\xi \neq 0$.

\proposition{A.1} Let $\theta$ be the Kodaira-Spencer map for the family
$\Cal W$ of vector bundles over $X$, so that $\theta$ is a map from the
tangent space of
$T$ at $0$ to $H^1(Hom(W_0, W_0))$. Let $\partial/\partial t$ be a normal
vector to $D$ at $0$. Then the image of $\theta(\partial/\partial t)$ in
$H^1(M^{-1}\otimes L)$ under the natural map $H^1(Hom(W_0, W_0)) \to
H^1(Hom(M,L)) =H^1(M^{-1}\otimes L)$ is independent mod scalars of the choice
of
$\partial/\partial t$ and is the extension class corresponding to $V_0$.
\endstatement
\proof Since $W_0$ is given as an extension, there is an open cover $\{U_i\}$
of $X$ and transition functions for $W_0$ with respect to the cover $\{U_i\}$
of the form
$$\bar{A}_{ij} = \pmatrix \lambda _{ij} & *\\ 0 &\mu _{ij} \endpmatrix.$$
Letting $t$ be a local equation for $D$ near $0$, we can then choose
transition  functions for $\Cal W$  of the form $A_{ij} = \bar{A}_{ij} +
tB_{ij}$. With these choices of trivialization, a basis of local sections for
$\Cal V$ on
$U_i\times T$ is of the form $\{e_1, te_2\}$. Thus the transition functions
for
$\Cal V$ are given by $$
\pmatrix 1& 0\\0 &t^{-1} \endpmatrix \cdot \bigl(\bar{A}_{ij} + tB_{ij}\bigr)
\cdot  \pmatrix 1& 0\\0 &t \endpmatrix.$$ If $B_{ij} = \pmatrix
a&b\\c&d\endpmatrix$, then a calculation shows that the  transition functions
are equal to
$$\pmatrix \lambda _{ij} & t* \\ 0 & \mu _{ij} \endpmatrix + \pmatrix ta &
t^2b\\ c & td \endpmatrix= \pmatrix \lambda _{ij} & 0 \\ c & \mu _{ij}
\endpmatrix + tB'_{ij}.$$ Here $c$ is a matrix coefficient which naturally
corresponds to the image of $B_{ij}$ in $Hom (M,L)$. The proposition is just
the intrinsic formulation of this local calculation.
\endproof

\noindent {\bf Note.} The proof shows that, if the extension does split, then
we  can repeat the process, viewing $V_0$ again as extension of $L$ by $M$.
Either this procedure will eventually terminate, creating a nonsplit
extension at the generic point of $D$, or $\Cal W$ was globally an extension
in a neighborhood of
$D$. \medskip

Let us give another proof for (A.1) in intrinsic terms which, although less
explicit, will generalize. There are canonical identifications $H^1(Hom(W_0,
W_0)) = \Ext ^1(W_0, W_0)$ and $H^1(Hom(M,L))=\Ext ^1(M,L)$. For simplicity
assume that
$\dim T = 1$. Note that if we restrict the defining exact sequence for $\Cal
V$ to $X\times C$, where $C$ is a smooth curve in $T$ transverse to $D$, then
the sequence remains exact (since $\operatorname{Tor}^R _1(R/tR, R/sR) =0$ if
$t$ and $s$ are relatively prime elements of the regular local ring $R$).
Thus we can always restrict to the case where $\dim T=1$. Now
$\Spec \Cee [t]/(t^2)$ is a subscheme of $T$, and we can restrict $\Cal W$ to
$\Spec \Cee [t]/(t^2)$ to get a bundle $\Cal W_\varepsilon$. The bundle $\Cal
W_\varepsilon$ is naturally an extension $$0 \to W_0 \to \Cal W_\varepsilon
\to W_0
\to 0,$$ and the associated class in $\Ext ^1(W_0, W_0)$ is the
Kodaira-Spencer class. The  natural map $\Ext ^1(W_0, W_0) \to \Ext ^1(M,L)$
is defined on the level of extensions as follows: given an extension $\Cal
W_\varepsilon$ of $W_0$ by
$W_0$, let $\Cal E$ be the preimage of $M$ in $\Cal W_\varepsilon$, so that
there is an exact sequence $$0 \to W_0 \to \Cal E \to M \to 0.$$ Given the
map $W_0 \to
\Cal E$, the quotient $\Cal F = \bigl(\Cal E \oplus L\bigr)/W_0$, where $W_0$
maps diagonally into each summand, surjects onto $M$ by taking the
composition of the projection to $\Cal E$ with the given map $\Cal E \to M$.
The kernel is naturally
$L$. Thus $\Cal F$ is an extension of $M$ by $L$, and it is easy to see that
$\Cal F$ corresponds to the image of the extension class for $\Cal
W_\varepsilon$ under the natural map. Finally note that, since $W_0 \to L$ is
surjective, there is a natural identification of $\Cal F= \bigl(\Cal E \oplus
L\bigr)/W_0$ with $\Cal E/M$ where we take the image of $M$ under the map
$M\to W_0 \to \Cal E$.

On the other hand, restricting the defining exact sequence for $\Cal V$ to
$\Spec \Cee [t]/(t^2)$ gives a new exact sequence
$$0 \to L \to \Cal V_\varepsilon \to \Cal W_\varepsilon \to L \to 0.$$ If we
set $\Cal E$ to be the image of $\Cal V_\varepsilon$ in $\Cal
W_\varepsilon$,  then it is clear that $\Cal E$ is the inverse image of
$M\subset W_0$ under the natural map. Now there is an isomorphism $\Cal
V_\varepsilon/L \cong
\Cal E$, and it is easy to see that this isomorphism identifies $V_0$ with
$\Cal E/M$ under the natural maps, compatibly with the extensions. Thus the
extension of
$M$ by $L$ defined by $V_0$ has an extension class equal to the image of the
extension class of $\Cal W_\varepsilon$ in $\Ext ^1(M,L)$ under the natural
map.

With this said, here is the promised generalization of (A.1):

\proposition{A.2} With notation at the beginning of this section, let $\Cal
W$ be  a rank two reflexive sheaf over $X\times T$, flat over $T$, let $D$ be
a reduced divisor on $T$, not necessarily smooth and let $i\: D \to T$ be the
inclusion. Suppose that $L$ is a line bundle on $X$ and that $\Cal Z$ is a
codimension two subscheme of $X\times D$, flat over $D$. Suppose further that
$\Cal W \to i_*\pi _1^*L\otimes I_{\Cal Z}$ is a surjection, and let $\Cal V$
be its kernel:  $$0 \to \Cal V \to \Cal W \to i_*\pi _1^*L\otimes I_{\Cal Z}
\to 0.$$ Then:
\roster
\item"{(i)}" $\Cal V$ is reflexive and flat over $T$.
\item"{(ii)}" For each $t\in D$, there are exact sequences
$$\align 0 \to M\otimes I_{Z'} \to &W_t \to L\otimes I_Z \to 0;\\ 0\to
L\otimes I_Z \to &V_t \to M\otimes I_{Z'} \to 0,
\endalign$$ where $Z$ is the subscheme of $X$ defined by $\Cal Z$ for the
slice $X\times
\{t\}$ and ${Z'}$ is a subscheme of $X$ of codimension at least two.
\item"{(iii)}" If $D$ is smooth, then the extension class corresponding to
$V_t$ in $\Ext ^1(M\otimes I_W,  L\otimes  I_Z)$ is defined by the image of
the normal vector to $D$ at $t$ under the composition of the Kodaira-Spencer
map from the tangent space of $T$ at $t$ to
$\Ext ^1(W_t, W_t)$, followed by the natural map $\Ext ^1(W_t, W_t) \to \Ext
^1(M\otimes I_{Z'},  L\otimes I_Z)$. \endroster
\endstatement
\proof First note that $\Cal V$ is a subsheaf of $\Cal W$ and is therefore
torsion  free. Given an open set $U$ of $X\times T$ and a closed subscheme
$Y$ of $U$ of codimension at least two, let $s$ be a section of $\Cal V$
defined on $U-Y$. Then
$s$ extends to a section $\tilde s$ of $\Cal W$ over $U$ since $\Cal W$ is
reflexive. Moreover the image of $\tilde s$ in $H^0(U\cap D; L\otimes I_{\Cal
Z})$ vanishes in codimension one and thus everywhere. Thus $\tilde s$ defines
a section of $\Cal V$ over $U$, and so $\Cal V$ is reflexive. That it is flat
over
$T$ follows from the next lemma:

\lemma{A.3} Let $R$ be a ring and $t$ an element of $R$ which is not a zero
divisor. Let $I$ be an $R/tR$-module which is flat over $R/tR$. For an
$R$-module
$N$, let $N_t$ be the kernel of multiplication by $t$ on $N$.
\roster
\item"{(i)}" For all $R$-modules $N$, $\operatorname{Tor}_1^R(I,N) = I\otimes
_{R/tR}N_t$, and $\operatorname{Tor}_i^R(I,N) = 0$ for all $i>1$.
\item"{(ii)}" Suppose that there is an exact sequence of $R$-modules
$$0 \to M_2 \to M_1 \to I \to 0,$$ where $M_1$ is flat over $R$. Then $M_2$
is flat over $R$ as well.
\endroster
\endstatement
\proof The statement (i) is easy if $I=R/tR$, by taking the free resolution
$$0\to R @>{\times t}>>R \to R/tR \to 0.$$ Thus it holds more generally if
$I$ is a free $R/tR$-module. In general, start with a free resolution
$F^\bullet$ of $I$. By standard homological algebra (see e.g\. [EGA III
6.3.2]) there is a spectral sequence with $E_1$ term
$\operatorname{Tor}_p^R(F^q, N)$ which converges to
$\operatorname{Tor}_{p+q}^R(I,N)$. The only nonzero rows correspond to
$p=0,1$ and the row for $p=1$ is the complex $F^\bullet \otimes _{R/tR}N_t$.
Since $I$ is flat, this complex is exact except in dimension zero and is a
resolution of
$I\otimes _{R/tR}N_t$. Since $F^q\otimes _RN = F^q\otimes _{R/tR}(N\otimes
_RR/tR)$, the flatness of $I$ over $R/tR$ implies  that the row for $p=0$ is
exact except in dimension zero. Thus $\operatorname{Tor}_1^R(I,N) = I\otimes
_{R/tR}N_t$.

The second statement now follows since, for every $R$-module $N$, the long
exact sequence for $\operatorname{Tor}$ defines an isomorphism, for all
$i\geq 1$, from
$\operatorname{Tor}_i^R(M_2, N)$ to $\operatorname{Tor}_{i+1}^R(I, N)=0$.
\endproof

Returning to (A.2), let us prove (ii). There is a surjection $W_t \to L\otimes
I_Z$ and the  kernel of this surjection is a rank one torsion free sheaf on
$X$, which is thus of the form  $M\otimes I_{Z'}$ for some subscheme ${Z'}$
of $X$ of codimension at least two. Now there is an exact sequence
$$Tor_1^{\scrO_{X\times T}}(i_*\pi _1^*L\otimes I_{\Cal Z}, \scrO _{X\times
\{t\}})
\to V_t \to W_t \to L\otimes I_Z \to 0.$$ In the $Tor_1$ term, the first
sheaf is an $\scrO_{X\times D}$-module, flat over $D$, and the second is an
$\scrO_D$-module. Using (i) of (A.3) identifies
$Tor_1^{\scrO_{X\times T}}(i_*\pi _1^*L\otimes I_{\Cal Z},\scrO _{X\times
\{t\}})$ with $L\otimes I_Z$.  Thus we obtain the exact sequence for $V_t$.

Finally, the identification of the extension class in (iii) is formally
identical  to the second proof of (A.1) given above and will not be repeated.
\endproof

Next we shall give some criteria for when the extension is nonsplit. The
simplest case is when the Kodaira-Spencer map is an isomorphism at $0$. In
this case we can check whether or not the extension is split by looking at
the map
$\Ext ^1(W_t, W_t) \to \Ext^1(M\otimes I_{Z'},  L\otimes I_Z)$. Thus the
problem  is essentially cohomological. A similar application concerns the
case where $T$ is the blowup of a universal family along the locus where the
sheaves  are extensions of  $L\otimes I_Z$ by $M\otimes I_{Z'}$. In our
applications, however, we shall need a more general situation and will have
to analyze some first order information about the family $\Cal W$. For
simplicity we shall assume that $\dim T=1$, with $t$ a coordinate. It is an
easy consequence of (A.3)(i) that the general case can be reduced to this
special case by taking a curve in $T$ transverse to $D$.

\proposition{A.4} In the notation of \rom{(A.3)}, let $\Cal W_\varepsilon$ be
the restriction of $\Cal W$ to $\Spec \Cee [t]/t^2$. Suppose that
\roster
\item"{(i)}" $\Hom (M\otimes I_{Z'}, L\otimes I_Z) =0$.
\item"{(ii)}" The map from
$\Ext ^1(M\otimes I_{Z'}, \Cal W_\varepsilon)/t\Ext ^1(M\otimes I_{Z'},  \Cal
W_\varepsilon)$ to $\Ext ^1(M\otimes I_{Z'}, \Cal W_\varepsilon)$ induced by
multiplication by $t$ has a one-dimensional kernel.
\endroster Then we may identify the kernel with a line in $\Ext ^1(M\otimes
I_{Z'}, W_0)$, and if the image of this line in $\Ext ^1(M\otimes I_{Z'},
L\otimes I_Z)$ is
$\Cee \cdot \xi$ then the corresponding extension class is $\xi$.
\endstatement
\proof From the first assumption $\dim \Hom (W_0, W_0) =1$. Thus if $\theta$
is the Kodaira-Spencer class, there is an exact sequence
$$0\to \Ext ^1(M\otimes I_{Z'}, W_0)/\Cee \cdot\theta \to \Ext ^1(M\otimes
I_{Z'},  \Cal W_\varepsilon)\to \Ext ^1(M\otimes I_{Z'}, W_0).$$
Multiplication by $t$ induces the natural map
$$\im(\Ext ^1(M\otimes I_{Z'},
\Cal W_\varepsilon))\subseteq \Ext ^1(M\otimes I_{Z'}, W_0) \to
\Ext ^1(M\otimes I_{Z'}, W_0)/\Cee \cdot\theta.$$ If this map has a kernel
then clearly $\theta \in \im(\Ext ^1(M\otimes I_{Z'},
\Cal W_\varepsilon))$ and the kernel is $\Cee \cdot \theta$. The image of the
kernel in $\Ext ^1(M\otimes I_{Z'}, L\otimes I_Z)$ is then just the image of
the Kodaira-Spencer class.
\endproof

Here is the typical way we will apply the above: suppose that $\Cal W$ is
locally free and that $Z' = \emptyset$. Then $\Ext ^1(M\otimes I_{Z'}, \Cal
W_\varepsilon) = R^1\pi _2{}_*(\Cal W_\varepsilon \otimes \pi _1^*M^{-1})$.
Suppose in addition that $\Cal W$ is globally an extension:
$$0 \to \pi _1^*\Cal L_1 \to \Cal W \to \pi _1^*\Cal L_2\otimes I_{\Cal Y} \to
0,$$ where $\Cal Y\subset X\times T$ is flat over $T$. Thus there is a map
$$R^0\pi _2{}_*\bigl(\pi _1^*\Cal L_2\otimes I_{\Cal Y}\otimes \pi
_1^*M^{-1}\bigr)
\to  R^1\pi _2{}_*\pi _1^*(\Cal L_1\otimes M^{-1})$$ whose cokernel sits
inside $R^1\pi _2{}_*(\Cal W\otimes \pi _1^*M^{-1})$. A similar statement is
true when we restrict to $\Spec \Cee [t]/(t^2)$. Now suppose that $\dim
H^0(X;\Cal L_2\otimes M^{-1}\otimes I_{Y_t})$ is independent of $t$. Then the
sheaves $R^0\pi _2{}_*\bigl(\pi _1^*\Cal L_2\otimes I_{\Cal Y}\otimes \pi
_1^*M^{-1}\bigr)$ and $R^1\pi _2{}_*\pi _1^*(\Cal L_1\otimes M^{-1})$ are
locally free and compatible with base change, by [EGA III, 7.8.3, 7.8.4,
7.7.5] so if we know that the map between them has a determinant which
vanishes simply along $D$ then the same will be true for the restrictions to
$\Spec \Cee [t]/(t^2)$.  The image in $R^1\pi _2{}_*(\Cal W\otimes \pi
_1^*M^{-1})$ is the direct image of a line bundle $\Cal K$ on $D$.
Furthermore suppose that $\dim H^1(X;\Cal L_2\otimes M^{-1} \otimes I_{Y_t})$
is independent of $t$. Then  $R^0\pi _2{}_*\bigl(\pi _1^*\Cal L_2\otimes
I_{\Cal Y}\otimes \pi _1^*M^{-1}\bigr)$ is locally free and compatible with
base change. If it is nonzero suppose further that
$R^2\pi _2{}_*\pi _1^*(\Cal L_1\otimes M^{-1}) =0$. Thus the torsion part of
$R^1\pi _2{}_*(\Cal W\otimes \pi _1^*M^{-1})$ is just $\Cal K$ and the
restriction of $\Cal K$ to $\Spec \Cee [t]/(t^2)$ gives the kernel of
multiplication by $t$ as in (ii). We can then take the map from the torsion
part of  $R^1\pi _2{}_*(\Cal W\otimes \pi _1^*M^{-1})_0$, namely the image of
$H^1(\Cal L_1\otimes M^{-1})$, to $H^1(L\otimes I_Z\otimes M^{-1}) = \Ext
^1(M, L\otimes I_Z)$ and this image gives the extension class.

\medskip We will also need to consider a slightly different situation.
Suppose that $\Cal W$ is a rank two vector bundle on $X\times T$, $E$ is a
smooth divisor on $X$ and
$L$ is a line bundle on $E\times T$. Let $j\: E\times T \to X\times T$ be the
inclusion and let $\Phi\:\Cal W \to j_*L$ be a morphism. We may think of
$\Phi$ as a family of morphisms parametrized by $T$. In local coordinates
$\Phi$ is given by two functions $f$, $g$ on $E\times T$, whose vanishing
defines a subscheme $Y$ of $E\times T$. Away from the projection $\pi _2(Y)$
of $Y$ to $T$, $\Phi$ defines a family of elementary modifications which
degenerates over $\pi _2(Y)$ at the points of $Y$.

\proposition{A.5} Let $\Phi\:\Cal W \to j_*L$ be a morphism as above and
suppose that the cokernel of $\Phi$ is supported on  a nonempty codimension
two subset $Y$ of $E\times T$, necessarily a local complete intersection.
Suppose further that, for each $t\in T$, the codimension of $Y\cap (X\times
\{t\})$ in $X\times \{t\}$ is at least two if $Y\cap (E\times \{t\})\neq
\emptyset$. Let $\Cal V$ be the kernel of $\Phi$. Then $\Cal V$ is a
reflexive sheaf, flat over $T$, and its restriction to each slice $X\times
\{t\}$ is a torsion free sheaf on $X$.
\endstatement
\proof The proof of (A.2)(i) shows that $\Cal V$ is reflexive. As for the
rest, the problem is local around a point of $Y$. Let $R$ be the local ring
of $X\times T$ at a point
$(x,t)$, $R'$ the local ring of $T$ at $t$,  and $S$  the local ring of
$X\times
\{t\}$ at $(x,t)$. Let $u$ be the local equation for $E$ in $X\times T$. Then
locally $\Phi$ corresponds to a map $R\oplus R \to R/uR$, necessarily given by
elements $\bar f, \bar g \in R/uR$. Lift $\bar f$ and $\bar g$ to elements $f,
g\in R$. Then $(u,f,g)R$ is the ideal of $Y$ in $R$, and $Y$ has codimension
three in $X\times T$. Thus $u, f,g$ is a regular sequence, any two of the
three are relatively prime, and necessarily $\dim R\geq 3$.

The kernel $M$ of the map $R\oplus R \to R/uR$ given by $(a,b)\mapsto a\bar f
+ b\bar g$ is clearly generated by $(-g,f)$, $(u,0$, and $(0, u)$. These three
elements define a surjection $R\oplus R\oplus R \to M$. The kernel of this
surjection is easily calculated to be $R\cdot(u, g, -f)$. Thus there is an
exact sequence $$0 \to R \to  R\oplus R\oplus R \to M \to 0.$$ This sequence
restricts to define
$$S\to S\oplus S\oplus S \to M\otimes _RS \to 0.$$ Here the image of $S$ in
$S\oplus S\oplus S$ is equal to $S\cdot (u, g, -f)$, where we denote the
images of $u,f,g$ in $S$ by the same letter. By hypothesis, not all of $u$,
$f$, $g$ vanish on $X\times \{t\}$ and so this map is injective. By the local
criterion of flatness  $M$ is flat over $R'$. Finally we must show that
$M\otimes _RS$ is a torsion free $S$-module. By hypothesis $u, g, -f$
generate the ideal of a subscheme of $\Spec S$ of codimension at least two
and thus $u$ does not divide both $f$ and $g$ in $S$.  Given $h\in S$ with
$h\neq 0$, suppose that $hm = 0$ for some $m\in M\otimes _RS$. Then there is
$(a,b,c)\in S\oplus S\oplus S$ such that $h(a,b,c) = \alpha (u, g, -f)$. We
claim that $u|a$. To see this, let $n$ be the largest integer such  that
$u^n|h$. Then $u^n|hb = \alpha g$ and likewise $u^n|\alpha f$. Since at least
one of $f$, $g$ is prime to $u$, $u^n|\alpha$. But then $u^{n+1}|\alpha u =
ha$, so that $u|a$. If $a=ua'$, then $\alpha = ha'$ and so $hb= ha'g$ and
$b=a'g$. Likewise $c=a'(-f)$. Thus $(a,b,c) = a'(u, g, -f)$ and its image in
$M\otimes _RS$ is zero. It follows that $M\otimes _RS$ is torsion free.
\endproof

Let us finally remark that we can calculate the class $p_1(\ad \Cal V)$, in
the above notation, by applying the formula (0.1), since $\Phi$ is surjective
in codimension two.

\enddocument

\bye